\documentclass[10pt]{amsart} 

\usepackage{microtype}
\usepackage[dvipsnames]{xcolor} 

\usepackage{amsmath, amssymb, mathrsfs} 
\usepackage{mathtools}
\usepackage[mathscr]{euscript} 
 
\newlength{\mylength}
\usepackage{enumitem}
\setlist{listparindent=\parindent, itemsep=0cm, parsep=\mylength, topsep=0cm}

\usepackage[breaklinks=true]{hyperref} 
\usepackage{comment} 
\usepackage{url}
\usepackage{tikz-cd}
\usepackage{amsthm}

\makeatletter
\renewenvironment{proof}[1][\proofname]{\par
	\pushQED{\qed}%
	\normalfont \topsep6\p@\@plus6\p@\relax
	\noindent\emph{#1.} 
	\ignorespaces
}{%
\popQED\endtrivlist\@endpefalse
}
\makeatother

\newtheoremstyle{mythm}% name of the style to be used
{\mylength}% measure of space to leave above the theorem. E.g.: 3pt
{0pt}% measure of space to leave below the theorem. E.g.: 3pt
{\itshape}% name of font to use in the body of the theorem
{0pt}% measure of space to indent
{}% name of head font
{. }% punctuation between head and body
{0em}% space after theorem head; " " = normal interword space
{\thmnumber{#2. }\bfseries{\thmname{#1}\thmnote{ (#3)}}}

\newtheoremstyle{mydef}% name of the style to be used
{\mylength}% measure of space to leave above the theorem. E.g.: 3pt
{0pt}% measure of space to leave below the theorem. E.g.: 3pt
{}% name of font to use in the body of the theorem
{0pt}% measure of space to indent
{}% name of head font
{. }% punctuation between head and body
{0em}% space after theorem head; " " = normal interword space
{\thmnumber{#2. }\bfseries{\thmname{#1}\thmnote{ (#3)}}}

\newtheoremstyle{myrmk}% name of the style to be used
{\mylength}% measure of space to leave above the theorem. E.g.: 3pt
{0pt}% measure of space to leave below the theorem. E.g.: 3pt
{}% name of font to use in the body of the theorem
{0pt}% measure of space to indent
{\itshape}% name of head font
{.\ }% punctuation between head and body
{ }% space after theorem head; " " = normal interword space
{\thmname{#1}\thmnumber{ #2}\thmnote{ (#3)}}

\theoremstyle{mythm} 
\newtheorem{thm}[subsubsection]{Theorem}
\newtheorem{lem}[subsubsection]{Lemma} 
\newtheorem{cor}[subsubsection]{Corollary}

\newtheorem{prop}[subsubsection]{Proposition}
\newtheorem*{thm*}{Theorem} 
\newtheorem*{lem*}{Lemma} 
\newtheorem*{cor*}{Corollary} 
\newtheorem*{claim*}{Claim} 
\newtheorem*{prop*}{Proposition} 

\theoremstyle{mydef}

\theoremstyle{myrmk}
\newtheorem*{defn*}{Definition}

\newtheorem*{rmk*}{Remark} 
\newtheorem*{rmks*}{Remarks} 
\newtheorem*{ex*}{Example} 

\newcommand{\on}{\operatorname}
\newcommand{\rnc}{\renewcommand} 

\rnc{\setminus}{\smallsetminus} 

\newcommand{\wt}{\widetilde}
\newcommand{\wh}{\widehat} 
\newcommand{\ol}{\overline}

\newcommand{\BZ}{\mathbb{Z}}

\newcommand{\BR}{\mathbb{R}}
\newcommand{\BC}{\mathbb{C}}
\newcommand{\BA}{\mathbb{A}}
\newcommand{\ul}{\underline}

\newcommand{\mb}{\mathbf}
\newcommand{\mf}{\mathfrak}
\newcommand{\mc}{\mathscr}

\newcommand{\hra}{\hookrightarrow}
\newcommand{\sra}{\twoheadrightarrow}

\newcommand{\Spec}{\on{Spec}}
\rnc{\Im}{\on{Im}}

\newcommand{\Hom}{\on{Hom}}
\newcommand{\Fun}{\on{Fun}}

\newcommand{\oh}{\mc{O}}
\newcommand{\D}{\mc{D}}
\DeclareMathOperator*\colim{colim}

\newcommand{\spl}{\mathcal{S}\kern -0.5pt p}

\newenvironment{cd}{\begin{equation*}\begin{tikzcd}}{\end{tikzcd}\end{equation*}\ignorespacesafterend}
\newcommand{\e}[1]{\begin{align*} #1 \end{align*}}

\usepackage[margin=1.5in]{geometry}

\makeatletter
\def\blfootnote{\gdef\@thefnmark{}\@footnotetext}
\makeatother

\title[A colimit presentation of $\D(G(K))$ via the Bott--Samelson hypercover]{A colimit presentation of $\D(G(K))$ \\ via the Bott--Samelson hypercover} 
\author{James Tao} 
\address{Massachusetts Institute of Technology, Cambridge, MA 02139, USA}
\email{jamestao@mit.edu}
\author{Roman Travkin}
\address{Skolkovo Institute of Science and Technology, Moscow, Russia}
\email{roman.travkin2012@gmail.com}
\date{March 25, 2021} 

\definecolor{myblue}{rgb}{0,0.1,0.4}
\newenvironment{myproof}{\color{myblue}\begin{proof}}{\end{proof}}

\usepackage{bm}

\newcommand{\alg}{\mathsf{Alg}}
\newcommand{\cat}{\mathsf{DGCat}}
\rnc{\mod}{\text{-}\mathsf{mod}}
\rnc{\bmod}{\mathsf{BMod}}

\newcommand{\subsubsectiona}{\subsubsection{}\hspace{-0.7em}}

\newcommand{\mr}{\mathrm}

\newcommand{\simp}{\mathsf{Simp}}
\newcommand{\ner}{\on{N}}

\newcommand{\xra}{\xrightarrow} 

\newcommand{\ms}{\mathsf}

\newcommand{\sd}{\on{Sd}}
\newcommand{\exx}{\on{Ex}}
\newcommand{\sset}{\ms{SSet}}
\newcommand{\set}{\ms{Set}}
\newcommand{\gal}{\ms{Gal}}
\newcommand{\LG}{\mc{L}G}

\newcommand{\pt}{\mr{pt}}
\newcommand{\op}{\mr{op}}
\newcommand{\Id}{\operatorname{Id}}
\newcommand{\pr}{\operatorname{pr}}

\newcommand{\fcsimpset}{\ms{SSet}_{\mr{fc}}}
\newcommand{\fsimpset}{\ms{SSet}_{\mr{f}}}

\begin{document}
	
	\begin{abstract}
		Let $G$ be a semisimple, simply connected algebraic group over an algebraically closed field of characteristic zero. We prove that the $\infty$-category of $\D$-modules on the loop group of $G$ is equivalent to the monoidal colimit of the $\infty$-categories of $\D$-modules on the standard parahoric subgroups. This also follows from~\cite{hecke}, but the present paper gives a simpler proof. The idea is to develop a combinatorial model for the path space of a simplicial complex, in which `paths' are sequences of adjacent simplices, and to use a generalized version of hyperdescent for $\D$-modules. We also give two more applications of this hyperdescent theorem: triviality of $\D$-modules on the `schematic Bruhat--Tits building,' which was first established by Varshavsky using a different method, and triviality of $\D$-modules on the `simplicial affine Springer resolution.' 
	\end{abstract}
	
	\maketitle
	
	\newpage 
	
	\tableofcontents
	
	\newpage
	
	\section{Introduction} 
	
	Our main theorem realizes the $\infty$-category of $\D$-modules on the loop group as a monoidal colimit of the $\infty$-categories of $\D$-modules on the standard parahoric subgroups. In~\ref{intro1}, we review some classical analogues of this statement. In~\ref{intro2}, we state the theorem and remark that it can be proved using the method of a different paper~\cite{hecke}. The purpose of the current paper is to give a much simpler proof which only uses well-understood ideas from topology and descent theory. These ideas are discussed in~\ref{intro3} and \ref{intro4}. In~\ref{intro5}, we explain two other applications of these ideas which do not depend on the main theorem. 
	
	\subsection{Some colimit presentations of groups} \label{intro1} 
	
	Let us review some known results which state that a group of Lie-theoretic nature is the colimit of its standard parabolic subgroups. 
	
	First, there is the case of Coxeter groups. A Dynkin diagram with vertex set $I$ determines a Coxeter group $W_I$ which is generated by a set of simple reflections indexed by $I$. For each $J \subseteq I$, we have the parabolic subgroup $W_J \subseteq W_I$ generated by the simple reflections in $J$. Since the Coxeter relations involve only two generators at a time, we have 
	\begin{equation}\tag{$1$} \label{e1} 
		W_I \simeq \colim_{\substack{J \subseteq I \\ |J| \le 2}} W_J
	\end{equation}
	where the colimit is taken in groups. 
	
	A subset $J \subseteq I$ is called \emph{finite type} if $W_J$ is finite. Unlike the previous statement, the following colimit presentation remains valid when the colimit is taken in the $\infty$-category of homotopy-coherent groups in $\ms{Spaces}$: 
	\begin{equation}\tag{$1'$} \label{e2} 
		W_I \simeq \colim_{\substack{J \subseteq I \\ J \text{ finite type}}} W_J 
	\end{equation}
	To prove this, first recall that the Coxeter complex is contractible:
	\[
		\pt \simeq \colim_{\substack{J \subseteq I \\ J \text{ finite type}}} W_I / W_J 
	\]
	(See the remark following Theorem 2.16 in~\cite{mitchell}.) Taking homotopy quotients by the $W_I$-action yields a colimit presentation of the classifying space $\pt / W_I$, and passing to loop spaces yields~(\ref{e2}). 
	
	Next, there is the case of Ka\v{c}--Moody groups. Let $\wh{G}$ be a Ka\v{c}--Moody group over $\BC$ associated to a Dynkin diagram $I$.  For $J \subseteq I$, let $\mb{P}_J \subseteq \wh{G}$ be the corresponding standard parabolic subgroup. The definition of Ka\v{c}--Moody group given in~\cite[6.1.16]{kumar} easily implies the following: 
	\begin{equation} \tag{$2$} \label{e3} 
		\wh{G} \simeq \colim_{\substack{J \subseteq I \\ |J| \le 2}} \mb{P}_J 
	\end{equation}
	where the colimit takes place in groups. Indeed, that definition constructs $\wh{G}$ as a colimit involving the groups $N_I$ (the normalizer of the torus) and $\mb{P}_{\{i\}}$ for $i \in I$. On the other hand,~\cite[Prop.\ 5.1.7]{kumar} expresses $\mb{P}_{\{i, j\}}$ as an analogous colimit involving the groups $N_{\{i, j\}} := N_I \cap \mb{P}_{\{i, j\}}$ and $\mb{P}_{\{i\}}, \mb{P}_{\{j\}}$. Now~(\ref{e3}) follows by comparing these colimit presentations and noting that 
	$
		N_I = \colim_{|J| \le 2} N_J, 
	$
	which follows from~(\ref{e1}). 
	
	As before, there is an analogous statement about topological groups: Theorem~4.2.3 in~\cite{nitu} says that  
	\begin{equation} \tag{$2'$} \label{e4} 
		\wh{G} \simeq \colim_{\substack{J \subseteq I \\ J \text{ finite type}}} \mb{P}_J, 
	\end{equation}
	where $\wh{G}$ and $\mb{P}_J$ are given the analytic topology, and the colimit takes place in the $\infty$-category of homotopy-coherent groups in $\ms{Spaces}$. The proof is entirely similar to~(\ref{e2}): Theorem 2.16 of~\cite{mitchell} establishes the contractibility of the `topological Bruhat--Tits building,' namely 
	\[
		\pt \simeq \colim_{\substack{J \subseteq I \\ J \text{ finite type}}} \wh{G} / \mb{P}_J. 
	\]
	Taking homotopy quotients by $\wh{G}$ and passing to loop spaces yields~(\ref{e4}). 
	
	\subsection{The main theorem} \label{intro2} 
	
	We now state an analogue of~\ref{intro1}(\ref{e4}) for $\D$-module $\infty$-categories on Ka\v{c}--Moody groups. In fact, to keep the statement nice, we restrict to a very special case. 
	
	Let $G$ be a semisimple simply-connected algebraic group over an algebraically closed field $k$ of characteristic zero. Let $\LG$ be its loop group, with affine Dynkin diagram $I$. Each subset $J \subset I$ indexes a standard parahoric subgroup $\mb{P}_J \subset \LG$. 
	
	Let $\D(\LG)$ denote the $\infty$-category of $\D$-modules on $\LG$, as in~\cite[3.4]{ber}. The group multiplication on $\LG$ defines a `convolution' monoidal structure on $\D(\LG)$ via $*$-pushforward. As explained in~\cite{ber}, the theory of module $\infty$-categories for $\D(\LG)$ provides an interesting categorification of representation theory for $\LG$. 
	
	\begin{thm*}[\ref{colim-thm}]
		We have 
		\[
			\D(\LG) \simeq \colim_{\substack{J \subset I \\ J \textnormal{ finite type}}} \D(\mb{P}_J), 
		\]
		where the colimit takes place in the $\infty$-category of monoidal stable $\infty$-categories (or monoidal DG-categories). 
	\end{thm*}
	
	This result is useful in practice. It says that, to construct a $\D(\LG)$-action on an $\infty$-category, it suffices to construct compatible $\D(\mb{P}_J)$-actions. This is often easier because each $\mb{P}_J$ is a scheme, while $\LG$ is only an ind-scheme. 
	
	\subsubsection{Remarks} \label{main-remarks} 
	
	\begin{enumerate}[label=(\arabic*)]
		\item The requirement that $G$ is semisimple and simply-connected is not essential. If one replaces the index diagram in Theorem~\ref{colim-thm} by Varshavsky's category of parahorics (defined in~\cite{varshavsky}), then the theorem holds for all reductive groups $G$. 
		
		The complication arises because, when $G$ is reductive, $(\LG, \mb{I})$ is not a $(B, N)$-pair, but rather a generalized $(B, N)$-pair. (Here $\mb{I} = \mb{P}_{\emptyset}$ is the Iwahori subgroup.) This causes some statements in `building theory' to break; for example, $\mb{I}$ is now merely one of the connected components of its normalizer. The category of parahorics puts in additional morphisms (between pairs of objects $J, J' \subseteq I$) which are used to keep track of the extra connected components that arise throughout. 
		
		In~\ref{intro5}, we will discuss a result which was first proved by Varshavsky in~\cite{varshavsky}.  
		\item The theorem should remain true when $\LG$ is replaced by an arbitrary Ka\v{c}--Moody group, as soon as one defines the $\infty$-category of $\D$-modules and the analogous `category of parahorics' which manages extra connected components. This is because our proof only uses buildings and flag varieties, and these are well-developed for Ka\v{c}--Moody groups. (See \cite[7.1.3]{kumar} which defines Bott--Samelson varieties in this setting.) We restricted to the affine case because that is our primary interest. 
		\item In another paper~\cite{hecke}, we proved the colimit presentation for the affine Hecke category using a more complicated method based on a cell-by-cell analysis. As remarked in~\cite[1.1]{hecke}, if one decomposes $\D(\LG)$ into cells based on the Bruhat stratification of $\LG$, then the same method also proves Theorem~\ref{colim-thm}. 
		
		Why did we write the present paper? Our goal is to provide a second proof of Theorem~\ref{colim-thm} which is simpler and more morally correct. This second proof is shorter and involves no tricks. Instead, it uses two very classical ideas: hyperdescent and the contractibility of the Bruhat--Tits building. The main takeaway of this paper is that Theorem~\ref{colim-thm} is a consequence of these two ideas.  
	\end{enumerate}
	
	\subsubsection{Outline of the proof} \label{outline} 
	We first construct a functor 
	\[
		\mc{X} : \mc{C} \to \ms{Sch}_{/\LG}
	\]
	which is the `Bott--Samelson hypercover' mentioned in the title. It will be the diagram consisting of all Bott--Samelson schemes with all group-theoretic maps between them. 
	
	The objects of $\mc{C}$ are finite sequences $(J_1, \ldots, J_n)$ of finite type subsets of $I$, and the functor sends each sequence to the corresponding Bott--Samelson scheme: 
	\[
		\mc{X}((J_1, \ldots, J_n)) := \mb{P}_{J_1} \overset{\mb{I}}{\times} \cdots \overset{\mb{I}}{\times} \mb{P}_{J_n}. 
	\]
	The arrows in $\mc{C}$ are generated by the following: 
	\begin{cd}[row sep = 4mm, column sep = 0.7in] 
		(\ldots, J, J, \ldots) \ar[d, "\text{merge}"] & (\ldots, J, \ldots) \ar[d, "\text{increase}"] & (\ldots, J_1, J_2, \ldots) \ar[d, "\text{insert}"] \\
		(\ldots, J, \ldots) & (\ldots, J', \ldots) & (\ldots, J_1, \emptyset, J_2 \ldots)
	\end{cd}
	where $J \subset J'$. The functor sends these to the following maps of Bott--Samelson schemes: 
	\begin{cd}[row sep = 4mm] 
		\cdots \overset{\mb{I}}{\times} \mb{P}_J \overset{\mb{I}}{\times} \mb{P}_J \overset{\mb{I}}{\times} \cdots \ar[d, "\text{multiply}"] & \cdots \overset{\mb{I}}{\times} \mb{P}_J \overset{\mb{I}}{\times} \cdots \ar[d, "\text{embed}"] & \cdots \overset{\mb{I}}{\times} \mb{P}_{J_1} \overset{\mb{I}}{\times} \mb{P}_{J_2} \overset{\mb{I}}{\times} \cdots \ar{d}{\text{insert}}[swap, rotate=90, anchor=south]{\sim} \\
		\cdots \overset{\mb{I}}{\times} \mb{P}_J \overset{\mb{I}}{\times} \cdots & \cdots \overset{\mb{I}}{\times} \mb{P}_{J'} \overset{\mb{I}}{\times} \cdots & \cdots \overset{\mb{I}}{\times} \mb{P}_{J_1} \overset{\mb{I}}{\times} \mb{I} \overset{\mb{I}}{\times} \mb{P}_{J_2} \overset{\mb{I}}{\times} \cdots
	\end{cd}
	
	To begin computing the colimit in Theorem~\ref{colim-thm}, we use a `bar construction' to rewrite it as a colimit over a larger diagram which takes place in (non-monoidal) $\infty$-categories: 
	\[
		\on{oblv}\left( \colim_{\substack{J \subset I \\ J \text{ finite type}}} \D(\mb{P}_J) \right) \simeq \colim_{(J_1, \ldots, J_n) \in \mc{C}} \D(\mb{P}_{J_1}) \underset{\D(\mb{I})}{\otimes} \cdots \underset{\D(\mb{I})}{\otimes} \D(\mb{P}_{J_n}). 
	\]
	Here `oblv' refers to forgetting the monoidal structure. This bar construction can be deduced from the material in~\cite{ha} concerning operadic colimits. We will refer to~\cite[5.2]{hecke} for the details of this deduction, but nothing hard (or original) occurs at this step. 
	
	Next, the bar construction for relative tensor products~\cite[4.4.2]{ha} rewrites the previous colimit as follows: 
	\[
		\colim_{(J_1, \ldots, J_n) \in \mc{C}} \D\Big( \mb{P}_{J_1} \overset{\mb{I}}{\times} \cdots \overset{\mb{I}}{\times} \mb{P}_{J_n} \Big). 
	\]
	As usual, we define the $\infty$-category of \emph{prestacks} to be $\on{Fun}(\ms{Sch}^{\mr{aff}, \op}_k, \ms{Spaces})$. The displayed colimit should be thought of as the (correctly normalized) $\infty$-category of $\D$-modules on the prestack $\colim \mc{X}$, where the latter colimit is taken in the $\infty$-category of prestacks. 
		
	Let us show that the geometric fibers of $\colim \mc{X} \to \LG$ are trivial in the following sense: for any geometric point $\eta : \Spec F \to \LG$, we claim that the space 
	\[
		\Hom_{\LG}(\Spec F, \colim \mc{X})
	\]
	of maps over $\LG$ is contractible. To show this, rewrite this space as follows: 
	\[
		\colim_{(J_1, \ldots, J_n) \in \mc{C}} \Hom_{\LG}\Big(\Spec F, \mb{P}_{J_1} \overset{\mb{I}}{\times} \cdots \overset{\mb{I}}{\times} \mb{P}_{J_n}\Big)
	\]
	Here, each Hom is the \emph{discrete} set of points in a Bott--Samelson fiber, and the colimit glues these points together according to the group-theoretic maps specified earlier. It is well-known that the points in a Bott--Samelson fiber correspond to galleries (i.e.\ connected sequences of top-dimensional simplices) in the Bruhat--Tits building, see~\cite{gallery}. In~\ref{intro3}, we will show that the space of all galleries has the same homotopy type as the path space of the (non-topological) Bruhat--Tits building. Now the claim follows from the fact that the Bruhat--Tits building is contractible. 
	
	Finally, since the geometric fibers of $\colim \mc{X} \to \LG$ are trivial in the above sense, and the Bott--Samelson maps
	\[
		\mb{P}_{J_1} \overset{\mb{I}}{\times} \cdots \overset{\mb{I}}{\times} \mb{P}_{J_n} \to \LG
	\]
	are proper, the generalized hyperdescent statement discussed in~\ref{intro4} yields the desired equivalence of $\D$-module $\infty$-categories $\D(\LG) \simeq \D(\colim \mc{X})$. 
	
	\subsubsection{What about Hecke categories?} \label{whatabout} 
	The reader may wonder whether this theorem implies the main theorem of~\cite{hecke} which asserts a similar colimit presentation for the affine Hecke category. This line of reasoning meets a fundamental obstruction, and we explain two incarnations of this obstruction below. 
	
	\noindent\emph{First approach}. The proof in~\ref{outline} implies that 
	\[
		\D(\mb{I} \backslash \LG / \mb{I}) \simeq \colim_{(J_1, \ldots, J_n) \in \mc{C}} \D\Big( \mb{I} \backslash \mb{P}_{J_1} \overset{\mb{I}}{\times} \cdots \overset{\mb{I}}{\times} \mb{P}_{J_n} / \mb{I} \Big),  
	\]
	where the colimit is taken in the $\infty$-category of stable presentable $\infty$-categories (or DG-categories). On the other hand, the colimit presentation of the affine Hecke category is equivalent to  
	\[
		\D'(\mb{I} \backslash \LG / \mb{I}) \simeq \colim_{(J_1, \ldots, J_n) \in \mc{C}} \D'\Big( \mb{I} \backslash \mb{P}_{J_1} \overset{\mb{I}}{\times} \cdots \overset{\mb{I}}{\times} \mb{P}_{J_n} / \mb{I} \Big), 
	\]
	where $\D'(-)$ refers to the full subcategory of $\D$-modules which are locally constant on each twisted product of Schubert cells, see~\cite[1.2.5]{hecke}. The left hand sides are equivalent, but the right hand sides are not obviously equivalent, because a colimit of full subcategories is not necessarily a full subcategory of the colimit. 
	
	\noindent\emph{Second approach}. Another tempting strategy is to look at $\infty$-categories of modules. Namely, if $\mc{A}$ is a monoidal DG-category, we define $\mc{A}\mod$ to be the $\infty$-category of left $\mc{A}$-module objects of $\ms{DGCat}$, as in~\cite[Def.\ 4.2.1.13]{ha}. To explain this approach (and why it meets a similar obstruction), we will take the following statements for granted: 
	\begin{enumerate}[label=(\roman*)]
		\item Let $\mc{C}$ be the $\infty$-category whose objects are (not necessarily small) $\infty$-categories which admit small colimits and whose morphisms are functors which preserve small colimits. The $\infty$-category $\mc{C}$ admits small limits. 
		\item Our main theorem implies that 
		\[
			\D(\LG)\mod \simeq \colim_{\substack{J \subset I \\ J \text{ finite type}}} \D(\mb{P}_J)\mod,
		\]
		where the colimit takes place in $\mc{C}$. 
		\item The colimit presentation of the affine Hecke category is equivalent to 
		\[
			\D(\mb{I} \backslash \LG / \mb{I})\mod \simeq \colim_{\substack{J \subset I \\ J \text{ finite type}}} \D(\mb{I} \backslash \mb{P}_J / \mb{I})\mod, 
		\]
		where the colimit takes place in $\mc{C}$. 
	\end{enumerate}
	
	As David Yang explained to us, it can be deduced from~\cite{gun} that, for any $J \subseteq J' \subseteq I$, there are strictly commutative diagrams
	\begin{equation}\tag{DY}\label{dy} 
	\begin{tikzcd}[column sep = 1.4in, row sep = 0.35in]
		\D(\mb{I}\backslash\mb{P}_J/\mb{I})\mod \ar[hookrightarrow]{r}{\D(\mb{P}_J / \mb{I}) \underset{\D(\mb{I}\backslash\mb{P}_J/\mb{I})}{\otimes} (-)} \ar[d, swap, "\D(\mb{I}\backslash\mb{P}_{J'}/\mb{I}) \underset{\D(\mb{I}\backslash\mb{P}_J/\mb{I})}{\otimes}(-)"] & \D(\mb{P}_J)\mod \ar[d, "\D(\mb{P}_{J'}) \underset{\D(\mb{P}_J)}{\otimes} (-)"] \\
		\D(\mb{I}\backslash\mb{P}_{J'}/\mb{I})\mod \ar[hookrightarrow]{r}{\D(\mb{P}_{J'} / \mb{I}) \underset{\D(\mb{I}\backslash\mb{P}_{J'}/\mb{I})}{\otimes} (-)} & \D(\mb{P}_{J'})\mod
	\end{tikzcd}
	\end{equation}
	where the upper horizontal functor is fully faithful, with essential image given by the $\D(\mb{P}_J)$-modules which are generated by their $\mb{I}$-invariant objects under the action of $\D(\mb{P}_J)$. A similar statement applies to the lower horizontal functor, and also to the functor 
	\begin{cd}[column sep = 1.4in] 
		\D(\mb{I} \backslash \LG / \mb{I})\mod \ar[hookrightarrow]{r}{\D(\LG / \mb{I}) \underset{\D(\mb{I}\backslash\LG/\mb{I})}{\otimes} (-)} & \D(\LG)\mod
	\end{cd}
	Taking the colimit (in $\mc{C}$) along the vertical arrows of (\ref{dy}) yields 
	\begin{cd}
		\displaystyle \colim_{\substack{J \subset I \\ J \text{ finite type}}} \D(\mb{I} \backslash \mb{P}_J / \mb{I})\mod \ar[r] \ar[d] & \displaystyle \colim_{\substack{J \subset I \\ J \text{ finite type}}} \D(\mb{P}_J)\mod \ar{d}[rotate = 90, anchor = north]{\sim} \\
		\D(\mb{I} \backslash \LG / \mb{I})\mod \ar[r, hookrightarrow] & \D(\LG)\mod
	\end{cd}
	By (ii), the right vertical arrow is an equivalence. David Yang's observation implies that the lower horizontal arrow is fully faithful. If we also knew that the upper horizontal arrow is fully faithful, then we could show that the left vertical arrow is an equivalence by comparing the essential images of the horizontal functors, which is easy. By (iii), this would prove the colimit presentation of the affine Hecke category. 

	Unfortunately, we do not know that the upper horizontal arrow is fully faithful, due to the same obstruction as before: a colimit of fully faithful functors is not \emph{a priori} fully faithful. 
	
	It is true that a limit of fully faithful functors is fully faithful. In addition, we expect the following statement to be true: 
	\begin{itemize}
		\item[(iv)] For any diagram $\mc{A}_{\bullet} : \mc{I} \to \ms{Alg}(\ms{DGCat})$ of monoidal DG-categories, we have 
		\[
			\lim_{i \in \mc{I}^{\op}} \mc{A}_i\mod \simeq \colim_{i \in \mc{I}} \mc{A}_i\mod, 
		\]
		where the limit and colimit take place in $\mc{C}$. The functors in the limit diagram are given by restriction of modules, while the functors in the colimit diagram are given by induction of modules. 
	\end{itemize}
	Therefore, it is tempting to express the aforementioned colimits of module $\infty$-categories as limits. However, this does not solve the problem, because the diagram obtained from (\ref{dy}) by replacing the vertical arrows by their right adjoints (i.e.\ restriction of modules) is only lax commutative. 
	
	Thus, we believe the main theorem of~\cite{hecke} is not really a hyperdescent theorem, but rather has more to do with the combinatorics of Schubert cells. In contrast, the message of the present paper is that the colimit theorem for $\D(\LG)$ is a hyperdescent theorem.

	\subsection{The path space of a simplicial complex} \label{intro3} 
	
	Let $K$ be a simplicial complex which is pure of dimension $d$. That is, we have a set (also denoted $K$) together with a downward-closed collection of subsets $\simp(K)$, called simplices, such that each simplex is contained in at least one size-$(d+1)$ simplex. Fix maximal simplices $a, b \in \simp(K)$. 
	
	We introduce the \emph{category of galleries from $a$ to $b$}, denoted $\ms{Gal}_K(a, b)$, whose objects are sequences of simplices 
	\[
		(a = C_0, F_1, C_1, F_2, \ldots, F_n, C_n = b)
	\]
	such that, for each $i$, the simplex $C_i$ is maximal and $C_{i-1} \supseteq F_i \subseteq C_i$. The arrows are generated by three basic classes which are analogous to those considered in~\ref{outline}: 
	\begin{cd}[row sep = 4mm, column sep = 0.4in]
		(\ldots, C, F, C', F, C'', \ldots) \ar[d, "\text{skip}"] & (\ldots, C, F, C', \ldots) \ar[d, "\text{specialize}"] & (\ldots, C, \ldots) \ar[d, "\text{stammer}"] \\
		(\ldots, C, F, C'', \ldots) & (\ldots, C, \bar{F}, C', \ldots) & (\ldots, C, C, C, \ldots)
	\end{cd}
	In the second arrow, $\bar{F} \subseteq F$. The word `stammer' comes from~\cite[Rmk.\ 5.2]{gallery}, where a gallery in which no two consecutive chambers are equal is called `non-stammering.' 
	
	Our main result concerning galleries is the following. 
	
	\begin{cor*}[\ref{gal00}] 
		The category $\ms{Gal}_K(a, b)$ is homotopy equivalent to the space of paths from $a$ to $b$ in the geometric realization of $K$.  
	\end{cor*}
	
	The proof is carried out in Sections~\ref{s-fiber} and \ref{s-path}. 
	
	\begin{rmks*}
		\begin{enumerate}[label=(\arabic*)] \item[ ] 
			\item In Section~\ref{s-fiber}, we prove a combinatorial analogue (Theorem~\ref{fiber0}) of the following (true) statement: if $f : X \to Y$ is a map of simplicial sets, and each fiber of the geometric realization $|f| : |X| \to |Y|$ is contractible, then $f$ is a homotopy equivalence. Analogous statements involving more general kinds of maps are well-known, see for example \cite{smale}, but we were not able to use them because they include compactness hypotheses. For our proof, the main idea is that the fiber of $|f|$ over the barycenter of a simplex $\Delta^n \to Y$ is naturally described as an $(n+1)$-simplicial set, i.e.\ a functor $(\mb{\Delta}^{\op})^{\times (n+1)} \to \set$. 
			\item The reader may wonder what happens for a simplicial complex which is not of pure dimension. In~\ref{ss-non-pure}, we give a similar model of the path space (Theorem~\ref{thm-nonpure}) which applies to an arbitrary simplicial complex. This result is not used in this paper.  
			\item For both of these path space models, composition of paths corresponds to concatenation of sequences. Composition of paths is not used in this paper. 
			
			The right way to interpret this structure is that it yields a functor 
			\[
			(\text{simplicial complexes}) \to (\text{simplicially enriched categories}), 
			\]
			where a simplicial complex $K$ goes to the enriched category in which objects are simplices of $K$, and $\Hom(a, b)$ is the nerve of our model for paths from $a$ to $b$. In~\cite[3.1]{dwyerkan}, Dwyer and Kan constructed an analogous functor 
			\[
			(\text{simplicial sets}) \to (\text{simplicially enriched groupoids}), 
			\]
			but we were not able to deduce our path space model from theirs. Their functor has the advantage of landing in groupoids (rather than categories), but it is less explicit than ours because it uses groupoids that are defined via generators and relations. 
		\end{enumerate}	
	\end{rmks*}

	\subsection{Generalized hyperdescent for \texorpdfstring{$\D$}{D}-modules} \label{intro4} 
	To go from topology to $\D$-module theory, we need the following statement: 
	\begin{thm*}[\ref{simp0}]
		Let $\mc{C}$ be a small $\infty$-category, let $\mc{X} : \mc{C} \to \ms{IndSch}_k$ be a diagram of ind-schemes, and let $Y$ be another ind-scheme of ind-finite type. Suppose we are given a map $\colim \mc{X} \to Y$
		such that $\mc{X}(c) \to Y$ is ind-proper for each $c \in \mc{C}$. Furthermore, suppose that the following holds: 
		\begin{itemize}
			\item[$(\textnormal{H}2'')$] For each geometric point $\eta : \Spec F \to Y$, the space 
			\[
				\Hom_Y(\Spec F, \colim \mc{X})
			\]
			is contractible. (The colimit is taken in the category of prestacks defined in~\ref{intro6}.) 
		\end{itemize}
		Then the pullback functor $\D(Y) \to \D(\colim \mc{X})$ is an equivalence. 
	\end{thm*}
	
	Compare this with the following hyperdescent theorem for $\D$-modules: 
	
	\begin{thm*}[\ref{hyper1}] 
		Let $\mc{X} : \mb{\Delta}^{\op} \to \ms{IndSch}_k$ be a simplicial ind-scheme, and let $Y$ be another ind-scheme of ind-finite type. Suppose we are given a map $\colim \mc{X} \to Y$ such that $\mc{X}(\Delta^n) \to Y$ is ind-proper for each $n$. Furthermore, suppose that the following holds: 
		\begin{itemize}
			\item[$(\textnormal{H}2)$] For each geometric point $\eta : \Spec F \to Y$, the simplicial set 
			\[
				\Delta^n \mapsto \Hom_Y(\Spec F, \mc{X}(\Delta^n))
			\]
			is an acyclic Kan complex. 
		\end{itemize}
		Then the pullback functor $\D(Y) \to \D(\colim \mc{X})$ is an equivalence. 
	\end{thm*}
	
	Sam Raskin taught us how to prove Theorem~\ref{hyper1}, and we record his proof in~\ref{ss-hyper}. Dennis Gaitsgory has pointed out that a simpler proof is possible if one assumes that singular homology (for ind-schemes) satisfies hyperdescent. Indeed, one first uses the projection formula for $\D$-modules to reduce the proof of descent to the case of the object $\omega_Y \in \D(Y)$. By passing to geometric fibers, this statement is further reduced to hyperdescent for homology. 
	
	Our idea for deducing Theorem~\ref{simp0} from Theorem~\ref{hyper1} is very simple. Given $\mc{X} : \mc{C} \to \ms{IndSch}_k$ which satisfies $(\textnormal{H}2'')$, it suffices to find another diagram $\wt{\mc{X}} : \mb{\Delta}^{\op} \to \ms{IndSch}_k$ which satisfies $(\textnormal{H}2)$ together with a map $\colim \mc{X} \to \colim \wt{\mc{X}}$ which is a sieve for the Zariski topology. To accomplish this, we first use~\cite[Ch.\ XII, \S 5]{bk} to replace $\mc{X}$ by a simplicial diagram with the same colimit. Then we get $\wt{\mc{X}}$ by adapting the usual `small object' fibrant replacement~procedure. 
	
	\begin{rmks*}
		\begin{enumerate}[label=(\arabic*)] \item[ ]
			\item Our motivation for Theorem~\ref{simp0} was the observation that Theorem~\ref{hyper1} is unnatural in the following sense: the conclusion refers only to the prestack $\colim \mc{X}$, but the hypothesis refers to a particular colimit presentation of it. Our Theorem~\ref{simp0} reduces this discrepancy but does not eliminate it. Indeed, the hypothesis that each map $\mc{X}(c) \to Y$ is ind-proper may become false upon passing from one colimit presentation of $\colim\mc{X}$ to another. 
			\item The reader may wonder how to apply Theorem~\ref{simp0} to prove the main theorem, since $\LG$ is not ind-finite type. In fact, we apply Theorem~\ref{simp0} only after taking the quotient of everything by an arbitrarily deep congruence subgroup.
		\end{enumerate}
	\end{rmks*}
	
	\subsection{Two more applications} \label{intro5} 
	
	In Section~\ref{s-app}, we also give two applications of hyperdescent (Theorem~\ref{simp0}) which do not depend on our main theorem about $\D(\LG)$. Both are inspired by Varshavsky's talks~\cite{varshavsky}. 
	
	Here is the first one: 
	
	\begin{thm*}[\ref{var0}] 
		We have $\displaystyle\colim_{\textnormal{finite type } J \subset I} \D(\LG/ \mb{P}_J) \simeq \ms{Vect}_k$. 
	\end{thm*}
	
	This result was first proved by Varshavsky, in~\cite{varshavsky}, using a different method. Namely, he decomposes each $\infty$-category $\D(\LG / \mb{P}_J)$ into cells using the Schubert stratification and then proceeds cell-by-cell to compute the colimit. Thus, his proof is parallel to the explicit simplex-by-simplex proof of contractibility of the Bruhat--Tits building (see~\cite[Thm.\ 2.16]{mitchell}). The proof of Theorem~\ref{var0} presented in this paper is simpler: it uses nothing but hyperdescent and the mere fact that the (non-topological) Bruhat--Tits building is contractible. 
	
	In~\cite{varshavsky}, Varshavsky also deduced from Theorem~\ref{var0} a statement expressing the $\LG$-invariants of a $\D(\LG)$-module $\infty$-category in terms of the $\mb{P}_J$-invariants for various $J$. For completeness, we record this as Corollary~\ref{var1} and sketch Varshavsky's proof of this result. Note, however, that this result also follows directly from our main theorem about $\D(\LG)$. 
	
	The second application concerns the affine Springer maps 
	\[
		\LG \overset{\mb{P}_J}{\times} \mb{P}_J \to \LG
	\]
	where the upper $\mb{P}_J$ acts on $\LG$ from the right and on $\mb{P}_J$ by conjugation. 
	
	\begin{thm*}[\ref{spr0}] 
		Let $\LG' \hra \LG$ be the (closed) union of the images of the affine Springer maps. Then $\displaystyle\colim_{\textnormal{finite type } J \subset I} \D\Big( \LG \overset{\mb{P}_J}{\times} \mb{P}_J  \Big) \simeq \D(\LG')$. 
	\end{thm*}
	
	This result did not appear in~\cite{varshavsky}, but it is inspired by Varshavsky's construction of a $W_I$-action on the affine Springer sheaf. Varshavsky's construction uses~\ref{intro1}(\ref{e2}) to reduce to constructing compatible actions by each finite standard subgroup $W_J$, and these actions are provided by finite type Springer theory. 
	
	To prove Theorem~\ref{spr0}, we first note that the affine Springer fibers over $g \in \LG$ are the fixed points of $g$ acting on the partial flag varieties $\LG / \mb{P}_J$. Therefore, to apply hyperdescent, we essentially have to show that the fixed subcomplex of $g$ acting on the Bruhat--Tits building is contractible, which follows from the theory of buildings. However, a complication arises. To apply hyperdescent, we must first reduce to a finite type situation; this requires considering `partial' affine Springer fibers, and a slightly stronger input from building theory will be needed to show the requisite contractibility. 
	
	\subsection{Notations and conventions} \label{intro6} 
	We use both simplicial sets and simplicial complexes (see~\ref{intro3}), and the reader is urged to keep these concepts separate. Let $\mb{\Delta}$ be the simplex category, i.e.\ the category of finite nonempty ordered sets. 
	
	We use the model of $\infty$-categories as quasicategories which is developed in~\cite{htt}. Following the terminology of that book, we say `$\infty$-category' whenever we mean a quasicategory. However, we do not notationally distinguish between a category (i.e.\ `1-category') and its simplicial nerve. In other words, a category is a simplicial set which admits \emph{unique} lifts along inner horn inclusions $\Lambda^n_k \to \Delta^n$ where $0 < k < n$, see~\cite[Prop.\ 1.1.2.2]{htt}. 
	
	Hence, we say that two $\infty$-categories are \emph{homotopy equivalent} when their simplicial nerves are weakly equivalent in the Kan model structure. In particular, a $\infty$-category is \emph{contractible} if its simplicial nerve is contractible. The word `equivalence' by itself always means `categorical equivalence' rather than `homotopy equivalence.' 
	
	For functors, we use the terms `initial/final' rather than `initial/cofinal.' If a pair of adjoint functors is denoted by $\rightleftarrows$, then the left adjoint is always at the top. 
	
	In the later sections, we fix an arbitrary field $k$ of characteristic zero and consider the category of affine $k$-schemes $\ms{Sch}^{\mr{aff}}_k$. We then have the $\infty$-category of prestacks
	\[
		\ms{PreStk}_k := \Fun(\ms{Sch}^{\mr{aff,op}}_k, \ms{Spaces})
	\]
	which contains the full subcategories of ind-schemes $\ms{IndSch}_k$ and schemes $\ms{Sch}_k$. For us, `finite type' includes a hypothesis of separability. This separability hypothesis is only used in two ways in Section~\ref{s-hyper}: a map from an affine scheme to a finite type scheme is affine, and a map from a proper scheme to a finite type scheme is proper.
	
	We use the theory of $\D$-modules developed in~\cite{crystals} and~\cite{raskin}. When discussing $\D$-module $\infty$-categories, we use the right $t$-structure with cohomological indexing conventions. 
	
	\subsection{Acknowledgments} \label{intro7} 
	This paper owes its existence to a number of people who have assisted and inspired us. We would like to thank Dennis Gaitsgory for asking a question which prompted the writing of this paper, and for numerous helpful conversations related to this problem. In addition, Sam Raskin and Dennis Gaitsgory taught us how to prove hyperdescent for $\D$-modules (see~\ref{ss-hyper}), David Yang taught us the material mentioned in~\ref{whatabout} and helped us think through the relation with Hecke categories, and Sanath Devalapurkar and Artem Prikhodko provided a helpful reference to~\cite{ha} during Gaitsgory's office hours. We are grateful to Yakov Varshavsky for his illuminating talks~\cite{varshavsky}, for generously providing materials related to those talks, and for patiently answering our questions. Finally, Roman Bezrukavnikov has been an indispensable source of motivation and useful conversation during the course of this project and the previous one~\cite{hecke} from which it developed. The first author is supported by the NSF GRFP, grant no.\ 1122374. 
	
	\section{A fiber theorem for simplicial sets} \label{s-fiber} 
	
	As explained in~\ref{intro3}(1), if a map $f : X \to Y$ of topological spaces is `nice enough,' and the fibers of $f$ are contractible, then $f$ is a homotopy equivalence, see e.g.\ \cite{smale}. Our goal in this section is to prove this statement when $f$ is the geometric realization of a map of simplicial sets. More precisely, we prove a combinatorial reformulation of this statement (Theorem~\ref{fiber0}) which will be applied in Section~\ref{s-path}. 
	
	\subsection{Edgewise subdivision} 
	
	\subsubsection{Motivation} \label{edge-motivate} 
	The goal of the present section is to study of the fibers of the geometric realization of a map of simplicial sets. In this subsection, we consider the case of the projection map 
	\[
		\pr_2 : \Delta^n \times \Delta^r \to \Delta^r. 
	\]
	The fiber of $|\mr{pr}_2|$ over the barycenter of $|\Delta^r|$ is obviously homeomorphic to $|\Delta^n|$. But the standard decomposition of $\Delta^n \times \Delta^r$ into $\binom{n+r}{r}$ nondegenerate simplices induces a nontrivial subdivision of the fiber $|\Delta^n|$ into products of simplices. One can furthermore decompose each product of simplices as a union of simplices in the standard way. The resulting subdivision of $|\Delta^n|$ is the $r$-th edgewise subdivision, which we will define more formally below. 
	
	\subsubsection{Definition} \label{edge-def} 
	We now recall the edgewise subdivision functor introduced in~\cite[Sect.\ 1]{edgewise}. Pictures of this subdivision can be found in~\cite{edgewise2}. 
	
	In this subsection, fix an integer $r > 0$. Let $\star^r : \mb{\Delta} \to \mb{\Delta}$ be the functor which sends a simplex $\Delta^n$ to its $r$-fold join $\Delta^n \star \cdots \star \Delta^n \simeq \Delta^{r(n+1)-1}$. In other words, it sends an ordered set $S$ to its $r$-fold disjoint union, ordered such that the $i$-th copy of $S$ precedes the $j$-th copy of $S$ when $i < j$. 
	
	The subdivision functor $\sd : \sset \to \sset$ is precomposition by $\star^r$, i.e.\ $\sd(K) := K \circ \star^r$. Its right adjoint is the extension functor $\exx : \sset \to \sset$, which can be characterized in the following way: 
	\[
		\Hom_{\sset}(\Delta^n, \exx(K)) \simeq \Hom_{\sset}(\sd(\Delta^n), K). 
	\]
	
	There is a natural transformation $\lambda : \sd \to \Id_{\sset}$ defined as follows. If $K \in \sset$, then any simplex $\sigma : \Delta^n \to \sd(K)$ corresponds to a map $\wt{\sigma} : \star^r(\Delta^n) \to K$. We define $\lambda(K)$ by requiring that it sends $\sigma$ to the simplex $\Delta^n \to K$ obtained by restricting $\wt{\sigma}$ to the $r$-th factor in the join. By adjunction, we get a natural transformation $\rho : \Id_{\sset} \to \exx$. 
	
	\begin{thm} \label{ex-thm} 
		For any $K \in \sset$, the map $\rho(K) : K \to \exx(K)$ is a homotopy equivalence. 
	\end{thm}
	
	\begin{myproof}
		In the course of developing the $\mr{Ex}^{\infty}$ fibrant replacement functor~\cite{kan-ex}, Kan proved an analogue of this result for barycentric subdivision. In fact, the proof of that result which appears in~\cite[\href{https://kerodon.net/tag/00Z1}{Tag 00Z1}]{kerodon} can be modified to apply to this situation. We explain how to prove the following statements, because their proofs in~\cite{kerodon} refer to the structure of the barycentric subdivision. 
		\begin{enumerate}[label=(\roman*)]
			\item \cite[Prop.\ 3.3.4.8]{kerodon} Let $K$ be a simplicial set. Then $\lambda(K) : \sd(K) \to K$ is a weak homotopy equivalence. 
			\item \cite[Prop.\ 3.3.5.3]{kerodon} If $f : K \to L$ is an anodyne morphism of simplicial sets, then the induced map $\sd(f) : \sd(K) \to \sd(L)$ is also anodyne. 
			\item The functor $\sd$ commutes with products. (Used in \cite[Prop.\ 3.3.5.6]{kerodon}.) 
			\item \cite[Prop.\ 3.3.5.8]{kerodon} Let $K$ be a simplicial set. Then the morphisms 
			\[
			\rho(\exx(K)), \exx(\rho(K)) : \exx(K) \to \exx(\exx(K))
			\]
			are homotopic. 
		\end{enumerate}
		
		The proof of (i) in~\cite{kerodon} reduces it to the statement that $\sd(\Delta^n)$ is contractible. For (ii), it suffices to show that $\sd(f)$ is an acyclic monomorphism. Both statements follow from~\cite[Lem.\ 1.1]{edgewise} since the geometric realization functor reflects homotopy equivalences and monomorphisms. 
		
		Point (iii) follows from the `adjunction' definition  
		\[
			\Hom_{\sset}(\Delta^n, \sd(K)) \simeq \Hom_{\sset}(\star^r(\Delta^n), K) 
		\]
		and the universal property of products. 
		
		The proof of (iv) in~\cite{kerodon} reduces it to the following statement: 
		\begin{itemize}
			\item Let $Q$ be any poset, and let $\on{N}(-)$ be the simplicial nerve functor. The maps 
			\[
				\sd(\lambda(\on{N}(Q))), \lambda(\sd(\on{N}(Q))) : \sd(\sd(\on{N}(Q))) \to \sd(\on{N}(Q))
			\]
			are homotopic. Moreover, this homotopy is constructed functorially in $Q$. 
		\end{itemize}
		
		To prove this, let us first describe these two maps explicitly. Giving a map 
		\[
			\sigma : \Delta^n \to \sd(\sd(\on{N}(Q))). 
		\]
		is equivalent to giving a map 
		\[
			\wt{\sigma} : \star^r(\star^r(\Delta^n)) = (\Delta^n \star \cdots \star \Delta^n) \star \cdots \star (\Delta^n \star \cdots \star \Delta^n) \to \on{N}(Q). 
		\]
		Then $\sd(\lambda(\on{N}(Q))) \circ \sigma$ is obtained by restricting $\wt{\sigma}$ to the join of the $r$-th factors within each parenthesis, while $\lambda(\sd(\on{N}(Q))) \circ \sigma$ is obtained by restricting $\wt{\sigma}$ to the last parenthesis. 
		
		Next, we introduce a sequence of functors 
		\[
			h^{(i)} : \sd(\sd(\on{N}(Q))) \to \sd(\on{N}(Q))
		\]
		which interpolate between these two, for $i = 0, \ldots, r$. Given a simplex $\sigma : \Delta^n \to \sd(\sd(\on{N}(Q)))$ as above, $h^{(i)}$ is required to send this to the simplex
		\[
			h^{(i)} \circ \sigma : \Delta^n \to \sd(\on{N}(Q))
		\]
		which corresponds to the map $\star^r \Delta^n \to \on{N}(Q)$ obtained by restricting $\wt{\sigma}$ to the join of the $r$-th factors within each of the first $i$ parentheses, together with the last $(r-i)$ factors in the last parenthesis. By construction, we have 
		\e{
			h^{(0)} &= \sd(\lambda(\on{N}(Q))) \\
			h^{(r)} &= \lambda(\sd(\on{N}(Q))). 
		} 
		
		To conclude, we fix $i \ge 0$ and construct an explicit homotopy 
		\[
			\eta^{(i)} : \Delta^1 \times \sd(\sd(\on{N}(Q))) \to \sd(\on{N}(Q))
		\]
		from $h^{(i)}$ to $h^{(i+1)}$. A simplex 
		\[
			\sigma' : \Delta^n \to \Delta^1 \times \sd(\sd(\on{N}(Q)))
		\]
		corresponds to a map $\Delta^n \to \Delta^1$ (which is encoded by the integer $a$ for which the first $(a+1)$ vertices of $\Delta^n$ map to $0 \in \Delta^1$) and a map $\sigma : \Delta^n \to \sd(\sd(\on{N}(Q)))$ as above. Now $\eta^{(i)}$ is required to send this to the simplex
		\[
			\eta^{(i)} \circ \sigma' : \Delta^n \to \sd(\on{N}(Q))
		\]
		which corresponds to the map $\star^r \Delta^n \to \on{N}(Q)$ given as follows: 
		\begin{itemize}
			\item On the first $(i-1)$ factors of $\Delta^n$, the map is obtained by restricting $\wt{\sigma}$ to the join of the $r$-th factors within each of the first $(i-1)$ parentheses. 
			\item On the subsimplex $\Delta^{\{0, \ldots, a\}}$ of the $i$-th factor of $\Delta^n$, the map  is obtained by restricting $\wt{\sigma}$ to the corresponding subsimplex of the $r$-th factor in the $i$-th parenthesis. 
			\item On the subsimplex $\Delta^{\{a+1, \ldots, n\}}$ of the $i$-th factor of $\Delta^n$, the map  is obtained by restricting $\wt{\sigma}$ to the corresponding subsimplex of the $i$-th factor in the last parenthesis. 
			\item On the last $(n-i)$ factors of $\Delta^n$, the map is obtained by restricting $\wt{\sigma}$ to the join of the last $(n-i)$ factors in the last parenthesis. 
		\end{itemize}
		It is straightforward to check that this is a valid homotopy from $h^{(i)}$ to $h^{(i+1)}$ which is functorial in the poset $Q$. This concludes the proof of (iv). 
	\end{myproof}
	
	\subsection{The fiber theorem} 
	Here is the main result of this section: 
	
	\begin{thm} \label{fiber0} 
		Let $f : K \to L$ be a map of simplicial sets. Assume that, for each simplex $\sigma : \Delta^n \to L$, the simplicial set 
		\[
		\underline{\Hom}_{\Delta^n}\Big(\Delta^n, \Delta^n \underset{L}{\times} K\Big)
		\]
		is contractible. Then $f$ is a homotopy equivalence. 
	\end{thm}
	
	To clarify the notation, we remark that the `relative internal Hom' has the following universal property: for $X \in \sset$, there is a natural bijection between the set 
	\[
		\Hom_{\sset}\Big(X, \underline{\Hom}_{\Delta^n}\Big(\Delta^n, \Delta^n \underset{L}{\times} K\Big)\Big)
	\]
	and the set of dotted arrows making this solid diagram commute: 
	\begin{cd}
		X \times \Delta^n \ar[r, dashed] \ar[d, "\pr_2"] & K \ar[d, "f"] \\
		\Delta^n \ar[r, "\sigma"] & L
	\end{cd}
	
	The rest of this subsection is devoted to proving the theorem. 
	
	\subsubsection{The `category of simplices' construction} The following idea goes back to~\cite[Ch.\ XII, \S 5]{bk}, where it was used to express a colimit over an arbitrary diagram as a colimit over a simplicial diagram. We will see this idea again in~\ref{ss-simp-repl}. 
	
	If $K$ is a simplicial set, let $\mb{\Delta}_{/K}$ be the category of simplices mapping to $K$. Note that the functor $K \mapsto \mb{\Delta}_{/K}$ preserves limits. There is also a `last vertex' map $\ell_K : \mb{\Delta}_{/K} \to K$ which is functorial in $K$. 
	
	\begin{prop*} \label{fiber1} 
		For every simplicial set $K$, the map $\ell_K$ is initial and final. In particular, it is a homotopy equivalence. 
	\end{prop*}
	\begin{myproof}
		This follows from~\cite[Prop.\ 7.1.10, Prop.\ 7.3.15]{cis}. It can also be deduced from~\cite[Prop.\ 4.2.3.8, Prop.\ 4.2.3.14]{htt}, but this is less convenient because the `initial' part of the statement requires some slight modifications. 
	\end{myproof}
	
	\begin{rmk*}
		If $K \in \sset$ is regular, then the nerve of the poset of nondegenerate simplices of $K$ identifies with the barycentric subdivision of $K$. If $K$ is not regular, then the poset of nondegenerate simplices of $K$ need not be homotopy equivalent to $K$. The category of simplices $\mb{\Delta}_{/K}$ and the barycentric subdivision of $K$ are both homotopy equivalent to $K$, but they are not isomorphic.
	\end{rmk*}
	
	\begin{prop}{\cite[Prop.\ 2.4.7]{hecke}} \label{suscat} 
		Let $\mc{C} \to \{0 \xra{\alpha} 1\}$ be a map of categories, and let $\mc{C}_0, \mc{C}_1$ be the fibers of this map. Let $\on{Arr}_\alpha(\mc{C})$ be the full subcategory of $\on{Arr}(\mc{C})$ consisting of arrows which map to $\alpha$. Then we have a homotopy pushout square 
		\begin{cd}
			\ner(\on{Arr}_{\alpha}(\mc{C})) \ar[r] \ar[d] & \ner(\mc{C}_0) \ar[d] \\
			\ner(\mc{C}_1) \ar[r] & \ner(\mc{C}) 
		\end{cd}
	\end{prop}
	
	\begin{lem} \label{fiber2} 
		Let $f : K \to L$ be a map of simplicial sets. Assume that, for each map $\sigma : \Delta^n \to L$, the full subcategory 
		\[
		\mb{\Delta}_{/\Delta^n \underset{L}{\times} K}^{\mr{surj}} \subset \mb{\Delta}_{/\Delta^n \underset{L}{\times} K},
		\]
		consisting of the simplices which project surjectively onto $\Delta^n$, is contractible. Then $f$ is a homotopy equivalence. 
	\end{lem}
	\begin{myproof}
		First, we prove the following statement: 
		\begin{itemize}
			\item[$(\star)$] Let $f : K \to L$ be a map of simplicial sets, let $L' \subset L$ be a simplicial subset, and let $\sigma : \Delta^n \to L$ be a nondegenerate simplex. Assume the following: 
			\begin{enumerate}[label=(\roman*)]
				\item $L$ is obtained from $L'$ by attaching the simplex $\sigma$, i.e.\ 
				\[
					\sigma(\partial \Delta^n) \subset L' \qquad \text{and} \qquad L = L' \cup_{\partial \Delta^n} \Delta^n.
				\]
				\item The map $(\Id_{L'} \times f) : L' \underset{L}{\times} K \to L'$ is a homotopy equivalence. 
				\item $\mb{\Delta}_{/\Delta^n \underset{L}{\times} K}^{\mr{surj}}$ is contractible. 
			\end{enumerate}
			Then $f$ is a homotopy equivalence. 
		\end{itemize}
		By Proposition~\ref{fiber1}, it suffices to show that the induced map $\wt{f} : \mb{\Delta}_{/K} \to \mb{\Delta}_{/L}$ is a homotopy equivalence. There is a functor $\mb{\Delta}_{/L} \to \{0 \xra{\alpha} 1\}$ which sends an object to $1$ if and only if it is a simplex which maps surjectively onto $\sigma$. Let $(\mb{\Delta}_{/L})_0, (\mb{\Delta}_{/L})_1$ be the fibers of this functor. The composition $\mb{\Delta}_{/K} \to \{0 \xra{\alpha} 1\}$ can be characterized in the same way. 
		
		By Proposition~\ref{suscat}, $\mb{\Delta}_{/K}$ and $\mb{\Delta}_{/L}$ are homotopy equivalent to the homotopy pushouts of the top and bottom rows in this diagram: 
		\begin{cd}
			(\mb{\Delta}_{/K})_0 \ar{d} & \on{Arr}_\alpha(\mb{\Delta}_{/K}) \ar[r, "e_K"] \ar[l] \ar[d] & (\mb{\Delta}_{/K})_1 \ar[d] \\
			(\mb{\Delta}_{/L})_0 & \on{Arr}_\alpha(\mb{\Delta}_{/L}) \ar[r, "e_L"] \ar[l] & (\mb{\Delta}_{/L})_1
		\end{cd}
		The left vertical arrow is a homotopy equivalence by assumption (ii), because $(\mb{\Delta}_{/L})_0 \simeq \mb{\Delta}_{/L'}$ and $(\mb{\Delta}_{/K})_0 \simeq \mb{\Delta}_{/L' \times_L K}$. The objects in the right column are contractible: $(\mb{\Delta}_{/K})_1 \simeq \mb{\Delta}_{/\Delta^n \underset{L}{\times} K}^{\mr{surj}}$ is contractible by (iii), and $(\mb{\Delta}_{/L})_1$ has a terminal object. Therefore, it suffices to show that the middle vertical arrow is a homotopy equivalence. 
		
		For a fixed $I := (\Delta^{n_I} \xra{\sigma_I} K) \in (\mb{\Delta}_{/K})_1$, let $Z \subset \Delta^{n_I}$ be the simplicial subset consisting of simplices which do not surject onto $\sigma$. By unwinding the definition of $\on{Arr}_\alpha(\mb{\Delta}_{/K})$, we obtain an equivalence of categories
		\[
		e_K^{-1}(I) \simeq \mb{\Delta}_{/Z}. 
		\]
		Therefore $e_K^{-1}(I)$ is homotopy equivalent to $Z$ and hence to $\partial \Delta^n$. Since $e_K$ is cocartesian, Thomason's theorem on homotopy colimits~\cite{thomason} implies that $\on{Arr}_{\alpha}(\mb{\Delta}_{/K})$ is homotopy equivalent to $\partial \Delta^n$. 
		
		Similarly, $\on{Arr}_\alpha(\mb{\Delta}_{/L})$ is homotopy equivalent to $\partial \Delta^n$. Tracing through the proof, it is clear that the vertical map is a homotopy equivalence. This proves~$(\star)$. 
		
		Finally, the desired statement results from applying transfinite induction to $(\star)$. 
	\end{myproof}
	
	\begin{lem}\label{fiber3} 
		In the situation of Lemma~\ref{fiber2}, assume that $L = \Delta^n$, and take $\sigma : \Delta^n \to L$ to be the identity. Then 
		$
			\underline{\Hom}_{\Delta^n}(\Delta^n, K)$ and $\mb{\Delta}_{/K}^{\mr{surj}}
		$
		are homotopy equivalent. 
	\end{lem}
	\begin{myproof}
		There is an equivalence 
		\[
			\mb{\Delta}_{/\Delta^n}^{\mr{surj}} \simeq \mb{\Delta}^{\times (n+1)}
		\]
		which sends a simplex $\sigma : \Delta^m \to \Delta^n$ to the $(n+1)$-tuple consisting of the fibers of $\sigma$ over the $(n+1)$ vertices of $\Delta^n$. The functor 
		\[
			\mb{\Delta}_{/K}^{\mr{surj}} \to \mb{\Delta}_{/\Delta^n}^{\mr{surj}} \simeq \mb{\Delta}^{\times (n+1)}
		\]
		is a cartesian fibration which is fibered in sets, so straightening yields an $(n+1)$-simplicial set, i.e.\ a functor  
		\[
			F : (\mb{\Delta}^{\times (n+1)})^{\op} \to \set. 
		\]
		By Thomason's theorem on homotopy colimits~\cite{thomason}, $\colim F$ is homotopy equivalent to $\mb{\Delta}_{/K}^{\mr{surj}}$. 
		
		Since $\mb{\Delta}^{\op}$ is sifted~\cite[Lem.\ 5.5.8.4]{htt}, the diagonal functor 
		\[
			\mr{diag} : \mb{\Delta}^{\op} \to (\mb{\Delta}^{\times (n+1)})^{\op}
		\]
		is final. This implies that $\colim F \simeq \colim (F \circ \mr{diag})$, so $\mb{\Delta}_{/K}^{\mr{surj}}$ is homotopy equivalent to the simplicial set $F \circ \mr{diag}$. 
		
		Next, we construct an isomorphism of simplicial sets 
		\begin{equation} \tag{$\ddag$} \label{ddag} 
			\ul{\Hom}_{\Delta^n}(\Delta^n, K) \simeq \exx(F \circ \mr{diag}). 
		\end{equation}
		This will finish the proof, because $\exx(F \circ \mr{diag})$ and $F \circ \mr{diag}$ are homotopy equivalent by Theorem~\ref{ex-thm}. 		
		
		Equivalently, for each $\Delta^a \in \sset$, we construct a bijection
		\[
			\Hom_{\Delta^n}(\Delta^a \times \Delta^n, K) \simeq \Hom_{\Delta^n} (\sd(\Delta^m), F \circ \mr{diag})
		\]
		which is functorial in $\Delta^a$. To do this, we will show that these sets are naturally in bijection: 
		\begin{enumerate}[label=(\Alph*)]
			\item Maps $\Delta^a \times \Delta^n \to K$ over $\Delta^n$. 
			\item Natural transformations 
			\[
				\Hom_{\Delta^n}(-, \Delta^a \times \Delta^n) \Rightarrow \Hom_{\Delta^n}(-, K)
			\]
			of functors $\mb{\Delta}_{/\Delta^n}^{\op} \to \set$. 
			\item Natural transformations 
			\[
				\Hom_{\Delta^n}(-, \Delta^a \times \Delta^n) \Rightarrow \Hom_{\Delta^n}(-, K)
			\]
			of functors $\mb{\Delta}_{/\Delta^n}^{\mr{surj}, \op} \to \set$. 
			\item Natural transformations 
			\[
				\Hom(\star^{n+1}(-), \Delta^a) \Rightarrow \Hom_{\Delta^n}(\star^{n+1}(-), K)
			\]
			of functors $\mb{\Delta} \to \set$. In the right hand side, the map from $\star^{n+1}(-)$ to the base $\Delta^n$ is obtained by sending the $(i+1)$-st factor in the join to $i \in \Delta^n$. 
			\item Natural transformations 
			$
				\Hom(-, \sd(\Delta^a)) \Rightarrow F(\mr{diag}(-))
			$
			of functors $\mb{\Delta} \to \set$. 
			\item Maps $\Delta^a \to \exx(F \circ \mr{diag})$. 
		\end{enumerate}
		Here are the proofs: 
		\begin{enumerate}
			\item[($\text{A} \simeq \text{B}$)] This follows from the definition of simplicial set. 
			\item[($\text{B} \simeq \text{C}$)] This follows from the fact that, in the category $\sset_{/\Delta^n}$, the object $\Delta^a \times \Delta^n$ is the colimit of objects of the form $\Delta^b \xra{p} \Delta^n$ where $p$ is surjective. In fact, $\Delta^a \times \Delta^n$ can be obtained by gluing $\binom{a+n}{a}$ simplices of dimension $a+n$ along codimension one faces; each simplex involved in this gluing procedure maps surjectively to $\Delta^n$. 
			\item[($\text{C} \simeq \text{D}$)] We have already seen that, since $\Delta^{\op}$ is sifted, the functor $\star^{n+1} : \mb{\Delta} \to \mb{\Delta}_{/\Delta^n}^{\mr{surj}}$ is initial and therefore preserves limits. 
			\item[($\text{D} \simeq \text{E}$)] This follows from the definitions of $\sd$ and $F$. 
			\item[($\text{E} \simeq \text{F}$)] Recall that $\exx$ is the right adjoint of $\sd$. 
		\end{enumerate}
		This concludes the construction of~(\ref{ddag}).		
	\end{myproof}
	
	\subsubsection{Proof of Theorem~\ref{fiber0}} \label{fiber0-proof} 
	
	In view of Lemma~\ref{fiber2}, it suffices to show that 
	\[
	\underline{\Hom}_{\Delta^n}\Big(\Delta^n, \Delta^n \underset{L}{\times} K\Big) \qquad \text{ and } \qquad \mb{\Delta}_{/\Delta^n \underset{L}{\times} K}^{\mr{surj}}
	\]
	are homotopy equivalent. This statement follows from Lemma~\ref{fiber3} applied to the projection map $\Delta^n \underset{L}{\times} K \to \Delta^n$. \hfill $\square$
	
	\section{The path space of a simplicial complex} \label{s-path} 
	
	We now establish the relation (Corollary~\ref{gal00}) between galleries in a pure simplicial complex and paths in its geometric realization which was mentioned in~\ref{intro3}. As promised  in~\ref{intro3}(2), we also prove an analogous result which does not require the simplicial complex to be pure (Theorem~\ref{thm-nonpure}). Our starting point is the `combinatorial' model of the path space which was developed in~\cite{stone}. 
	
	In this section, $K$ is a simplicial complex, and $a, b \in \simp(K)$ are two fixed simplices. 
	
	\subsection{Stone's model for the path space} 
	
	\subsubsection{Definition} For each $n \ge 0$, we consider the set $\bm{t}_n := \frac{1}{2^n} \BZ \cap [0, 1]$ of \emph{binary fractions of level $n$}. Let $\bm{t} = \cup_n\, \bm{t}_n$ be the set of all binary fractions. 
	
	Let $\mathsf{Path}_K(a, b)_n \subset \simp(K)^{\times (2^n+1)}$ be the induced subposet consisting of functions 
	\[
		f : \bm{t}_n \to \simp(K)
	\]
	satisfying $f(0) = a$, $f(1) = b$, and $f(\frac{i}{2^n}) \cup f(\frac{i+1}{2^n}) \in \simp(K)$ for all $i$. There is a full embedding of posets
	\[
		\imath_n : \mathsf{Path}_K(a, b)_n \hra \mathsf{Path}_K(a, b)_{n+1}
	\]
	which is defined by `linear interpolation' in the following sense: 
	\[
		(\imath_n f)(\tfrac{i}{2^{n+1}}) = \begin{cases}
			f(\tfrac{i}{2^{n+1}}) & \text{ if $i$ is even} \\
			f(\tfrac{i-1}{2^{n+1}}) \cup f^(\tfrac{i+1}{2^{n+1}}) & \text{ if $i$ is odd}
		\end{cases}
	\]
	We also consider the colimit taken along the full embeddings $\imath_n$: 
	\[
		\mathsf{Path}_K(a, b) := \colim_n \mathsf{Path}_K(a, b)_n
	\] 
	
	\begin{thm} \label{stone-thm} 
		Fix points $p_a \in |a|$ and $p_b \in |b|$. The poset $\mathsf{Path}_K(a, b)$ is homotopy equivalent to the path space $\Omega(p_a, p_b, |K|)$. 
	\end{thm}
	\begin{myproof}
		This is a consequence of~\cite[Thm.\ 1]{stone}. For the reader's convenience, we provide some comments on how the constructions in that paper match up with ours. Stone begins with a sequence of compatible subdivisions of $[0, 1]$ into intervals, which we take to be $\bm{t}_0, \bm{t}_1$, etc. Given a simplicial complex $K$ and an integer $n$, Stone considers the space $N(n)$ of piecewise-linear paths whose vertices lie in $\bm{t}_n$. This space has an obvious CW complex structure whose cells are indexed by $\ms{Path}_K(a, b)_n$. Namely, a function $f \in \ms{Path}_K(a, b)_n$ corresponds to the space of piecewise-linear paths for which the $i$-th vertex lies in the interior of $f(\frac{i}{2^n})$. Since this CW complex structure is regular, and the poset structure of $\ms{Path}_K(a, b)_n$ agrees with the closure relations between cells, $N(n)$ is homotopy equivalent to $\ms{Path}_K(a, b)_n$. 
		
		There is an obvious cellular map $N(n) \to N(n+1)$ whose map on cells is given by the full embedding 
		\[
			\ms{Path}_K(a, b)_n \hra \ms{Path}_K(a, b)_{n+1}. 
		\]
		Now~\cite[Thm.\ 1]{stone} says that the path space $\Omega(p_a, p_b, |K|)$ is homotopy equivalent to $\colim_n N(n)$. By the previous paragraph, this is equivalent to $\colim_n \ms{Path}_K(a, b)_n$, as desired. 
	\end{myproof}
	
	\subsubsectiona  \label{fun-interpret} 
	Here is an alternative description of $\mathsf{Path}_K(a, b)$. It consists of functions 
	\[
		f: [0,1] \to \simp(K) \qquad (\text{equivalently, } f : \bm{t} \to \simp(K))
	\]
	which can be constructed in the following way: 
	\begin{itemize}
		\item Choose a weakly increasing sequence $0 = t_0 \le t_1 \le \cdots \le t_\ell = 1$ of binary fractions satisfying $t_i < t_{i+2}$ for all $i$. 
		\item Choose a sequence $a = \sigma_1, \sigma_2, \ldots, \sigma_\ell = b$ of elements of $\simp(K)$ which satisfies the following for each $i$: 
		\begin{itemize}
			\item $\sigma_i \subsetneq \sigma_{i+1}$ or $\sigma_i \supsetneq \sigma_{i+1}$
			\item If $t_{i-1} = t_i$, then $\sigma_{i-1} \supsetneq \sigma_i \subsetneq \sigma_{i+1}$. 
		\end{itemize}
		\item Construct the function  
		\[
			f(t) := \begin{cases}
				\sigma_i & \text{ if $t \in (t_{i-1}, t_i)$} \\
				\min(\sigma_i, \sigma_{i+1}) & \text{ if $t = t_i$}
			\end{cases}
		\]
	\end{itemize}
	In fact, each element of $\mathsf{Path}_K(a, b)$ arises from this construction in exactly one way. 
	
	\begin{defn*}
		Given an element of $\ms{Path}_K(a, b)$, we refer to the corresponding sequences $(t_0, \ldots, t_\ell)$ and $(\sigma_1, \ldots, \sigma_\ell)$ as its \emph{transition points} and \emph{values}, respectively. 
	\end{defn*}

	\subsection{Fat paths in a pure simplicial complex} \label{ss-fat} 
	In this subsection, assume that $K$ is pure of dimension $d$ and that $a, b \in \simp(K)$ are maximal (i.e.\ $d$-dimensional).
	
	\subsubsection{Definition} Let $\mathsf{Path}_K^{\mr{fat}}(a, b) \subset \mathsf{Path}_K(a, b)$ be the full subposet consisting of functions $\bm{t} \to \simp(K)$ whose sequence of values $\sigma_1, \ldots, \sigma_\ell$ satisfies the following: 
	\begin{itemize}
		\item For each $i$, exactly one of $\sigma_i, \sigma_{i+1}$ is maximal. 
	\end{itemize}
	Since $a$ and $b$ are maximal, this is equivalent to requiring that $\ell$ is odd and that $\sigma_i$ is maximal if and only if $i$ is odd.
	
	\begin{thm} \label{fat0} 
		The embedding $\mathsf{Path}_K^{\mr{fat}}(a, b) \hra \mathsf{Path}_K(a, b)$ is a homotopy equivalence. 
	\end{thm}
	
	The rest of this subsection is devoted to proving the theorem.
	
	\subsubsection{Definition} 
	Let $\mathsf{Path}_K^{\mr{reg}}(a, b) \subset \mathsf{Path}_K(a, b)$ be the subposet consisting of paths whose sequence of transition points is strictly increasing. We also define 
	\[
		\mathsf{Path}_K^{\mr{fat, reg}}(a, b) := \mathsf{Path}_K^{\mr{fat}}(a, b) \cap \mathsf{Path}_K^{\mr{reg}}(a, b).
	\]
	\begin{lem} \label{lem-reg} 
		The embeddings $\mathsf{Path}_K^{\mr{reg}}(a, b) \hra \mathsf{Path}_K(a, b)$ and $\mathsf{Path}_K^{\mr{fat, reg}}(a, b) \hra \mathsf{Path}_K^{\mr{fat}}(a, b)$ are homotopy equivalences. 
	\end{lem}
	\begin{myproof}
		We will prove the first statement. Let $\imath$ be the embedding in question, and fix a path $p \in \mathsf{Path}_K(a, b)_n$. In view of Quillen's Theorem A, it suffices to show that the poset $(p \downarrow \imath)$ is contractible. 
		
		In fact, we will show that $(p \downarrow \imath)$ is cofiltered. Let $(t_0, \ldots, t_\ell)$ and $(\sigma_1, \ldots, \sigma_\ell)$ be the transition points and values of $p$. For each binary fraction $0 < r < \frac{1}{2^{n+1}}$, let $p(r)$ be the path obtained by modifying the transition points as follows: 
		\begin{itemize}
			\item For each $i$ such that $t_{i-1} = t_i$, we make the replacements $t_{i-1} \rightsquigarrow t_{i-1} - r$ and $t_i \rightsquigarrow t_i + r$. 
		\end{itemize}
		It is clear that $p(r) \in \mathsf{Path}_K^{\mr{reg}}(a, b)$, and there is a map $p \to p(r)$. 
		
		If $q_1, \ldots, q_k$ are a finite collection of objects in $(p \downarrow \imath)$, i.e.\ paths dominated by $p$, then for sufficiently small $r$, the path $p(r)$ also dominates $q_1, \ldots, q_k$. This shows that $(p \downarrow \imath)$ is cofiltered, as desired. 
		
		The second statement is proved in the same way. (If $p$ is fat, then so is $p(r)$.)
	\end{myproof}
	
	\subsubsectiona For the rest of this subsection, fix a path $p \in \mathsf{Path}_K^{\mr{reg}}(a, b)_n$ which has transition points $(t_0, \ldots, t_\ell)$ and values $(\sigma_1, \ldots, \sigma_\ell)$. Let $\imath^{\mr{fat}} : \mathsf{Path}_K^{\mr{fat, reg}}(a, b) \hra \mathsf{Path}_K^{\mr{reg}}(a, b)$ be the embedding from before. Fix an integer $m \ge n$, and let $(\imath^{\mr{fat}} \downarrow p)_m \subset (\imath^{\mr{fat}} \downarrow p)$ consist of paths $q$ lying over $p$ such that $q \in \mathsf{Path}_K^{\mr{fat,reg}}(a, b)_m$. 
	
	\subsubsectiona \label{fat1} 
	Let 
	\[
		G_1 : (\imath^{\mr{fat}} \downarrow p)_m \to (\imath^{\mr{fat}} \downarrow p)_{m+2}
	\]
	be the functor which sends a path $q \in (\imath^{\mr{fat}} \downarrow p)_m$ to the path defined as follows: 
	\begin{itemize}
		\item For each $i = 1, \ldots, \ell$, if $t_i$ is also a transition point of $q$, then that transition point is modified as follows: 
		\begin{itemize}
			\item If $q$ increases at $t_i$, then the transition point is shifted left by $\frac{1}{2^{m+2}}$. 
			\item If $q$ decreases at $t_i$, then the transition point is shifted right by $\frac{1}{2^{m+2}}$. 
		\end{itemize}
	\end{itemize} 
	There is a natural transformation $\Id_{(\imath^{\mr{fat}} \downarrow p)_m} \to G_1$. 
	
	\subsubsectiona \label{fat2} 
	Let $I \subset \{1, \ldots, \ell\}$ be the subset of indices such that $i \in I$ if and only if neither $\sigma_i$ nor $\sigma_{i+1}$ is maximal. (These are the places where $p$ is `not fat.') Let $\mc{C} \subset (\imath^{\mr{fat}} \downarrow p)$ be the full subcategory such that $q \in \mc{C}$ if and only if there exists a path $q' \in \Im(G_1)$ such that $q$ can be obtained from $q'$ by the following procedure: 
	\begin{itemize}
		\item For each $i \in I$, change a nonzero set of values on the domain 
		\[
			(t_i-\tfrac{1}{2^{m+3}}, t_i-\tfrac{1}{2^{m+4}})
		\]
		to maximal simplices, and change a nonzero set of values on the domain 
		\[
			(t_i+\tfrac{1}{2^{m+4}}, t_i+\tfrac{1}{2^{m+3}})
		\]
		to maximal simplices. 
	\end{itemize} 
	For a fixed path $q \in (\imath^{\mr{fat}} \downarrow p)$, such a path $q'$ is unique if it exists. 
	
	Let 
	\[
		G_2 : \mc{C} \to \Im(G_1)
	\]
	be the functor which sends $q \mapsto q'$. It is clear that $G_2$ is a cocartesian fibration. In addition, there is a natural transformation $\Id_{\mc{C}} \to G_2$ of functors $\mc{C} \to (\imath^{\mr{fat}} \downarrow p)$. 
	\begin{lem*}
		The functor $G_2 : \mc{C} \to \Im(G_1)$ is a homotopy equivalence. 
	\end{lem*}
	\begin{myproof}
		Since $G_2$ is a cocartesian fibration, it suffices to show that, for each path $q' \in \Im(G_1)$, the fiber $G_2^{-1}(q')$ is contractible. 
		
		For any set $S$, let $\mc{J}(S)$ denote the poset of nonzero functions $[0, 1] \to \{0\} \sqcup S$ which can be constructed in the following way: 
		\begin{itemize}
			\item Choose a sequence $0 = c_0 \le c_1 < c_2 < \cdots < c_{2r-1} < c_{2r} \le 1$ of binary fractions. 
			\item Choose a sequence $s_1, \ldots, s_r$ in $S$. 
			\item Construct the function 
			\[
				f(t) := \begin{cases}
					s_k &\text{ if } t \in (c_{2k-1}, c_{2k}) \\
					0 & \text{ if } t \in [c_{2k}, c_{2k+1}]. 
				\end{cases}
			\]
		\end{itemize}
		For each $i \in I$, let $S_i$ be the set of maximal simplices containing $q'(t_i)$ in their closures. Then we have an isomorphism of posets
		\[
			G_2^{-1}(q') \simeq \prod_{i \in I} \mc{J}(S_i)^{\times 2}. 
		\]
		Indeed, each $q \in G_2^{-1}(q')$ is specified by choosing, for each interval in the domain considered at the start of this section, a (nontrivial) way of increasing some of the values to be maximal simplices. Thus, it suffices to show that $\mc{J}(S)$ is contractible, for any nonempty set $S$. 
		
		Now we show that $\mc{J}(S)$ is contractible. For any integer $x \ge 1$, define the full subcategory $\mc{J}(S)_x \subset \mc{J}(S)$ to consist of functions whose transition points $(c_0, \ldots, c_{2r})$ are contained in $\bm{t}_x$. We have 
		\[
			\mc{J}(S) \simeq \colim_x \mc{J}(S)_x, 
		\]
		so it suffices to show that each embedding $\imath_x :  \mc{J}(S)_x \hra  \mc{J}(S)_{x+2}$ is homotopy equivalent to a constant functor. 
		
		Choose an arbitrary element $s \in S$. Then define the functors 
		\[
			\imath_x^{(2)}, \imath_x^{(3)}, \imath_x^{(4)} : \mc{J}(S)_x \hra  \mc{J}(S)_{x+2}
		\]
		as follows. 
		\begin{itemize}
			\item Given a function $f \in \mc{J}(S)_x$, the function $\imath_x^{(2)}(f)$ is obtained by redefining $f|_{[0, \frac{1}{2^{x+1}}]}$ to be the zero function. 
			\item Given a function $f \in \mc{J}(S)_x$, the function $\imath_x^{(3)}(f)$ is obtained by redefining $f|_{[0, \frac{1}{2^{x+1}}]}$ to be the function which sends $(0, \frac{1}{2^{x+2}})$ to $s$ and is zero elsewhere. 
			\item The functor $\imath_x^{(4)}$ the constant functor whose output is the function $[0, 1] \to \{0\} \sqcup S$ which sends $(0, \frac{1}{2^{x+2}})$ to $s$ and is zero elsewhere. 
		\end{itemize}
		There are natural transformations $\imath_x \rightarrow \imath_x^{(2)} \leftarrow \imath_x^{(3)} \rightarrow \imath_x^{(4)}$. Therefore $\imath_x$ is homotopy equivalent to a constant functor, as desired.  
	\end{myproof}
	
	\subsubsectiona  \label{fat3} 
	Define the functor 
	\[
		G_3 : \mc{C} \to (\imath^{\mr{fat}} \downarrow p)
	\]
	as follows: if $q \in \mc{C}$ is a path, then $G_3(q)$ is the path obtained via this modification: 
	\begin{itemize}
		\item For $i \in I$, let $(c_i, d_i)$ be the largest interval containing $t_i$ on which $q$ is constant. By construction, this interval contains $(t_i - \frac{1}{2^{m+4}}, t_i + \frac{1}{2^{m+4}})$. Redefine $q|_{(c_i, d_i)}$ to have constant value $\max(\sigma_i, \sigma_{i+1})$. 
	\end{itemize}
	There is a natural transformation $\Id_{\mc{C}} \to G_3$. 
	
	\subsubsectiona  \label{fat4} 
	For each $i \in I$, choose a maximal simplex $s_i \in \simp(K)$ which contains $\max(\sigma_i, \sigma_{i+1})$. Define the functor 
	\[
		G_4 : \Im(G_3) \to (\imath^{\mr{fat}} \downarrow p)
	\]
	as follows: if $q \in \Im(G_3)$ is a path, then $G_4(q)$ is the path obtained via this modification: 
	\begin{itemize}
		\item For $i \in I$, the restriction $q|_{(t_i - \frac{1}{2^{m+4}}, t_i + \frac{1}{2^{m+4}})}$ is constant with value $\max(\sigma_i, \sigma_{i+1})$, by construction. Redefine the smaller restriction $q|_{(t_i - \frac{1}{2^{m+5}}, t_i + \frac{1}{2^{m+5}})}$ to be constant with value $s_i$. 
	\end{itemize}
	There is a natural transformation $G_4 \to \Id_{\Im(G_3)}$. 
	
	\subsubsectiona \label{fat5} 
	Define the functor 
	\[
		G_5 : \Id_{\Im(G_4)} \to (\imath^{\mr{fat}} \downarrow p)
	\]
	to be constant with value equal to the path $q_0$ defined as follows: 
	\begin{itemize}
		\item For $i \in I$, the restriction $q_0|_{(t_i - \frac{1}{2^{m+5}}, t_i + \frac{1}{2^{m+5}})}$ is constant with value $s_i$. Everywhere else, $q_0$ agrees with $p$. 
	\end{itemize}
	There is a natural transformation $\Id_{\Im(G_4)} \to G_5$. 
	
	\subsubsection{Proof of Theorem~\ref{fat0}} 
	In view of Lemma~\ref{lem-reg}, it is equivalent to show that 
	\[
		\imath^{\mr{fat}} : \ms{Path}_K^{\mr{fat, reg}}(a, b) \hra \ms{Path}_K^{\mr{reg}}(a,b)
	\]
	is a homotopy equivalence. By Quillen's Theorem A, it suffices to show, for any $p \in \ms{Path}_K^{\mr{reg}}(a,b)$, that the category $(\imath^{\mr{fat}} \downarrow p)$ is contractible. If $n$ is such that $p \in \ms{Path}_K^{\mr{reg}}(a,b)_n$, then we have 
	\[
		(\imath^{\mr{fat}} \downarrow p) \simeq \colim_{m \ge n} (\imath^{\mr{fat}} \downarrow p)_m. 
	\]
	Therefore, it suffices to show that each embedding 
	\[
		(\imath^{\mr{fat}} \downarrow p)_m \hra (\imath^{\mr{fat}} \downarrow p)
	\]
	is homotopic to a constant map, i.e.\ that this subcategory can be contracted to a point inside $(\imath^{\mr{fat}} \downarrow p)$. By~\ref{fat1}, this embedding is homotopic to the functor $G_1$. Therefore, it suffices to show that $\Im(G_1)$ can be contracted to a point inside $(\imath^{\mr{fat}} \downarrow p)$. 
	
	In~\ref{fat2}, the natural transformation $\Id_{\mc{C}} \to G_2$ yields a lax-commutative diagram 
	\begin{cd}
		\mc{C} \ar[rr, hookrightarrow] \ar[rd, swap, "G_2"] \ar[bend right = 20, draw=none]{rr}[anchor=center]{\Downarrow} & & (\imath^{\mr{fat}} \downarrow p) \\
		& \Im(G_1) \ar[ru, hookrightarrow] 
	\end{cd}
	Lemma~\ref{fat2} says that $G_2$ is a homotopy equivalence. Therefore, it suffices to show that the embedding $\mc{C} \hra (\imath^{\mr{fat}} \downarrow p)$ is homotopic to a constant functor. 
	
	In the same way as before, \ref{fat3}, \ref{fat4}, and \ref{fat5} collectively show that this embedding is homotopic to the constant functor with value $q_0$, as desired. \hfill $\square$

	\subsection{Galleries in pure simplicial complexes} \label{s-gal} 
	As in~\ref{ss-fat}, we assume that $K$ is pure of dimension $d$ and that $a, b \in \simp(K)$ are maximal. 
	
	\subsubsection{Definition} \label{gal-def} 
	In~\ref{intro3}, we sketched how to define a category $\ms{Gal}_K(a, b)$. Now we give the precise definition. The objects, which are called \emph{galleries}, are sequences of simplices 
	\[
		(a=C_0, F_1, C_1, F_2, \ldots, F_n, C_n = b)
	\]
	satisfying the following properties: 
	\begin{itemize}
		\item Each $C_i$ is a maximal simplex. 
		\item For each $i \ge 1$, we have $C_{i-1} \supseteq F_i \subseteq C_i$. 
	\end{itemize}
	A morphism of galleries 
	\[
		(C_0, F_1, C_1, F_2, \ldots, F_n, C_n) \to (C_0', F_1', C_1', F_2', \ldots, F_m', C_m')
	\]
	is a weakly increasing map $f : \{1, \ldots, n\} \to \{1, \ldots, m\}$ which satisfies the following: 
	\begin{itemize}
		\item For each $j \in \{1, \ldots, m\}$, we have $F_j' \subseteq \bigcap_{i \in f^{-1}(j)} F_i$. 
		\item For each $i \in \{0, 1, \ldots, n\}$ and $j \in \{F(i), F(i)+1, \ldots, F(i+1)-1\}$, we have $C_i = C_j'$. (By convention, $F(0) = 0$.) 
	\end{itemize}
	
	\subsubsectiona \label{gal-functor} 
	Define the functor 
	\[
		G : \mathsf{Path}_K^{\mr{fat}}(a, b) \to \gal_K(a, b)
	\]
	as follows: 
	\begin{itemize}
		\item If a path $p$ has sequence of values $(\sigma_1, \ldots, \sigma_\ell)$, then $G(p)$ equals $(\sigma_1, \ldots, \sigma_\ell)$ interpreted as a gallery. 
		\item Given a map of paths $p \to p'$, where $p$ has transition points $(t_0, \ldots, t_\ell)$ and values $(\sigma_1, \ldots, \sigma_\ell)$, and $p'$ has transition points $(t_0', \ldots, t_r')$ and values $(\sigma_1', \ldots, \sigma_r')$, define $G(p \to p')$ to be the map of galleries 
		\[
			(\sigma_1, \ldots, \sigma_\ell) \to (\sigma_1', \ldots, \sigma_\ell')
		\]
		specified by the weakly increasing map $f : \{1, \ldots, \frac{\ell-1}{2}\} \to \{1, \ldots, \frac{r-1}{2}\}$ which is given as follows:  
		\begin{itemize}
			\item For each $i$, the interval $[t_{2i-1}, t_{2i}]$ is contained in $[t_{2j-1}, t_{2j}]$ for a unique $j$. Define $f(i) := j$. 
		\end{itemize}
	\end{itemize} 
	
	\begin{thm} \label{gal0} 
		The functor $G : \mathsf{Path}_K^{\mr{fat}}(a, b) \to \gal_K(a, b)$ is a homotopy equivalence. 
	\end{thm}
	
	\begin{cor} \label{gal00} 
		The categories $\mathsf{Path}_K(a, b)$ and $\gal_K(a, b)$ are homotopy equivalent. 
	\end{cor}
	\begin{myproof}
		In the diagram 
		\[
			\mathsf{Path}_K(a, b) \hookleftarrow \mathsf{Path}_K^{\mr{fat}}(a, b) \xra{G} \gal_K(a, b),
		\]
		the left and right arrows are homotopy equivalences by Theorems~\ref{fat0} and~\ref{gal0}, resp. 
	\end{myproof}
	
	The rest of this subsection is devoted to proving the theorem.
	
	\subsubsection{Definition} Define the full subcategory $\gal_K^{\mr{nu}}(a, b) \subset \gal_K(a, b)$ to consists of galleries $(C_0, F_1, C_1, \ldots, F_n, C_n)$ such that each $F_i$ is non-maximal. It is easy to see that 
	\[
		\Im(G) = \gal_K^{\mr{nu}}(a, b). 
	\]	
	\begin{lem} \label{gal1} 
		The embedding $\gal_K^{\mr{nu}}(a, b) \hra \gal_K(a, b)$ is a homotopy equivalence. 
	\end{lem}
	\begin{myproof}
		This embedding has a right adjoint 
		\[
			R : \gal_K(a, b) \to \gal_K^{\mr{nu}}(a, b)
		\]
		defined as follows. If $(C_0, F_1, \ldots, F_n, C_n)$ is a gallery, then $R$ sends it to the gallery obtained by performing these modifications: 
		\begin{itemize}
			\item For each $i \ge 1$, if $F_i$ is maximal, then $C_i = F_i = C_{i+1}$. Delete $F_i$ and $C_{i+1}$ from the sequence. (If $F_i$ is not maximal, do nothing at the $i$-th step.) 
		\end{itemize}
		Quillen's Theorem A implies that any adjoint functor is a homotopy equivalence. 
	\end{myproof}
	
	\subsubsectiona  \label{gal2} 
	Let $n > 0$ be an integer, and choose tuples $(d_i)_{i=1}^n \in \{\mr{strict}, \mr{weak}\}^{\times n}$ and $(e_i)_{i=1}^{n-1} \in \{-1, +1\}^{\times (n-1)}$. Define the poset $\mc{P}$ as follows: 
	\begin{itemize}
		\item The objects are tuples $(0 = t_0, t_1, \ldots, t_n = 1)$ of binary fractions such that, for each $i$, we have $t_{i-1} < t_i$ if $d_i = \mr{strict}$ and $t_{i-1} \le t_i$ if $d_i = \mr{weak}$. 
		\item There is an arrow 
		\[
			(t_0, t_1, \ldots, t_n) \to (t_0', t_1', \ldots, t_n')
		\]
		if and only if, for each $i$, we have $t_i' \le t_i$ if $e_i = -1$ and $t_i' \ge t_i$ if $e_i = +1$. 
	\end{itemize}
	\begin{lem*}
		The poset $\mc{P}$ is contractible. 
	\end{lem*}
	\begin{myproof}
		Define $\mc{P}_m := \mc{P} \cap (\bm{t}_m)^{\times (n+1)}$, so that $\mc{P} \simeq \colim_m \mc{P}_m$. We will show that each embedding $\imath_m : \mc{P}_m \hra \mc{P}$ is homotopic to a constant functor. 
		
		Define the functor 
		\[
			\imath_m' : \mc{P}_m \hra \mc{P}_{m+1}
		\]
		as follows: 
		\[
			\imath_m'((t_0, t_1, t_2, \ldots, t_n)) := (t_0, t_1+\tfrac{1}{2^{m+1}}, t_2, \ldots, t_n). 
		\]
		If $e_i = -1$, there is a natural transformation $\imath_m \leftarrow \imath_m'$. If $e_i = +1$, there is a natural transformation $\imath_m \rightarrow \imath_m'$. In either case, $\imath_m$ and $\imath_m'$ are homotopic, so it suffices to show that the embedding $\Im(\imath_m') \hra \mc{P}$ is homotopic to a constant functor. 
		
		Let $\mc{P}_{m+1}^{\circ} \subset \mc{P}_{m+1}$ be the full subposet consisting of tuples $(t_0, t_1, \ldots, t_n)$ satisfying $t_1 > t_0$. We have $\Im(\imath_m') \subset \mc{P}_{m+1}^{\circ}$, so it suffices to show that the embedding $\mc{P}_{m+1}^{\circ} \hra \mc{P}$ is homotopic to a constant functor. 
		
		Let $r$ be an integer such that $2^r > n$. For each $j = 0, 1, \ldots, (n-1)$, define the functor 
		\[
			\eta_j : \mc{P}_{m+1}^{\circ} \hra \mc{P}
		\]
		as follows: 
		\[
			\imath_m^{(j)}((t_0, \ldots, t_n)) := (\tfrac{0}{2^{m+r+1}}, \tfrac{1}{2^{m+r+1}}, \ldots, \tfrac{j}{2^{m+r+1}}, t_{j+1}, \ldots, t_n). 
		\]
		To check that this is a functor, use the fact that $t_{j+1} \ge \frac{1}{2^{m+1}}$. (The purpose of the two previous paragraphs was to ensure that this inequality is true.) 
		
		Fix $j \ge 1$. If $e_j = -1$, we have a natural transformation $\eta_{j-1} \to \eta_j$. If $e_j = +1$, we have a natural transformation $\eta_{j-1} \leftarrow \eta_j$. Thus the $\eta_j$ are homotopic to one another. This finishes the proof because $\eta_0$ is the embedding while $\eta_{n-1}$ is a constant functor. 
	\end{myproof}
	
	\begin{rmk*}
		We have stated this lemma in great generality because we will use it two more times in~\ref{ss-non-pure}. In some sense, this lemma is the crux of how we pass from Stone's model for the path space to more flexible models involving sequences of simplices. It may be helpful to compare this lemma with the following statements: 
		\begin{enumerate}[label=(\roman*)]
			\item The poset of closed intervals in $S^1$ is homotopy equivalent to $S^1$. 
			\item Let $M$ be a manifold. The poset of closed balls embedded into $M$ is homotopy equivalent to $M$. 
			\item Let $K$ be a regular simplicial set. The poset of nondegenerate simplices of $K$ is homotopy equivalent to $K$. (See~\ref{fiber1}.) 
		\end{enumerate}
		These statements say that \emph{movement} in a space can be captured by \emph{containment} of subsets of that space -- even though the subsets do not move. Similarly, Lemma~\ref{gal2} will be used to `move' transition points using the partial order on functions. Once the transition points can move, the distinction between Stone's model and ours is blurred. 
		
		Points (i) and (ii) are proved using a hypercover theorem which is analogous to the one we develop in Theorem~\ref{simp0}, and a proof of (iii) can be deduced from~\cite[Prop.\ 4.2.3.8]{htt}, which uses a similar `contractibility of fibers' idea. In Segal's framework for conformal field theory, (ii) can be used to obtain the parallel transport of point operators; this lends further credence to the intuition that these statements encode movement in a space. 		
	\end{rmk*}
	
	\subsubsection{Proof of Theorem~\ref{gal0}} 
	
	In view of Lemma~\ref{gal1}, it suffices to show that the functor 
	\[
		G : \ms{Path}^{\mr{fat}}_K(a, b) \to \gal^{\mr{nu}}_K(a, b)
	\]
	is a homotopy equivalence. By Theorem~\ref{fiber0}, it suffices to prove the following statement: 
	\begin{itemize}
		\item[$(\dagger)$] For any simplex $\sigma : \Delta^n \to \gal^{\mr{nu}}_K(a, b)$, the simplicial set 
		\[
			\underline{\Hom}_{\Delta^n}\Big(\Delta^n, \Delta^n \underset{\gal^{\mr{nu}}_{K}(a, b)}{\times} \ms{Path}^{\mr{fat}}_K(a, b) \Big)
		\]
		is contractible. 
	\end{itemize}
	Since the domain and target of $G$ are 1-categories, this simplicial set is the 1-category $\mc{C}$ defined as follows: 
	\begin{itemize}
		\item We think of $\sigma$ as being a composable sequence of maps $g_0 \to \cdots \to g_n$ in $\gal^{\mr{nu}}_K(a, b)$. 
		\item An object of $\mc{C}$ is a composable sequence of maps $p_0 \to \cdots \to p_n$ in $\ms{Path}^{\mr{fat}}_K(a, b)$ whose image under $G$ is identified with $\sigma$. 
		\item A morphism in $\mc{C}$ is a commutative diagram 
		\begin{cd}
			p_0 \ar[r] \ar[d] & p_1 \ar[r] \ar[d] & \cdots \ar[r] & p_n \ar[d] \\
			p_0' \ar[r] & p_1' \ar[r] & \cdots \ar[r] & p_n'
		\end{cd}
		whose image under $G$ equals the identity natural transformation $\sigma \xra{\Id_{\sigma}} \sigma$. 
	\end{itemize}
	Since the domain of $G$ is a poset, $\mc{C}$ is also a poset. In fact, it arises from the construction of~\ref{gal2}, because a composable sequence $p_0 \to \cdots \to p_n$ whose image under $G$ is identified with $\sigma$ can be completely characterized by the multiset of the transition values of the paths $p_0, \ldots, p_n$. Thus, Lemma~\ref{gal2} implies that $\mc{C}$ is contractible, so $(\dagger)$ is true. \hfill $\square$ 
	
	\subsection{Non-pure simplicial complexes} \label{ss-non-pure} 
	We no longer assume that the simplicial complex $K$ is of pure dimension. The material in this subsection is not used in the rest of the paper. 
	
	\subsubsection{Definition} We define the category $\ms{Comb}_K(a, b)$ which replaces the category of galleries studied in~\ref{s-gal}. The objects are sequences of simplices 
	\[
		(a = F_0, F_1, \ldots, F_n = b)
	\]
	satisfying that, for each $i$, we have $F_i \subseteq F_{i+1}$ or $F_i \supseteq F_{i+1}$. A morphism 
	\[
		(F_0, F_1, \ldots, F_n) \to (F_0', F_1', \ldots, F_m') 
	\]
	is a weakly increasing surjective map $f : \{1, \ldots, n\} \to \{1, \ldots, m\}$ satisfying that, for each $i \in \{1, \ldots, n\}$, we have $F'_{f(i)} \subset F_i$. (This is equivalent to requiring that, for each $j \in \{1, \ldots, m\}$, we have $F_j' \subseteq \bigcap_{i \in f^{-1}(j)} F_i$.) 
	
	\begin{rmk*}
		We caution the reader that this does not agree with the notion of `generalized gallery' from~\cite[Def.\ 5.1]{gallery}. That definition requires the containment relations between adjacent simplices to alternate, but we impose no such restriction. 
	\end{rmk*}
	
	\begin{thm} \label{thm-nonpure} 
		The categories $\ms{Path}_K(a, b)$ and $\ms{Comb}_K(a, b)$ are homotopy equivalent. 
	\end{thm}
	
	\begin{rmk*} Although we could not deduce this theorem from the literature, here are two earlier works which are related. 
		\begin{enumerate}[label=(\arabic*)] 
			\item In~\cite[Def.\ 2.2]{segal}, Segal defined \emph{partial monoids} and their classifying spaces. Briefly, a partial monoid is a space with a partially defined associative multiplication, and its classifying space is the universal homotopy-coherent topological monoid which admits a map from it. There is an obvious variant of this notion in which `monoid' is replaced by `semigroup.' 
			
			To relate this to our construction, view $\simp(K)$ as a \emph{discrete} partial semigroup, where the partial multiplication is given by union of simplices. The classifying space (in Segal's definition) is a simplicial set $X$; let $\mb{\Delta}_{/X}$ be its category of simplices. There is a fully faithful embedding 
			\[
				\ms{Comb}_K(a, b) \hra \mb{\Delta}_{/X}. 
			\]
			This is not an equivalence because the right hand side does not incorporate our requirement that adjacent simplices are comparable. 
			\item Let $a$ and $b$ be vertices of $K$. In~\cite[2.2]{grandis}, one can find the following model for the space of paths from $a$ to $b$. Equip $\BZ$ with the simplicial complex structure that makes it look like a line. A \emph{Moore path} is a map $p : \BZ \to K$ which is eventually constant for $t \ll 0$ with value $a$, and for $t \gg 0$ with value $b$. The desired model for the path space is the set of Moore paths equipped with the simplicial complex structure in which paths $p_0, \ldots, p_n$ form an $n$-simplex if and only if, for all $t \in \BZ$, the set $\cup_{i=0}^n \{p_i(t), p_i(t+1)\}$ is a simplex of $K$. (See also~\cite[\S 5]{paths}.) 
			
			Although this model also has a combinatorial flavor, it is essentially different from ours because it makes reference to the time coordinate. As explained in Remark~\ref{gal2}, the main difficulty in developing our path space models is to make them agnostic with respect to the time coordinate by allowing the transition points to `move.' 
			
			On the other hand, Segal's construction in (1) makes no reference to a time coordinate, but he uses partial monoids in a very different way; in the framework of~\cite{segal}, it would be unnatural to consider a discrete partial monoid. 
		\end{enumerate}
	\end{rmk*}
	
	The rest of this subsection is devoted to proving the theorem.
	
	\subsubsection{Definition} We define a quotient 
	\[
		\ms{Path}_K(a, b) \sra \ol{\ms{Path}}_K(a, b)
	\]
	of simplicial sets as follows. According to~\ref{fun-interpret}, we can interpret objects of $\ms{Path}_K(a, b)$ as functions on $\mb{t}$, so an $n$-simplex of $\ms{Path}_K(a, b)$ is a certain $n$-tuple of functions on $\mb{t}$. For any two simplices $\sigma_1, \sigma_2 : \Delta^n \to \ms{Path}_K(a, b)$, they are identified in the quotient if and only if there exists an increasing piecewise-linear homeomorphism $\mb{t} \to \mb{t}$ which turns $\sigma_1$ into $\sigma_2$. 
	
	More concretely, $\sigma_1$ and $\sigma_2$ are identified if one is obtained from the other by keeping the values of the functions the same while moving the transition points. 
 
	\begin{lem} \label{lifts} 
		The simplicial set $\ol{\ms{Path}}_K(a, b)$ satisfies the following: 
		\begin{enumerate}[label=(\roman*)]
			\item Existence of right lifts for $\Lambda^2_1 \hra \Delta^2$. 
			\item Uniqueness of right lifts for $\partial \Delta^2 \hra \Delta$. 
			\item $(n \ge 3)$ Existence and uniqueness of right lifts for $\partial \Delta^n \hra \Delta^n$, i.e.\ it is 2-cotruncated.
			\item Let $E := \Delta^2 / \Delta^{\{1, 2\}}$, where we collapse the last edge of $\Delta^2$. 
			
			Existence of right lifts for $E \sra \Delta^1$. (Uniqueness is trivial.) 
			\item Let $F := \Delta^2 / \Delta^{\{0, 2\}}$, where we collapse the middle edge of $\Delta^2$.
			
			Existence of right lifts for $F \sra \on{pt}$. (Uniqueness is trivial.) 
		\end{enumerate}
	\end{lem}
	\begin{myproof}
		Here is the key idea. A vertex of $\ol{\ms{Path}}_K(a, b)$ is specified by a collection of values, along with a `complete and satisfiable' collection of constraints (equalities and strict inequalities) relating the transition points. Here `complete' means that any pair of transition points is related by a constraint, and `satisfiable' means that the constraints do not imply a contradiction. Similarly, an edge $v_1 \to v_2$ is specified by a `complete and satisfiable' collection of constraints relating the transition points of $v_1$ to those of $v_2$. Finally, a map $\mr{sk}_1(\Delta^n) \to \ol{\ms{Path}}_K(a, b)$ extends to a (necessarily unique) map $\Delta^n \to \ol{\ms{Path}}_K(a, b)$ if and only if the corresponding constraints are simultaneously satisfiable. 
		
		With this in mind, everything is straightforward except for the existence statement in (iii). To prove this, suppose we are given a map $\partial \Delta^n \to \ol{\ms{Path}}_K(a, b)$. This determines a collection of constraints (on transition points) for each vertex and edge of $\Delta^n$, subject to the requirement that each set of constraints corresponding to a facet of $\Delta^n$ is simultaneously satisfiable. We want to prove that all the constraints are simultaneously satisfiable. 
		
		Suppose not. Then there is a cycle of constraints which gives a contradiction. This cycle of constraints corresponds to a cycle of edges in $\Delta^n$.
		
		Suppose that 
		\begin{cd}
			v_0 \ar[r, dash, "e_1"] & v_1 \ar[r, dash, "e_2"] & v_2
		\end{cd}
		is a consecutive pair of edges in this cycle. This part of the cycle implies a constraint $C$ between a transition point of $v_0$ and a transition point of $v_2$. Let $e$ be the edge of $\Delta^n$ which joins $v_0$ and $v_2$. Since $n \ge 3$, the triangle 
		\begin{cd}
			v_0 \ar[r, dash, "e_1"] \ar[rr, bend right = 30, dash, swap, "e"] & v_1 \ar[r, dash, "e_2"] & v_2
		\end{cd}
		is filled in $\partial \Delta^n$, so the corresponding constraints are satisfiable. But $v_0, e, v_2$ give a complete collection of constraints relating the transition points of $v_0$ and $v_2$. Therefore, the constraint $C$ is also part of those associated to $e$. By replacing $e_1, e_2$ by $e$, we obtain a shorter cycle which is still contradictory. However, any cycle of length $0$ is not contradictory. Therefore there cannot exist a cycle of constraints which gives a contradiction, so the desired lift along $\partial \Delta^n \hra \Delta^n$ exists. 
	\end{myproof}
	
	\begin{cor} \label{lifts2} 
		Suppose we are given maps 
		\begin{cd}
			\Delta^n \ar[r, "\sigma"] & \ol{\ms{Path}}_K(a, b) \ar[r, leftarrow, "f"] & \mc{C}
		\end{cd}
		where $\mc{C}$ is a 1-category. Then 
		\[
			\ul{\Hom}_{\Delta^n}\Big(\Delta^n, \Delta^n \underset{\ol{\ms{Path}}_K(a, b)}{\times} \mc{C}\Big)
		\]
		is also a 1-category.  
	\end{cor}
	\begin{myproof}
		From our point of view, a 1-category is a simplicial set which satisfies existence and uniqueness of right lifts for $\Lambda^m_k \hra \Delta^m$, for all $0 < k < m$, see~\ref{intro6}. Upon tracing through the definition of the internal Hom, we reduce to showing the following: 
		\begin{itemize}
			\item The map $g$ given by 
			\begin{cd}
				\Lambda^m_k \times \Delta^n \ar[r, "\pr_2"] & \Delta^n \ar[r, "\sigma"] & \ol{\ms{Path}}_K(a, b)
			\end{cd}
			admits a unique extension along $\imath : \Lambda^m_k \times \Delta^n \hra \Delta^m \times \Delta^n$. 
		\end{itemize}
		Certainly the formula $\sigma \circ \pr_2$ gives one such extension, so we only have to show that any extension must equal this one. Moreover, agreement on the vertices is automatic, and we only have to show agreement on the edges. Then Lemma~\ref{lifts}(ii, iii) implies agreement on the higher-dimensional simplices. 
		
		If $m \ge 3$, then every edge of $\Delta^m \times \Delta^n$ already lies in $\Lambda^m_k \times \Delta^n$, so we are done.  
		
		The only remaining case is $m = 2$. Every edge of $\Delta^2 \times \Delta^n$ factors through $\Delta^2 \times e$ where $e$ is an edge of $\Delta^n$. Hence, we may assume that $n = 1$. Here are the only edges in $\Delta^2 \times \Delta^1$ which do not lie in $\Lambda^2_1 \times \Delta^1$: 
		\e{
			 e_1 & : (0,0) \to (2, 0)  \\
			 e_2 & : (0, 0) \to (2, 1)  \\
			 e_3 & : (0, 1) \to (2, 1) 
		}
		We have labeled the vertices of $\Delta^2 \times \Delta^1$ by pairs consisting of a vertex of $\Delta^2$ and a vertex of $\Delta^1$, in that order. 
		
		Note that $e_1$ is part of the 2-simplex $\Delta^2 \times \{0\}$, and the map $g$ sends the other two edges to degenerate edges. Thus, Lemma~\ref{lifts}(iv) implies that $e_1$ must map to a degenerate edge as well. This agrees with our desired extension. By the same reasoning, $e_3$ must map to a degenerate edge, which agrees with our desired extension. 
		
		Lastly, in the 2-simplex with vertices $(0, 0), (0, 1), (2, 1)$, the first edge is part of $\Lambda^2 \times \Delta^1$, the middle edge is $e_2$, and the last edge is $e_3$. Since $e_3$ must map to a degenerate edge, Lemma~\ref{lifts}(iv) tells us that $e_2$ and the middle edge must map to the same edge. This agrees with our desired extension. We have shown that the extension is uniquely determined on $e_1, e_2, e_3$, so the proof is finished. 
	\end{myproof}
	
	\begin{lem} \label{lifts3} 
		The map $\ms{Path}_K(a, b) \to \ol{\ms{Path}}_K(a, b)$ is a homotopy equivalence. 
	\end{lem}
	\begin{myproof}
		In view of the fiber theorem for simplicial sets (Theorem~\ref{fiber0}), it suffices to show that, for every $\sigma : \Delta^n \to \ol{\ms{Path}}_K(a, b)$, the simplicial set 
		\[
			\mc{C} := \ul{\Hom}_{\Delta^n}\Big(\Delta^n, \Delta^n \underset{\ol{\ms{Path}}_K(a, b)}{\times} \ms{Path}_K(a, b) \Big)
		\]
		is contractible. By Corollary~\ref{lifts2}, this is a 1-category. 
		
		To proceed, let us describe this 1-category more concretely. By the proof of Lemma~\ref{lifts}, we may think of $\sigma$ as a collection of constraints (equalities and strict inequalities) relating the transition points that would appear in a composable sequence
		\begin{cd}
			v_0 \ar[r] &   v_1 \ar[r] &  \cdots \ar[r] &  v_n
		\end{cd}
		in the poset $\ms{Path}_K(a, b)$. In addition, $\sigma$ specifies the values which occur. Now, an object in $\mc{C}$ is a composable sequence in $\ms{Path}_K(a, b)$ which actually satisfies the constraints imposed by $\sigma$. A morphism in $\mc{C}$ is a commutative diagram 
		\begin{cd}
			v_0 \ar[r]  \ar[d]&  v_1 \ar[r] \ar[d] &  \cdots \ar[r] &   v_n  \ar[d]\\
			v_0' \ar[r] & v_1' \ar[r] & \cdots \ar[r] & v_n'
		\end{cd}
		in $\ms{Path}_K(a, b)$. 
		
		This description makes it clear that $\mc{C}$ arises from the construction of~\ref{gal2}: a composable sequence (as above) is sent to the tuple consisting of the transition points of $v_0, \ldots, v_n$. Now Lemma~\ref{gal2} implies that $\mc{C}$ is contractible, as desired. 		
	\end{myproof}
	
	\subsubsection{Definition} \label{path-ol-bar} 
	Let 
	\[
		\ol{\ms{Path}}_K^{\mr{reg}}(a, b) \subset \ol{\ms{Path}}_K(a, b)
	\]
	be the image of $\ms{Path}_K^{\mr{reg}}(a, b)$ under the quotient map $\ms{Path}_K^{\mr{reg}}(a, b) \sra \ol{\ms{Path}}_K(a, b)$. This is a full embedding of simplicial subsets. From the point of view of the proof of Lemma~\ref{lifts}, it consists of all vertices which correspond to a set of `constraints' which are all strict inequalities (no equalities). 
	
	\begin{lem*}
		The embedding $\ol{\ms{Path}}_K^{\mr{reg}}(a, b) \hra \ol{\ms{Path}}_K(a, b)$ is a homotopy equivalence. 
	\end{lem*}
	\begin{myproof}
		Consider the commutative diagram 
		\begin{cd}
			\ms{Path}_K^{\mr{reg}}(a, b) \ar[r, hookrightarrow] \ar[d] & \ms{Path}_K(a, b) \ar[d] \\
			\ol{\ms{Path}}_K^{\mr{reg}}(a, b) \ar[r, hookrightarrow] & \ol{\ms{Path}}_K(a, b)
		\end{cd}
		By Lemma~\ref{lem-reg}, the upper horizontal map is a homotopy equivalence. By Lemma~\ref{lifts3}, the right vertical map is a homotopy equivalence. Therefore, it suffices to show that the left vertical map is a homotopy equivalence. This can be proved in exactly the same way as Lemma~\ref{lifts3}; the reason is that Lemma~\ref{gal2} allows us to freely choose, in the definition of $\mc{P}$, whether the various inequalities are weak or strict. 
	\end{myproof}
	
	\subsubsection{} \label{lifts4} 
	We will define a map 
	\[
		G : \ms{Comb}_K(a, b) \to \ol{\ms{Path}}^{\mr{reg}}_K(a, b). 
	\]
	A simplex $\sigma : \Delta^n \to \ms{Comb}_K(a, b)$ is a composable sequence 
	\[
		v_0 \to v_1 \to \cdots \to v_n, 
	\]
	where each object $v_i = (F_{i, 1}, \ldots, F_{i, \ell_i})$ is a sequence of simplices of $K$, and each arrow is specified by a surjective map of index sets: 
	\[
		f_i : \{(i-1, 1), \ldots, (i-1, \ell_{i-1})\} \to \{(i, 1), \ldots, (i, \ell_{i})\}. 
	\]
	We have written these elements as ordered pairs to emphasize that the index sets are disjoint. Assign a nonzero binary fraction to each element of these index sets as follows: 
	\begin{itemize}
		\item First, assign nonzero binary fractions to $(1, 1), \ldots, (1, \ell_1)$ such that the sum of these fractions is 1. 
		\item For $i > 1$, we inductively assign to $(i, j)$ the sum of the fractions which are assigned to the elements of $f_{i}^{-1}((i, j))$. 
	\end{itemize}
	The fractions assigned to $(i, 1), \ldots, (i, \ell_i)$ will sum to 1. 
	
	Next, convert each $v_i$ into an object $p_i \in \ms{Path}_K(a, b)$ as follows. According to~\ref{fun-interpret}, it suffices to give sequences of transition points $(t_0, \ldots, t_{\ell_i})$ and values $(\sigma_1, \ldots, \sigma_{\ell_i})$. We determine the transition points by requiring that $t_j - t_{j-1}$ equals the fraction assigned to $(i, j)$. We determine the values by requiring that $\sigma_j = F_{i, j}$. Lastly, if $\sigma_j = \sigma_{j+1}$, then the transition point $t_j$ is redundant, so we delete it. 
	
	The assignment $i \mapsto p_i$ yields a simplex $\eta : \Delta^n \to \ms{Path}_K(a, b)$ which depends on the choices of binary fractions made above. But the image $\ol{\eta} : \Delta^n \to \ol{\ms{Path}}^{\mr{reg}}_K(a, b)$ does not depend on these choices. We characterize $G$ by requiring that $G \circ \sigma = \ol{\eta}$. 
	
	\begin{lem} \label{lifts5} 
		The map $G : \ms{Comb}_K(a, b) \to \ol{\ms{Path}}^{\mr{reg}}_K(a, b)$ is a homotopy equivalence. 
	\end{lem}
	\begin{myproof}
		As in Lemma~\ref{lifts3}, it suffices to show that, for any simplex $\sigma : \Delta^n \to \ol{\ms{Path}}^{\mr{reg}}_K(a, b)$, the 1-category 
		\[
		\mc{C}' := \ul{\Hom}_{\Delta^n}\Big(\Delta^n, \Delta^n \underset{\ol{\ms{Path}}^{\mr{reg}}_K(a, b)}{\times} \ms{Comb}_K(a, b) \Big)
		\]
		is contractible. Again, the simplex $\sigma$ specifies a collection of constraints relating the transition points that would appear in a composable sequence 
		\[
			v_0 \to v_1 \to \cdots \to v_n
		\]	
		in $\ms{Path}_K^{\mr{reg}}(a, b)$. To show that $\mc{C}'$ is contractible, we will construct a terminal object.

		Choose an arbitrary lift $\wt{\sigma} : \Delta^n \to \ms{Path}_K^{\mr{reg}}(a, b)$ of $\sigma$, i.e.\ a map which sends $i \in \{0, \ldots, n\}$ to a function $f_i : \bm{t} \to \simp(K)$. Denote the sequence of transition points of $f_i$ by $(t_{i, 0}, \ldots, t_{i, \ell_1})$. Then $\sigma$ is obtained from $\wt{\sigma}$ by retaining only the sequences of values of the $f_i$ as well as the equalities and inequalities which are satisfied by the transition points $t_{i, j}$. 
		
		For each $i$, we call $t \in \bm{t}$ a \emph{level-$i$ break point} if $t = t_{j, k}$ for some $j \ge i$ and any $k$. Let $(0 = b_{i, 0}, b_{i, 1}, \ldots, b_{1, m_i} = 1)$ be the set of all level-$i$ break points sorted in increasing order. Define the sequence of simplices $(F_{i, 1}, F_{i, 2}, \ldots, F_{i, m_i}) \in \ms{Comb}_K(a, b)$ by 
		\[
			F_{i, j} := (\text{the value of $f_i$ on the open interval $(b_{j-1}, b_j)$}). 
		\]
		For any $i < i'$, there is a morphism 
		\[
			(F_{i, 1}, F_{i,2}, \ldots, F_{i, m_i}) \to (F_{i', 1}, F_{i',2}, \ldots, F_{i', m_{i'}})
		\]
		specified by a map of index sets which is given as follows: if $(b_{i, j-1}, b_{i, j}) \subseteq (b_{i', j'-1}, b_{i', j'})$ for some $j$ and $j'$, then we require $j \mapsto j'$. The assignment 
		\[
			i \mapsto (F_{i, 1}, F_{i,2}, \ldots, F_{i, m_i})
		\]
		together with the above morphisms determines an $n$-simplex $\sigma_{\mr{min}} : \Delta^n \to \ms{Comb}_K(a, b)$ which specifies the desired object of $\mc{C}'$. This object does not depend on the choice of lift $\wt{\sigma}$. 
		
		Next, we prove that $\sigma_{\mr{min}}$ is terminal. Any other object of $\mc{C}'$ is represented by a simplex $\sigma' : \Delta^n \to \ms{Comb}_K(a, b)$ and is therefore specified by an assignment
		\[
			i \mapsto (F'_{i, 1}, F'_{i, 2}, \ldots, F'_{i, m'_i}) \in \ms{Comb}_K(a, b)
		\]
		together with surjective maps between the index sets. We have to show that there is a unique map $\sigma' \to \sigma_{\mr{min}}$ which lies over the identity map $\sigma \to \sigma$ in $\ol{\ms{Path}}^{\mr{reg}}_K(a, b)$. We do this inductively, starting from the $n$-th vertex. 
		
		\noindent \emph{The base case}. The only map
		\[
			\sigma'(i) = (F'_{n, 1}, F'_{n, 2}, \ldots, F'_{n, m'_n}) \to (F_{n, 1}, F_{n,2}, \ldots, F_{n, m_n})
		\]
		in $\ms{Comb}_K(a, b)$ which lies over the identity map $v_n \to v_n$ in $\ol{\ms{Path}}_K^{\mr{reg}}(a, b)$ is the one given as follows. Partition $\sigma'(i)$ into maximal constant substrings. If 
		\[
		(F'_{n, j}, F'_{n, j+1}, \ldots, F'_{n, j'}) \text{ satisfying } F'_{n, j} = F'_{n, j+1} = \cdots = F'_{n, j'}
		\]
		is the $k$-th maximal constant substring, then the indices $\{(n, j),(n, j+1), \ldots, (n, j')\}$ are sent to the index $(n, k)$ corresponding to $F_{n, k}$. 
		
		The uniqueness follows from the fact that no two consecutive $F_{n, k}$'s are equal. The latter is true because the level-$n$ break points are just the transition points for some function satisfying the constraints for $v_n$. 
		
		\noindent \emph{The inductive step}. Now assume that the map $\sigma' \to \sigma_{\mr{min}}$ is uniquely determined at levels $i+1, i+2, \ldots, n$. To show uniqueness at level $i$, we have to show that there is a unique dashed arrow which makes this solid diagram in $\ms{Comb}_K(a, b)$ commute: 
		\begin{cd}[column sep = 1in] 
			(F'_{i, 1}, F'_{i, 2}, \ldots, F'_{i, m'_i}) \ar[r, dashed] \ar[d, swap, "\sigma'(i \to i+1)"] & (F_{i, 1}, F_{i, 2}, \ldots, F_{i, m_i}) \ar[d, "\sigma_{\mr{min}}(i \to i+1)"] \\
			(F'_{i+1, 1}, F'_{i+1, 2}, \ldots, F'_{i+1, m'_{i+1}}) \ar[r, "(\sigma' \to \sigma_{\mr{min}})(i+1)"] & (F_{{i+1}, 1}, F_{{i+1}, 2}, \ldots, F_{i+1, m_{i+1}})
		\end{cd}
		Fix an index $(i+1, j)$ in the bottom-right. We have to determine the map on fibers 
		\begin{cd}[column sep = 0.5in]
			\big((\sigma' \to \sigma_{\mr{min}})(i+1) \circ \sigma'(i \to i+1)\big)^{-1}((i+1, j)) \ar[rd, bend right = 10, swap, "(\sigma' \to \sigma_{\mr{min}})(i+1) \,\circ\, \sigma'(i \to i+1)"] \ar[r, dashed] & (\sigma_{\min}(i\to i+1))^{-1}((i+1, j')) \ar[d, "\sigma_{\mr{min}}(i \to i+1)"] \\
			& \{(i+1, j')\}
		\end{cd}
		
		This dashed map is uniquely determined for the same reason as before. Indeed, by the construction of $\sigma_{\mr{min}}$, no two consecutive $F_{i, \bullet}$ among those indexed by $(\sigma_{\min}(i\to i+1))^{-1}((i+1, j'))$ can be equal. (The level-$i$ break points which lie strictly between two consecutive level-$(i+1)$ break points must be transition points for some function satisfying the constraints for $v_i$. In particular, since these are transition points, the value of this function must change at those points.) Then the requirement that $(\sigma' \to \sigma_{\mr{min}})(i)$ lies over the identity map $v_{i} \to v_{i}$ now implies that, working in the substring of $\sigma'(i)$ indexed by 
		\[
			\big((\sigma' \to \sigma_{\mr{min}})(i+1) \circ \sigma'(i \to i+1)\big)^{-1}((i+1, j)), 
		\]
		the indices in the $k$-th maximal constant substring must map to the $k$-th element of $(\sigma_{\min}(i\to i+1))^{-1}((i+1, j'))$. 
	\end{myproof}
	
	\subsubsection{Proof of Theorem~\ref{thm-nonpure}} \label{lifts6} 
	By Lemma~\ref{lifts3}, $\ms{Path}_K(a, b)$ is homotopy equivalent to $\ol{\ms{Path}}_K(a, b)$. By Lemma~\ref{path-ol-bar}, this is homotopy equivalent to $\ol{\ms{Path}}^{\mr{reg}}_K(a, b)$. Lastly, by Lemma~\ref{lifts5}, this is homotopy equivalent to $\ms{Comb}_K(a, b)$, and the result follows. \hfill $\square$

	\section{Generalized hyperdescent for \texorpdfstring{$\D$}{D}-modules} \label{s-hyper}
	
	We will prove the generalized hyperdescent statement (Theorem~\ref{simp0}) which we need for our main theorem. As explained in~\ref{intro4}, hyperdescent for $\D$-modules (Theorem~\ref{hyper1}) is known to the experts, and our main task is to modify the statement to allow `generalized hypercovers' which are indexed by an arbitrary small $\infty$-category instead of $\mb{\Delta}^{\op}$. 
	
	We work over an arbitrary field $k$ of characteristic zero and adopt the conventions in~\ref{intro6}. In particular,  `finite type' includes a hypothesis of separability. 
	
	\subsection{Hyperdescent} \label{ss-hyper} 
	The material in this subsection was taught to us by Sam Raskin.

	\subsubsection{Simplicial indschemes} \label{move}
	
	Let $X_\bullet \to Y$ be an augmented simplicial indscheme. We think of $X_\bullet$ as a simplicial object in the category $\ms{IndSch}_{/Y}$, i.e.\ as a functor 
	\[
		X_\bullet : \mb{\Delta}^{\op} \to \ms{IndSch}_{/Y}
	\]
	Let $\fsimpset$ be the category of contractible simplicial sets with finitely many nondegenerate simplices. Since the category $\ms{IndSch}_{/Y}$ admits finite limits, we can form the following right Kan extension: 
	\[
		\ul{\Hom}(-, X_\bullet) := \mr{RKE}_{\mb{\Delta}^{\op} \hra \fsimpset^{\op}} \big( X_\bullet \big) : \fsimpset^{\op} \to \ms{IndSch}_{/Y}
	\]
	Abstractly, this construction shows that simplicial indschemes are cotensored over $\fsimpset$. In addition, for any $K \in \fsimpset$ and $S \in \ms{IndSch}_{/Y}$, we have 
	\begin{equation} \tag{$\spadesuit$} \label{swap} 
		\Hom_{\sset}(K, \Hom_{\ms{IndSch}_{/Y}}(S, X_\bullet)) \simeq \Hom_{\ms{IndSch}_{/Y}}(S, \ul{\Hom}(K, X_\bullet)). 
	\end{equation}
	
	If the maps $X_n \to Y$ are ind-proper, and the indschemes $X_n$ and $Y$ are separated, then $\ul{\Hom}(K, X_\bullet) \to Y$ is ind-proper for any $K \in \fsimpset$, and the functor $\ul{\Hom}(-, X_\bullet)$ sends each morphism in $\fsimpset$ to an ind-proper map. 
	
	As a special case, the indscheme of $n$-simplices of the $(n-1)$-dimensional coskeleton of $X_\bullet \to Y$ can be realized as follows: 
	\[
		\ul{\Hom}(\partial \Delta^n, X_\bullet) \in \ms{IndSch}_{/Y}
	\]
	For the purposes of this paper, one may take this statement as a definition. 
	
	\subsubsectiona \label{hyper0} 
	An ind-proper hypercover is an augmented simplicial indscheme $X_\bullet \to Y$ with ind-proper transition maps satisfying the following equivalent `hypercover' properties: 
	\begin{enumerate}[label=(H\arabic*)] 
		\item For each $n \ge 0$, the map 
		\[
			X_n \to \ul{\Hom}(\partial \Delta^n, X_\bullet) \in \ms{IndSch}_{/Y}
		\]
		is surjective on geometric points.\footnote{As noted in~\ref{move}, if $Y$ and $X_n$ are separable, and the maps $X_n \to Y$ are ind-proper, this map is automatically ind-proper.} Note that $\partial \Delta^0 := \emptyset$, so $\ul{\Hom}(\partial \Delta^0, X_\bullet) = Y$.
		\item For each geometric point $\eta : \Spec F \to Y$, the simplicial set 
		\[
			\Delta^n \mapsto \Hom_Y(\Spec F, X_n)
		\]
		is an acyclic Kan complex.
	\end{enumerate}
	The equivalence between these properties is a consequence of the fact that a simplicial set is an acyclic Kan complex if and only if it satisfies the right lifting property with respect to the maps $\partial \Delta^n \hra \Delta^n$ for all $n \ge 0$. 
	
	\begin{rmk*}
		Even if $k = \BC$, it is not enough to consider surjectivity at the level of $\BC$-points. For example, consider the ind-proper map of indschemes 
		\[
			\coprod_{z \in \BC} \Spec \BC \to \BA^1_{\BC}, 
		\]
		where the factor of the disjoint union indexed by $z$ maps to the closed point $z \in \BA^1_{\BC}$. This map is surjective at the level of $\BC$-points, but not at the level of geometric points. Moreover, the 
		\v{C}ech nerve of this map becomes constant upon discarding infinitesimal structure, so $\D$-modules do not satisfy descent for this cover. 
	\end{rmk*}
	
	\begin{thm} \label{hyper1} 
		Let $X_\bullet \to Y$ be an ind-proper hypercover of ind-finite type. The pullback functor 
		$\D(Y) \to \D(\colim X_\bullet)$
		is an equivalence.\footnote{Here and in what follows, $\colim X_\bullet$ is evaluated in the $\infty$-category of prestacks. Thus, we have $\D(\colim X_\bullet) \simeq \lim_n \D(X_n)$, where the limit diagram uses $!$-pullbacks.} 
	\end{thm}
	
	The rest of this subsection is devoted to proving the theorem. The crux is to prove Lemma~\ref{hyper2}, which says that the pullback functor is fully faithful. Then, we use a trick to deduce essential surjectivity from fully faithfulness; the basic idea is explained in~\ref{surj1}. 
	
	\begin{lem} \label{bound1} 
		Let $X$ be a finite type scheme, and let
		\[
			U \overset{j}{\hra} X \overset{i}{\hookleftarrow} Z
		\]
		be the embedding of an open subscheme and its closed complement. 
		\begin{enumerate}[label=(\textrm{\roman*})] 
			\item If $j$ is affine, then $j_*$ is $t$-exact, and $i^!$ has cohomological amplitude bounded by $[0, 1]$. 
			\item Assume that each $p \in Z$ admits a neighborhood on which $Z \hra X$ is set-theoretically cut out by $d$ equations. Then $j_*$ has cohomological amplitude bounded by $[0, \dim X - 1]$, and $i^!$ has cohomological amplitude bounded by $[0, \dim X]$. 
		\end{enumerate} 
	\end{lem}
	\begin{myproof}
		We only prove (ii), since the proof of (i) is entirely similar. Choose a point $p \in Z$. By hypothesis, there is an open neighborhood $p \in V$ and functions $f_1, \ldots, f_d \in \oh_V$ such that 
		\[
			Z \cap V = \mc{V}(f_1, \ldots, f_{\dim X})
		\]
		as sets. Upon shrinking $V$, we can find a closed embedding 
		\[
			\imath : V \hra S
		\]
		where $S$ is smooth, such that $f_1, \ldots, f_{\dim X}$ extend to regular functions $\wt{f}_1, \ldots, \wt{f}_{\dim X}$ on $S$. Define 
		\[
			\wt{Z} := \mc{V}(\wt{f}_1, \ldots, \wt{f}_{\dim X}) \subset S, 
		\]
		so we have a cartesian diagram 
		\begin{cd}
			U \cap V \ar[r, "j"] \ar[d, "i"] & V \ar[d, "\imath"] \\
			S \setminus \wt{Z} \ar[r, "\jmath"] & S
		\end{cd}
		This gives rise to a commutative diagram of $\infty$-categories: 
		\begin{cd}
			\D(U \cap V) \ar[r, "j_*"] \ar[d, "i_*"] & \D(V) \ar[d, "\imath_*"] \\
			\D(S \setminus \wt{Z}) \ar[r, "\jmath_*"] \ar[d, "\mathbf{oblv}_{S \setminus \wt{Z}}"]  & \D(S) \ar[d, "\mathbf{oblv}_S"]  \\
			\ms{IndCoh}(S \setminus \wt{Z}) \ar[r, "\jmath_*"] & \ms{IndCoh}(S)
		\end{cd}
		Recall that $\imath_*$ is $t$-exact by~\cite[Prop.\ 4.2.5]{crystals}, and $\textbf{oblv}_S$ is $t$-exact by~\cite[Prop.\ 4.2.11(a)]{crystals} since $S$ is smooth. Since these functors are also conservative, to show the desired bound on the cohomological amplitude of $j_*$, it suffices to show a similar bound on $\mathbf{oblv}_S\circ \imath_* \circ j_*$. 
		
		Let us rewrite this composition as $\jmath_* \circ \mb{oblv}_{S \setminus \wt Z} \circ i_*$. The last two functors are $t$-exact by the same reasoning as before. Furthermore, $\jmath_*$ has cohomological amplitude bounded by $[0, d]$ because it be computed by a \v{C}ech complex living in degrees $[0, \dim X -1]$, coming from the Zariski open cover of $S \setminus \wt Z$ given by the equations $\wt{f}_1, \ldots, \wt{f}_{\dim X}$. 	
		
		Finally, the statement about $i^!$ follows from the Cousin exact triangle 
		\[
			i_*i^! \to \Id_{\D(X)} \to j_*j^! \to 
		\]
		Indeed, since $i_*$ is $t$-exact and conservative, $i^!$ has cohomological amplitude bounded by $[0, \dim X]$ if and only if $i_*i^!$ does. Since $j^!$ is $t$-exact, $j_*j^!$ has cohomological amplitude bounded by $[0, \dim X - 1]$, and the result follows from the long exact sequence for cohomology.
	\end{myproof} 
	
	\begin{lem} \label{bound2} 
		Let $f : X \to Y$ be an ind-proper map, and assume that $Y$ is a finite type scheme. Let $(f^!)^R$ be the non-continuous right adjoint of $f^!$. Then the functor
		\[
		(f^!)^R \circ f^! [-\dim Y] : \D(Y) \to \D(Y)
		\]
		is left $t$-exact. 
	\end{lem}
	\begin{myproof}
		Choose an ind-scheme presentation $X = \colim_a X_a$ such that the resulting maps $f_a : X_a \to Y$ are proper. We have 
		\[
		(f^!)^R f^! \simeq \lim_a (f_a^!)^R f_a^!. 
		\]
		Since filtered limits are left $t$-exact, it suffices to show that each $(f_a^!)^R f_a^! [-d_Y]$ is left $t$-exact. In other words, we may assume that $f$ is proper, so $X$ is also a finite type scheme. 
		
		Assume without loss of generality that $Y$ is affine. We inductively construct compatible stratifications $X = \bigcup_b X_b$ and 
		\[
			Y = Y_0 \cup Y_1 \cup \cdots \cup Y_{\dim Y} 
		\]
		into locally closed subschemes such that the following statements are true: 
		\begin{enumerate}[label=(\roman*)]
			\item Each $Y_i$ is affine. For all $i$, we have $\dim Y_i = i$. For all $i < j$, the closed embedding $\ol{Y}_i \hra \ol{Y}_j$ is set-theoretically cut out by $j-i$ equations. 
			\item For each $b$, restriction of $f$ to $X_b$ is smooth and surjective onto the corresponding stratum of $Y$, and the embedding $X_b \hra f^{-1}(f(X_b))$ is affine. 
		\end{enumerate}
		To construct this stratification, start by finding a dense open subscheme $U \subset Y$ such that the map $f^{-1}(U) \to U$ admits such a stratification for which $U$ is the only stratum of the target. This can be done as long as $U$ is sufficiently small. By shrinking $U$ even more, we may assume that it is a distinguished affine open subscheme, meaning that it is cut out by one equation, say $g \in \oh_Y$. Since $\dim V(g) = \dim Y - 1$, we may assume by induction that the map $f^{-1}(V(g)) \to V(g)$ admits such a stratification. Finally, taking the union of these two stratifications yields the desired stratification for $f$. The key point is that property (i) remains true because, if $\ol{Y}_i$ is cut out from $\ol{Y}_{\dim Y - 1}$ by $\dim Y - 1 - i$ equations, then it is cut out from $\ol{Y}_{\dim Y}$ by $\dim Y - i$ equations, since we can add in the equation $g$. 
		
		For each $b$, let $\imath_b : X_b \hra X$ denote the embedding of the corresponding stratum. It suffices to show that, for $\mc{N} \in \D(Y)^{\le 0}$ and $\mc{M} \in \D(Y)^{\ge 0}$, we have 
		\[
			\Hom_{\D(X)}(f^! \mc{N}, f^!\mc{M}) \in \ms{Vect}^{\ge -\dim Y}. 
		\]
		By the Cousin exact triangle, this Hom complex is an iterated extension of complexes
		\[
			\Hom_{\D(X)}((\imath_{c})_* \imath_{c}^! f^! \mc{N}, (\imath_{b})_* \imath_{b}^! f^!\mc{M})
		\]
		for various indices $b, c$. Our goal is to show that each of these complexes lies in $\ms{Vect}^{\ge -\dim Y}$. An easy base-change argument shows that this complex is zero unless $X_{b} \subset \ol{X}_{c}$, so we assume that this is the case. 
		
		Let $b' \le c'$ be the indices for which $f(X_b) = Y_{b'}$ and $f(X_c) = Y_{c'}$. Let $X_{[b, c]}$ obtained by taking the union of all strata $X_{d}$ such that $\ol{X}_{b} \subset X_{d} \subset \ol{X}_{c}$. Define $Y_{[{b'}, {c'}]}$ similarly. Consider the following morphisms: 
		\begin{cd}
			& X_b \ar[ld, swap, "i_b"] \ar[rd, "\imath_b"] \\
			X_{[b, c]} \ar[rr, swap, "\dot{\epsilon}"] \ar[d, "\dot{f}"] & & X \ar[d, "f"] \\
			Y_{[{b'}, {c'}]} \ar[rr, swap, "\epsilon"] & & Y
		\end{cd}
		We have $(\imath_b)_*\imath_b^!f^! \simeq \dot{\epsilon}_* (i_b)_* i_b^! \dot{f}^! \epsilon^!$, and a similar relation with $c$ in place of $b$. By Kashiwara's lemma~\cite[Prop.\ 2.5.6]{crystals}, the functor $\dot{\epsilon}_*$ is fully faithful since $\dot{\epsilon}$ is a locally closed embedding. It follows that 
		\[
			\Hom_{\D(X)}((\imath_c)_* \imath_c^! f^! \mc{N}, (\imath_{b})_* \imath_{b}^! f^!\mc{M}) \simeq \Hom_{\D(X_{[b, c]})}((i_c)_* i_c^! \dot{f}^! (\epsilon^! \mc{N}), (i_{b})_* i_{b}^! \dot{f}^! (\epsilon^! \mc{M})). 
		\]
		By (i), the locally closed embedding $Y_{[{b'}, {c'}]} \hra Y$ is set-theoretically cut out by $\dim Y - c'$ equations. Lemma~\ref{bound1}(ii) implies that $\epsilon^! \mc{N} \in \D(Y_{[{b'}, {c'}]})^{\le \dim Y - c'}$, and we clearly also have $\epsilon^! \mc{M} \in \D(Y_{[{b'}, {c'}]})^{\ge 0}$. Therefore, it suffices to show that, for all $\dot{\mc{N}} \in \D(Y_{[{b'}, {c'}]})^{\le 0}$ and $\dot{\mc{M}} \in \D(Y_{[{b'}, {c'}]})^{\ge 0}$, we have 
		\[
			\Hom_{\D(X_{[b, c]})}((i_c)_* i_c^! \dot{f}^! \dot{\mc{N}}, (i_{b})_* i_{b}^! \dot{f}^! \dot{\mc{M}}) \in \ms{Vect}^{\ge -c'}. 
		\]
		
		We will prove this by bounding the cohomological amplitudes of the two objects in the Hom complex. Let $d_1$ be the relative dimension of $X_b \to Y_{b'}$, and let $d_2$ be the relative dimension of $X_{c} \to Y_{c'}$. Since $X_b \subset \ol{X}_c$, we have 
		\[
			d_1 \le d_2 + \mr{codim}(Y_{b'}, \ol{Y}_{c'}) \le d_2 + c'. 
		\] 		
		\begin{itemize}
			\item We claim that $(i_{c})_* i_{c}^! \dot{f}^!\dot{\mc{N}} \in \D(X)^{\le -d_2}$. There is a commutative diagram 
			\begin{cd}
				X_{c} \ar[r, "i_{c}"] \ar[d, "\ddot{f}"] & X_{[b, c]} \ar[d, "\dot{f}"] \\
				Y_{c'} \ar[r, "\gamma"] & Y_{[{b'}, {c'}]}
			\end{cd}
			which implies that $i_{c}^! \dot{f}^! \simeq \ddot{f}^! \gamma^!$. Since $\gamma$ is an open embedding, $\gamma^!$ is $t$-exact. By (ii), $\ddot{f}$ is smooth of relative dimension $d_2$, so it is $t$-exact up to a shift of $-d_2$. Finally, to see that $(i_{c})_*$ is $t$-exact, split into two cases: 
			\begin{itemize}
				\item Assume that $X_{b} \subset f^{-1}(f(X_{c}))$, so $b' = c'$. Then (ii) implies that $i_{c}$ is an affine open embedding, so $(i_{c})_*$ is $t$-exact by Lemma~\ref{bound1}(i). 
				\item Assume that $X_{b} \cap f^{-1}(f(X_{c})) = \emptyset$, so $b' \neq c'$. We factor $i_{c}$ as follows: 
				\[
					X_{c} \hra \dot{f}^{-1}(Y_{c'}) \hra X_{[b, c]}
				\]
				By (ii), the first map is an affine open embedding, since $Y_{c'} = f(X_{c})$. By (i), $Y_{c'} \hra Y_{[b', c']}$ is an affine open embedding, so the same is true of $\dot{f}^{-1}(Y_{c'}) \hra \dot{f}^{-1}(Y_{[{b'}, {c'}]}) = X_{[b, c]}$. Therefore, $i_{c}$ is an affine open embedding, so Lemma~\ref{bound1} implies that $(i_{c})_*$ is $t$-exact. 
			\end{itemize} 
			\item We claim that $(i_b)_* i_b^! \dot{f}^! \dot{\mc{M}} \in \D(X)^{\ge -d_1}$. This time, we consider the commutative diagram 
			\begin{cd}
				X_{b} \ar[r, "i_{b}"] \ar[d, "\dddot{f}"] & X_{[b, c]} \ar[d, "\dot{f}"] \\
				Y_{b'} \ar[r, "\dot{\gamma}"] & Y_{[{b'}, {c'}]}
			\end{cd}
			which implies that $i_{b}^! \dot{f}^! \simeq \dddot{f}^! \dot{\gamma}^!$. Since $\dot{\gamma}$ is a closed embedding, $\dot{\gamma}^!$ is left $t$-exact. By (ii), $\dddot{f}$ is smooth of relative dimension $d_1$, so it is $t$-exact up to a shift of $-d_1$. Finally, $i_{b}$ is a closed embedding, so $(i_b)_*$ is $t$-exact. 
		\end{itemize}
		This implies that the Hom complex lies in $\ms{Vect}^{\ge -d_1 + d_2}$. As noted above, this is contained in $\ms{Vect}^{\ge -c'}$, as desired. 
	\end{myproof}
	
	\begin{lem} \label{hyper2} 
		Let $X_\bullet \to Y$ be an ind-proper hypercover of ind-finite type. The pullback functor 
		$
			\D(Y) \to \D(\colim X_\bullet)
		$
		is fully faithful. 
	\end{lem}
	\begin{myproof}
		First, assume that $Y$ is finite type, and fix two objects $\mc{F}, \mc{G} \in \D(Y)$. For each $n$, let $f_n : X_n \to Y$ be the map in the hypercover. We want to show that the map 
		\[
			\Hom_{\D(Y)}(\mc{F}, \mc{G}) \to \lim_{\Delta^n \in \mb{\Delta}} \Hom_{\D(X_n)}(f_n^! \mc{F}, f_n^! \mc{G})
		\]
		is an equivalence. 
		
		By~\cite[4.5]{crystals}, the (right) $t$-structure on $\D(Y)$ is left and right complete. It follows that 
		\e{
			\mc{F} & \simeq \colim_m \tau^{\le m} \mc{F} \\
			\mc{G} & \simeq \lim_m \tau^{\ge m} \mc{G}, \\
			\Hom_{\D(Y)}(\mc{F}, \mc{G}) &\simeq \lim_{m_1, m_2} \Hom_{\D(Y)}(\tau^{\le m_1} \mc{F}, \tau^{\ge m_2} \mc{G}) \\
			\Hom_{\D(X_n)}(f_n^!\mc{F}, f_n^!\mc{G}) &\simeq \lim_{m_1, m_2} \Hom_{\D(X_n)}(f_n^! \tau^{\le m_1} \mc{F},  f_n^! \tau^{\ge m_2} \mc{G}). 
		} 
		In the last line, we have used that $f_n^!$ is a continuous right adjoint (since $f_n$ is ind-proper), so it commutes with limits and colimits. In other words, it suffices to prove the original statement under the assumption that $\mc{F} \in \D(Y)^{\le m_1}$ and $\mc{G} \in \D(Y)^{\ge m_2}$ for some $m_1, m_2$.	
		
		Next, consider the category  
		\[
			\mc{C} := \ms{Sch}^{\mr{aff}}_{/Y}, 
		\]
		so that 
		\[
			\Fun(\mc{C}^{\op}, \ms{Spaces}) \simeq \ms{PreStk}_{/Y}. 
		\]
		In this prestack $\infty$-category, we consider two objects: 
		\begin{itemize}
			\item Define $\mc{X} := \colim X_\bullet$ where the colimit is evaluated in this prestack $\infty$-category. 
			\item Define a prestack $F$ by the formula 
			\[
				F(Y' \xra{p} Y) := \Hom_{\D(Y')}(p^!\mc{F}, p^! \mc{G}).
			\]
			By Lemma~\ref{bound2} and our assumption that $\mc{F} \in \D(Y)^{\le m_1}$ and $\mc{G} \in \D(Y)^{\ge m_2}$, we find that $F$ has no homotopy in (cohomological) degree $< m_2 - m_1 - \dim Y$. 
		\end{itemize}
		The desired fully faithfulness statement is that the functor $\Hom_{\ms{PreStk}_{/Y}}(-, F)$ sends the map
		$
			\mc{X} \to Y
		$
		to an equivalence. 
		
		We will deduce this statement from~\cite[Cor.\ A.9]{dugger}, which articulates the general idea that descent implies hyperdescent under a coconnectivity hypothesis. The proof of~\cite[Prop.\ 3.2.2]{crystals} implies that $\D$-modules satisfy descent with respect to (\v{C}ech nerves of) ind-proper maps which are surjective on geometric points. (The main ingredient is the analogous statement for ind-coherent sheaves, which is stated in \cite[Lem.\ 2.10.3]{indschemes}.) The `fully faithfulness' part of this descent statement implies an analogous descent statement for the prestack $F$. Now~\cite[Cor.\ A.9]{dugger} together with the bound on homotopy of $F$ implies that $F$ satisfies descent for all hypercovers, which includes the map $\mc{X} \to Y$ from above. This concludes the proof when $Y$ is finite type. 
		
		For the general case, let $Y = \colim_\alpha Y_\alpha$ be a presentation of $Y$ as an ind-scheme of ind-finite type. Applying what we have already proved to 
		\[
			Y_\alpha \underset{Y}{\times} X_\bullet \to Y_\alpha
		\]
		and taking a filtered limit with respect to $\alpha$ yields the desired statement. 
	\end{myproof}
	
	\subsubsection{Motivation} \label{surj1} 
	Let $p : \colim X_\bullet \to Y$ be the hypercover, which includes the data of maps $p_n : X_n \to Y$. Lemma~\ref{hyper2} says that $p^! : \D(Y) \to \D(\colim X_\bullet)$ is fully faithful, so it remains to prove essential surjectivity. We will accomplish this by applying Lemma~\ref{hyper2} to various hypercovers which are constructed from the original one. To illustrate the method, we sketch a proof that $\mc{F}_0$ lies in the essential image of $p_0^! : \D(Y) \to \D(X_0)$. 
	
	We need to show that $\mc{F}_0$ descends along $p_0 : X_0 \to Y$. To construct a descent datum, the first step is to consider the diagram 
	\begin{cd}
		X_0 \underset{Y}{\times} X_0 \ar[r, "\pr_2"] \ar[d, "\pr_1"] & X_0 \ar[d, "p_0"] \\
		X_0 \ar[r, "p_0"] & Y
	\end{cd}
	and construct an isomorphism $\pr_1^! \mc{F}_0 \simeq \pr_2^! \mc{F}_0$. For this, we expand the diagram as follows: 
	\begin{cd}
		X_1 \ar[rd] \ar[rrd] \ar[rdd]  \\
		& X_0 \underset{Y}{\times} X_0 \ar[r, swap, "\pr_2"] \ar[d, "\pr_1"] & X_0 \ar[d, "p_0"] \\
		& X_0 \ar[r, "p_0"] & Y
	\end{cd}
	Since $\mc{F}_0$ comes from an object $\mc{F}_\bullet$, the two pullbacks of $\mc{F}_0$ to $X_1$ are canonically identified. In order to descend this identification along $X_1 \to X_0 \underset{Y}{\times} X_0$, we want to realize $X_1$ as the zeroth ind-scheme in a hypercover: 
	\begin{cd}
		\cdots \ar[r, shift left = 3] \ar[r, shift left = 1] \ar[r, shift right = 1] \ar[r, shift right = 3] & (?) \ar[r, shift left = 2] \ar[r] \ar[r, shift right = 2] & (?) \ar[r, shift left = 1] \ar[r, shift right = 1] & X_1 \ar[r] & X_0 \underset{Y}{\times} X_0
	\end{cd}
	To guess the higher terms in the hypercover, we observe that, for any geometric point $q : \Spec F \to X_0 \underset{Y}{\times} X_0$, which corresponds to a pair of points $q_1, q_2 \in \Hom_{\ms{IndSch}_{/Y}}(\Spec F, X_0)$, the set $\Hom_{/(X_0 \times_Y X_0)}(\Spec F, X_1)$ corresponds to 0-simplices in the simplicial mapping space from $q_1$ to $q_2$ in the simplicial set $\Hom_{\ms{IndSch}_{/Y}}(\Spec F, X_\bullet)$. Inspired by this observation, we take the $n$-th term in our hypercover to be an ind-scheme corresponding to $n$-simplices in the (left) simplicial mapping space from $q_1$ to $q_2$ in the Kan complex $\Hom_{\ms{IndSch}_{/Y}}(\Spec F, X_\bullet)$, as defined in~\cite[1.2.2]{htt}. These ind-schemes can be constructed as fibered products using the maps in $X_\bullet \to Y$, and one can show directly that $\mc{F}_\bullet$ canonically determines (via $!$-pullback) a compatible family of $\D$-modules on them. In particular, there is a canonical isomorphism between the two pullbacks of $\mc{F}_0$ to this hypercover. By Lemma~\ref{hyper2}, this descends to an isomorphism on $X_0 \underset{Y}{\times} X_0$, as desired. 
	
	For the general case, replacing $X_0$ by $X_n$, the first step is to generalize the `simplicial mapping space' construction to apply to pairs of simplices rather than pairs of points. The actual proof will proceed not by constructing descent data along each $p_n : X_n \to Y$ but by directly showing that $\mc{F}_\bullet$ is isomorphic to $p^!p_*(\mc{F}_\bullet)$. 
		
	\subsubsectiona \label{surj2} 
	Define a functor 
	\[
		\mc{M} : \mb{\Delta} \times \sset \times \mb{\Delta} \to \sset
	\]
	via the formula 
	\[
		\mc{M}(\Delta^a, K, \Delta^c) := \colim \left( \begin{tikzcd}
		K \times \Delta^c \ar[r,hookrightarrow] \ar[d, twoheadrightarrow, "\pr_2"] & \Delta^a \star (K \times \Delta^c) \\
		\Delta^c
		\end{tikzcd} \right)
	\]
	As a special case, note that $\mc{M}(\Delta^a, \emptyset, \Delta^c) \simeq \Delta^a \sqcup \Delta^c$. 
	
	\begin{lem*}
		For fixed $a, c \ge 0$, the functor 
		\[
			\mc{M}(\Delta^a, -, \Delta^c) : \sset \to \sset_{(\Delta^a \sqcup \Delta^c)/}
		\]
		commutes with colimits. 
	\end{lem*}
	\begin{myproof}
		This is a consequence of the following three facts. First, the functor
		\[
			(-) \times \Delta^c : \sset \to \sset
		\]
		commutes with colimits. Second,~\cite[Rmk.\ 1.2.8.2]{htt} says that the functor 
		\[
			\Delta^a \star (-) : \sset \to \sset_{\Delta^a/}
		\]
		commutes with colimits. Third, the constant functor 
		\[
			\sset \to \sset_{\Delta^c/}
		\]	
		with value $(\Delta^c \xra{\Id} \Delta^c)$	commutes with colimits. 
	\end{myproof}
	
	\begin{lem} \label{surj3} 
		Let $X_\bullet \to Y$ be an ind-proper hypercover of ind-finite type. For any integers $a, c \ge 0$, the augmented simplicial diagram 
		\[
			\ul{\Hom}\big(\mc{M}(\Delta^a \times (-) \times \Delta^c), X_\bullet\big) : \mb{\Delta}_{+}^{\op} \to \ms{IndSch}
		\]
		is an ind-proper hypercover of $X_a \underset{Y}{\times} X_c$ of ind-finite type. 
	\end{lem}
	\begin{myproof}
		We have to check~\ref{hyper0}(H2). By~\ref{move}(\ref{swap}), it suffices to check the following statement about simplicial sets: 
		\begin{itemize}
			\item Let $K \in \sset$ be an acyclic Kan complex, and fix a map $\varphi : \Delta^a \sqcup \Delta^c \to K$. Then the simplicial set 
			\[
				\mc{Z} : \Delta^b \mapsto \Hom_{\sset}\big(\mc{M}(\Delta^a, \Delta^b, \Delta^c), K\big) \underset{\Hom_{\sset}(\Delta^a \sqcup \Delta^c, K)}{\times} \{\varphi\}
			\]
			is also an acyclic Kan complex. 
		\end{itemize}
		The rest of the proof is devoted to checking this statement. 
		
		For any $L \in \sset$, we have 
		\e{
			\Hom_{\sset}(L, \mc{Z}) &\simeq \lim_{\Delta^{b} \in (\mb{\Delta}_{/L})^{\op}} \Hom_{\sset}\big( \mc{M}(\Delta^a, \Delta^b, \Delta^c), K\big) \underset{\Hom_{\sset}(\Delta^a \sqcup \Delta^c, K)}{\times} \{\varphi\} \\
			&\simeq \lim_{\Delta^{b} \in (\mb{\Delta}_{/L})^{\op}} \Hom_{\sset_{(\Delta^a \sqcup \Delta^c)/}}\big(\mc{M}(\Delta^a, \Delta^b, \Delta^c), K \big) \\
			&\simeq \Hom_{\sset_{(\Delta^a \sqcup \Delta^c)/}}\big( \colim_{\Delta^{b} \in \mb{\Delta}_{/L}} \mc{M}(\Delta^a, \Delta^b, \Delta^c), K\big) \\
			&\simeq \Hom_{\sset_{(\Delta^a \sqcup \Delta^c)/}}\big( \mc{M}(\Delta^a, L, \Delta^c), K\big)
		}
		In the second line, we have interpreted $K$ as an object of $\sset_{(\Delta^a \sqcup \Delta^c)/}$ using the map $\varphi$. The fourth line uses Lemma~\ref{surj2}. 
		
		Thus, the simplicial set $\mc{Z}$ has the right lifting property with respect to a map $L \to L'$ if and only if the following lifting property holds: 
		\begin{cd}
			\mc{M}(\Delta^a, L, \Delta^c) \ar[r] \ar[d, "f"] & K \\
			\mc{M}(\Delta^a, L', \Delta^c) \ar[ru, dashed] 
		\end{cd}
		This diagram \emph{a priori} takes place in $\sset_{(\Delta^a \sqcup \Delta^c)/}$, but the lifting property is unchanged if we interpret it as a diagram in $\sset$. 
		
		To show that $\mc{Z}$ is an acyclic Kan complex, we will show that this lifting property holds whenever $L \to L'$ is a monomorphism of simplicial sets. Since $K$ is an acyclic Kan complex by hypothesis, it suffices to show that the left vertical map $f$ is a monomorphism. This follows from the definition of $\mc{M}$. Indeed, for any simplex $\Delta^t$, there is a diagram of sets 
		\begin{cd}
			\Hom_{\sset}(\Delta^t, \Delta^c) \ar[r, leftarrow] \ar{d}[rotate=90, anchor=north]{\sim} & \Hom_{\sset}(\Delta^t, L \times \Delta^c) \ar[r, hookrightarrow] \ar[d, hookrightarrow] & \Hom_{\sset}(\Delta^t, \Delta^a \star (L \times \Delta^c)) \ar[d, hookrightarrow] \\
			\Hom_{\sset}(\Delta^t, \Delta^c) \ar[r, leftarrow] & \Hom_{\sset}(\Delta^t, L' \times \Delta^c) \ar[r, hookrightarrow] & \Hom_{\sset}(\Delta^t, \Delta^a \star (L' \times \Delta^c)) 
		\end{cd}
		We have indicated arrows which are obviously injections or bijections. Furthermore, it is easy to see that the right square is cartesian. Finally, the value of $f$ at $\Delta^t$ is obtained by taking colimits along the horizontal arrows. The aforementioned properties of this diagram imply that this is an injection. 
	\end{myproof}
	
	\begin{lem}\label{surj5} 
		The embedding $\mb{\Delta} \hra \fcsimpset$ is initial. 
	\end{lem}
	\begin{myproof}
		By~\cite[Thm.\ 4.1.3.1]{htt}, it suffices to show that, for every $K \in \fcsimpset$, the category $\mb{\Delta}_{/K}$ is contractible. Proposition~\ref{fiber1} gives a homotopy equivalence $\mb{\Delta}_{/K} \to K$, so this follows from the assumption that $K$ is contractible. 
	\end{myproof}
	
	\subsubsection{Proof of Theorem~\ref{hyper1}} \label{hyper1-proof} 
	
	Lemma~\ref{hyper2} says that $\D(Y) \to \D(\colim X_\bullet)$ is fully faithful. To show essential surjectivity, we will fix an object $\mc{F}_\bullet \in \D(\colim X_\bullet)$ and show that it lies in the essential image. 
	
	First, we produce an object of $\D(Y)$ which should map to $\mc{F}_\bullet$. By~\cite[Cor.\ 5.5.3.4]{htt},
	there is a pair of adjoint functors 
	\begin{cd}
		\D(\colim X_\bullet) \ar[r, shift left = 1.5, "p_*"] \ar[r, hookleftarrow, swap, shift right = 1, "p^!"] & \D(Y) 
	\end{cd}
	The desired object is $p_*(\mc{F}_\bullet)$. For each $a \ge 0$, let $p_a : X_a \to Y$ be the projection map. We will construct an isomorphism $\mc{F}_\bullet \simeq p^!p_*(\mc{F}_\bullet)$ by constructing isomorphisms $\mc{F}_a \simeq p_a^!p_*(\mc{F}_\bullet)$ in each $\D(X_a)$ which are compatible as $\Delta^a \in \mb{\Delta}$ varies. 
	
	For each $a \ge 0$, we produce an ind-proper hypercover of $X_a$ via pullback along $X_a \to Y$. 
	\begin{cd}
		\colim X_a \underset{Y}{\times} X_\bullet \ar[r, "\dot{p}_a"] \ar[d, "q_a"] & \colim X_\bullet \ar[d, "p"] \\
		X_a \ar[r, "p_a"] & Y
	\end{cd}
	By ind-proper base change~\cite[Prop.\ 2.9.2]{indschemes}, we have $p_a^!p_*(\mc{F}_\bullet) \simeq (q_a)_* \dot{p}_a^! (\mc{F}_\bullet)$. We will show that the latter is isomorphic to $\mc{F}_a$. 
	
	Since $q_a^!$ is fully faithful (Lemma~\ref{hyper2}), we know that $\mc{F}_a \simeq (q_a)_*q_a^! \mc{F}_a$, so it suffices to show that $q_a^! \mc{F}_a \simeq \dot{p}_a^! (\mc{F}_\bullet)$. (Then applying $(q_a)_*$ will yield the result.) For each $c \ge 0$, we introduce notations for the following maps: 
	\begin{cd}
		X_a \underset{Y}{\times} X_c \ar[r, "\dot{p}_{a, c}"] \ar[d, "q_{a, c}"] & X_c \ar[d] \\
		X_a \ar[r] & Y
	\end{cd}
	Constructing the isomorphism $q_a^! \mc{F}_a \simeq \dot{p}_a^! (\mc{F}_\bullet)$ is equivalent to constructing isomorphisms $q_{a,c}^!\mc{F}_a \simeq \dot{p}_{a,c}^! \mc{F}_c$ which are compatible as $\Delta^c \in \mb{\Delta}$ varies. 
	
	In order to construct $q_{a,c}^!\mc{F}_a \simeq \dot{p}_{a,c}^! \mc{F}_c$, we first explain how to expand `for free' the collection of ind-schemes on which $\mc{F}_\bullet$ lives. By Lemma~\ref{surj5}, we have 
	\[
		\D(\colim X_\bullet) \simeq \lim_{\Delta^n \in \mb{\Delta}} \D(X_n) \simeq \lim_{K \in \fcsimpset} \D(\ul{\Hom}(K, X_\bullet)). 
	\]
	Therefore, we can think of $\mc{F}_\bullet$ as an object of the right hand side, i.e.\ a compatible choice of $\mc{F}_K \in \D(\ul{\Hom}(K, X_\bullet))$ for all $K \in \fcsimpset$. 
	
	Now, for fixed $a, c \ge 0$, we consider the functor 
	\[
		\mc{M}(\Delta^a, -, \Delta^c) : \mb{\Delta} \to \fcsimpset 
	\]
	defined in~\ref{surj2}. The following maps specify a diagram in $\fsimpset$: 
	\begin{cd}
		\displaystyle\lim_{\Delta^b \in \mb{\Delta}} \mc{M}(\Delta^a, \Delta^b, \Delta^c) \ar[rd, leftarrow] \ar[rrd, leftarrow] \ar[rdd, leftarrow]  \\
		& \Delta^a \sqcup \Delta^c \ar[r, leftarrow] \ar[d, leftarrow] & \Delta^c \\
		& \Delta^a
	\end{cd}
	(The limit is interpreted formally, i.e.\ it is notational shorthand for the limit diagram.) When the object $\Delta^a \sqcup \Delta^c$ is removed, this becomes a diagram in $\fcsimpset$. Applying $\D(\ul{\Hom}(-, X_\bullet))$, we obtain a diagram of $\infty$-categories: 
	\begin{cd}
		\displaystyle\lim_{\Delta^b \in \mb{\Delta}} \D(\ul{\Hom}(\mc{M}(\Delta^a, \Delta^b, \Delta^c), X_\bullet)) \ar[rd, hookleftarrow] \ar[rrd, leftarrow] \ar[rdd, leftarrow] \\
		& \D\big(X_a \underset{Y}{\times} X_c\big) \ar[r, leftarrow, swap, "\dot{p}_{a, c}^!"] \ar[d, leftarrow, "q_{a,c}^!"] & \D(X_c) \\
		& \D(X_a)
	\end{cd}
	Lemma~\ref{hyper2} and Lemma~\ref{surj3} imply that the `hook' arrow is fully faithful. The previous paragraph implies that $\mc{F}_\bullet$ determines an object in the limit of the following subdiagram, where $\D\big(X_a \underset{Y}{\times} X_c \big)$ has been excluded: 
	\begin{cd}
		\displaystyle\lim_{\Delta^b \in \mb{\Delta}} \D(\ul{\Hom}(\mc{M}(\Delta^a, \Delta^b, \Delta^c), X_\bullet)) \ar[rrd, leftarrow] \ar[rdd, leftarrow] \\
		&  & \D(X_c) \\
		& \D(X_a)
	\end{cd}
	In other words, $\mc{F}_\bullet$ determines a canonical isomorphism between the images of $\mc{F}_a \in \D(X_a)$ and $\mc{F}_c \in \D(X_c)$ in $\lim_b \D(\ul{\Hom}(\mc{M}(\Delta^a, \Delta^b, \Delta^c), X_\bullet))$. Since the `hook' arrow is fully faithful, this yields an isomorphism between the images of $\mc{F}_a$ and $\mc{F}_c$ in $\D\big(X_a \underset{Y}{\times} X_c \big)$, which is the desired isomorphism $q_{a,c}^!\mc{F}_a \simeq \dot{p}_{a,c}^! \mc{F}_c$. Moreover, the isomorphism is functorial with respect to $\Delta^a, \Delta^c \in \mb{\Delta}$ because $\mc{M}(-, -, -)$ is also functorial in its first and third arguments. This concludes the proof of essential surjectivity. \hfill $\square$ 
	
	\subsection{Fibrant replacement} The hypothesis \ref{hyper0}$(\textnormal{H}2)$ includes a Kan fibrancy requirement. In fact, this requirement can be removed because ind-proper maps are big enough: they permit one to carry out a version of the `small object' Kan fibrant replacement procedure on the simplicial indscheme. 
	
	\begin{prop} \label{fib0} 
		Let $X_\bullet \to Y$ be an augmented simplicial indscheme. Assume that the transition maps are ind-proper, and that the following property is satisfied: 
		\begin{itemize}
			\item[$(\textnormal{H}2')$] For each geometric point $\eta : \Spec F \to Y$, the simplicial set 
			\[
				\Delta^n \mapsto \Hom_Y(\Spec F, X_n)
			\]
			is contractible. 
		\end{itemize}
		Then the pullback functor $\D(Y) \to \D(\colim X_\bullet)$ is an equivalence. 
	\end{prop}
	
	\begin{rmks*}
		\begin{enumerate}[label=(\arabic*)] \item[ ]
			\item At a technical level, the proof of this proposition is similar to Varshavsky's proof of Theorem~\ref{var0} (see~\ref{intro5}), in the sense that both use an acyclic cofibration $L \hra \Delta^n$ to get an equivalence of $\D$-module $\infty$-categories. From a homotopical point of view, however, these results serve opposite purposes: Varshavsky directly uses the acyclic~cofibration
			\[
			\pt \overset{\sim}{\hra} (\text{Bruhat--Tits building}), 
			\]
			while this proposition creates an acyclic fibration: 
			\[
			(\text{Bruhat--Tits building})\overset{\sim}{\hra} (\text{fibrant replacement}) \overset{\sim}{\sra} \pt. 
			\]
			It is easier to generalize the latter approach to arbitrary simplicial indschemes. Indeed, if $X_\bullet$ is a simplicial indscheme for which the simplicial set defined in $(\textnormal{H}2')$ is contractible, the simplicial set can be built from horn inclusions $\Lambda^n_k \hra \Delta^n$ via compositions of pushouts and retracts, but it is not clear that this procedure can be carried out compatibly with the indscheme structure. On the other hand, the `small object' Kan fibrant replacement procedure can always be carried out compatibly with the indscheme structure, and that is the main idea of this~proposition.
			\item The reader who wonders why this proposition is not phrased as a general statement about simplicial presheaves should see~\ref{fib33} where two facts specific to our situation are used: $\D$-modules satisfy Zariski descent, and $\Spec F$ is connected. 
		\end{enumerate} 
	\end{rmks*}
	
	The rest of this subsection is devoted to proving the proposition.
	
	\subsubsectiona \label{fib1} 
	Fix integers $0\le k \le n$, and consider the following pushout in $\sset$: 
	\begin{cd}
		\Lambda^n_k \ar[r, hookrightarrow] \ar[d] & \Delta^n \ar[d] \\
		\pt \ar[r, hookrightarrow, "\imath"] & \Delta^n/\Lambda^n_k
	\end{cd}
	We will use the functors 
	\begin{cd}
		\mb{\Delta} \ar[r, hookrightarrow, "\imath_*"] & \mb{\Delta}_{/(\Delta^n / \Lambda^n_k)}  \ar[r, "\mr{oblv}"] & \mb{\Delta} 
	\end{cd}
	where $\imath_* : \mb{\Delta}_{/\pt} \to \mb{\Delta}_{/(\Delta^n / \Lambda^n_k)}$ is induced by the map $\imath$, and $\mr{oblv}$ forgets the map to $\Delta^n / \Lambda^n_k$. 
	
	\subsubsectiona \label{fib2} 
	Let $X_\bullet \to Y$ be an augmented simplicial indscheme. We think of $X_\bullet$ as a functor $\mb{\Delta}^{\op} \to \ms{PreStk}_{/Y}$, and we construct another such functor $\wt{X}_\bullet$ which is an object-wise Kan fibrant replacement. 
	
	First define a functor $\ol{X}_\bullet : \mb{\Delta}_{/(\Delta^n / \Lambda^n_k)}^{\op} \to \ms{PreStk}_{/Y}$ as follows: 
	\[
		\ol{X}_\bullet(\Delta^m \xra{p} \Delta^n / \Lambda^n_k) := \begin{cases}
		X_m & \text{ if } \Im(p) = \pt \\
		\underline{\Hom}_Y(\Lambda^n_k, X_\bullet) & \text{ otherwise}. 
		\end{cases}
	\]
	The behavior on morphisms is defined using the following functoriality: a map of simplices $\Delta^m \to \Lambda^n_k$ induces a map of indschemes $\underline{\Hom}_Y(\Lambda^n_k, X_\bullet) \to X_m$. 
	
	Next, we define $\wt{X}_\bullet := \on{LKE}_{\mr{oblv}^{\mr{op}}} \ol{X}_\bullet$. 
	
	To see that $\wt{X}_\bullet$ has the claimed property, let us first describe it more concretely. Since $\mr{oblv}$ from~\ref{fib1} is a cartesian fibration which is fibered in sets, the functor $\mr{oblv}^{\op}$ is a cocartesian fibration which is fibered in sets. Therefore, the preceding left Kan extension is computed as the coproduct (in~$\ms{PreStk}_{/Y}$) of values on the fiber: 
	\[
		\wt{X}_\bullet (\Delta^m) \simeq \coprod_{p \in \Hom_{\mb{\Delta}}(\Delta^m, \Delta^n / \Lambda^n_k)}^{\ms{PreStk}} \ol{X}_\bullet(\Delta^m \xra{p} \Delta^n / \Lambda^n_k). 
	\]
	All but one of the factors of this coproduct look like $\ul{\Hom}_Y(\Lambda^n_k, X_\bullet)$. The remaining one, which is indexed by the constant map $p : \Delta^m \to \pt$, looks like $X_m$. 
	
	There is an evident map 
	\[
		\psi : X_\bullet \to \wt{X}_\bullet
	\]
	coming from the first case of the definition of $\ol{X}_\bullet$. The value of $\psi$ on any simplex $\Delta^m \in \mb{\Delta}$ is the insertion of the factor of the coproduct indexed by the constant map $p$. 
	
	For any $S \in \ms{Sch}^{\mr{aff}}_{/Y}$, the concrete description of $\wt{X}_\bullet$ implies the following isomorphism of simplicial sets: 
	\begin{align} \tag{$\diamond$} \label{small} 
	\begin{split} 
		&\Hom_Y(S, \wt{X}_\bullet) \\
		&\qquad \simeq \colim
		\left( 
		\begin{tikzcd} [ampersand replacement=\&]
		\underset{\Hom_{\sset}(\Lambda^n_k, \Hom_{Y}(S, X_\bullet))}{\displaystyle\coprod} \Lambda^n_k \ar[r] \ar[d] \& \Hom_Y(S, X_\bullet) \\
		\underset{\Hom_{\sset}(\Lambda^n_k, \Hom_{Y}(S, X_\bullet))}{\displaystyle\coprod} \Delta^n
		\end{tikzcd} 
		\right) 
	\end{split} 
	\end{align}
	To obtain this, we have used that colimits in an $\infty$-category of prestacks are evaluated object-wise. This is the desired property of $\wt{X}_\bullet$. The map $\psi$ corresponds to the insertion of $\Hom_Y(S, X_\bullet)$ into this pushout. 
	
	\subsubsectiona \label{fib22} 
	We will deduce from~\ref{fib2}(\ref{small}) that the map $\colim X_\bullet \to \colim \wt{X}_\bullet$ is an equivalence. For this, it suffices to show that 
	\[
		\colim \Hom_Y(S, X_\bullet) \simeq \colim \Hom_Y(S, \wt{X}_\bullet) 
	\]
	for all $S \in \ms{Sch}^{\mr{aff}}_{/Y}$, since the colimits are taken in the category $\ms{PreStk}_{/Y}$. This follows from~\ref{fib2}(\ref{small}) because the simplicial set described by the pushout is homotopy equivalent to $\Hom_Y(S, X_\bullet)$ since each horn inclusion $\Lambda^n_k \hra \Delta^n$ is an acyclic cofibration.

	\subsubsectiona \label{fib33} 
	Unfortunately, $\wt{X}_\bullet$ does not land in $\ms{IndSch}_{/Y}$ because the coproduct in~\ref{fib2} was taken in $\ms{PreStk}_{/Y}$. To fix this, we define a functor $\wh{X}_\bullet : \mb{\Delta}^{\op} \to \ms{IndSch}_{/Y}$ as follows: 
	\[
		{\wh{X}}_\bullet (\Delta^m) \simeq \coprod_{p \in \Hom_{\mb{\Delta}}(\Delta^m, \Delta^n / \Lambda^n_k)}^{\ms{IndSch}} \ol{X}_\bullet(\Delta^m \xra{p} \Delta^n / \Lambda^n_k), 
	\]
	where this coproduct takes place in $\ms{IndSch}_{/Y}$. This nontrivially alters the simplicial colimit in prestacks that we are interested in:  
	\[
		\colim \wt{X}_\bullet \not\simeq \colim {\wh{X}}_\bullet. 
	\]
	However, we still have an equivalence of $\D$-module $\infty$-categories 
	\[
		\D\big(\colim \wt{X}_\bullet\big) \simeq \D\big(\colim {\wh{X}}_\bullet\big)
	\]
	because $\D$-modules satisfy Zariski descent. Indeed, the map 
	\[
		\coprod^{\ms{PreStk}}_{p \in \Hom_{\mb{\Delta}}(\Delta^m, \Delta^n / \Lambda^n_k)} \ol{X}_\bullet(\Delta^m \xra{p} \Delta^n / \Lambda^n_k) \to \coprod^{\ms{IndSch}}_{p \in \Hom_{\mb{\Delta}}(\Delta^m, \Delta^n / \Lambda^n_k)} \ol{X}_\bullet(\Delta^m \xra{p} \Delta^n / \Lambda^n_k)
	\]
	is a sieve for the Zariski topology. Moreover, for any geometric point $\eta : \Spec F \to Y$, this map also becomes an equivalence upon applying $\Hom_{Y} (\Spec F, -)$, because $\Spec F$ is connected. As a consequence, although~\ref{fib2}(\ref{small}) does not generally hold when $\wt{X}_\bullet$ is replaced by $\wh{X}_\bullet$, it does hold when $S = \Spec F$. 
	
	\subsubsection{Proof of Proposition~\ref{fib0}} \label{fib4} 
	Starting from $X_\bullet \to Y$, we define a countable sequence 
	\[
		X_\bullet = X_\bullet^{(0)} \xra{\psi^{(1)}} X_\bullet^{(1)}\xra{\psi^{(2)}} \cdots 
	\]
	of simplicial indschemes over $Y$ as follows: for each $i \ge 1$, we choose a pair of integers $0 < k < n$ specifying a horn inclusion, and then we define $X_\bullet^{(i)} := \wh{X}_\bullet^{(i-1)}$ and $\psi^{(i)} := \psi$ as in~\ref{fib2}. Upon taking the colimit of this sequence, we obtain a map 
	\[
		X_\bullet \xra{\psi^{(\infty)}} X_\bullet^{(\infty)} 
	\]
	of simplicial indschemes over $Y$. Each value $X_\bullet^{(\infty)}(\Delta^n)$ is an indscheme because the maps $\psi^{(i)}$ are closed embeddings. Now~\ref{fib22} and \ref{fib33} imply that 
	\[
		\D(\colim X_\bullet) \simeq \D(\colim X_\bullet^{(\infty)}). 
	\]
	
	Let $\eta : \Spec F \to Y$ be any geometric point. Equation~\ref{fib2}(\ref{small}) says that, if the integers $0 \le k \le n$ in the previous paragraph were chosen appropriately, then $\Hom_Y(\Spec F, X_\bullet^{(\infty)})$ is the simplicial set obtained by using the small object argument to construct a fibrant replacement of $\Hom_Y(\Spec F, X_\bullet)$. Therefore $\Hom_Y(\Spec F, X_\bullet^{(\infty)})$ is an acyclic Kan complex, so hyperdescent (Theorem~\ref{hyper1}) says that the pullback functor 
	\[
		\D(Y) \to \D(\colim X_\bullet^{(\infty)})
	\]
	is an equivalence. The desired equivalence is the composition of this equivalence with the previous one. \hfill $\square$ 

	\subsection{Simplicial replacement} \label{ss-simp-repl}
	We now show that the index category $\mb{\Delta}^{\op}$ can be replaced by any small $\infty$-category. See~\ref{fiber1} for more discussion of this idea. 
	
	\begin{thm} \label{simp0} 
		Let $\mc{C}$ be a small $\infty$-category, and let $\mc{X} : \mc{C}^{\op} \to \ms{PreStk}_k$ be a diagram of ind-schemes, and let $Y$ be another ind-scheme of ind-finite type. Suppose we are given a map $\colim \mc{X} \to Y$ such that $\mc{X}(c) \to Y$ is ind-proper for each $c \in \mc{C}$. Furthermore, suppose that the following holds: 
		\begin{itemize}
			\item[$(\textnormal{H}2'')$] For each geometric point $\eta : \Spec F \to Y$, the space 
			\[
				\Hom_Y(\Spec F, \colim \mc{X})
			\]
			is contractible. 
		\end{itemize}
		Then the pullback functor $\D(Y) \to \D(\colim \mc{X})$ is an equivalence. 
	\end{thm}
	\begin{myproof}
		From~\ref{fiber1}, we have a diagram of functors 
		\begin{cd}
			\mb{\Delta}_{/\mc{C}} \ar[r, "\ell"] \ar[d, "\mr{oblv}"] & \mc{C} \\
			\mb{\Delta}
		\end{cd}
		where $\ell$ is initial. Thus $\ell^\op$ is final, which implies that 
		\e{
			\colim_{\mc{C}^{\op}} \mc{X} &\simeq \colim_{\mb{\Delta}_{/\mc{C}}^{\op}} \mc{X} \circ \ell^{\op} \\
			&\simeq \colim_{\mb{\Delta}^{\op}} \on{LKE}_{\mr{oblv}^{\op}} (\mc{X} \circ \ell^{\op}). 
		} 
		As in~\ref{fib2}, the simplicial prestack $\mr{LKE}_{\mr{oblv}^{\op}}(\mc{X} \circ \ell^{\op})$ can be concretely described as follows: 
		\[
			\mr{LKE}_{\mr{oblv}^{\op}}(\mc{X} \circ \ell^{\op})(\Delta^n) \simeq \coprod_{g : \Delta^n \to \mc{C}^{\op}}^{\ms{PreStk}} \mc{X}(g(\{n\})), 
		\]
		where the coproduct takes place in $\ms{PreStk}_{/Y}$. The expression $g(\{n\})$ appears because $\ell$ is defined using the last vertex. 
		
		As in~\ref{fib33}, we define a simplicial indscheme $\wh{\mc{X}}$ via 
		\[
			\wh{\mc{X}}(\Delta^n) \simeq \coprod_{g : \Delta^n \to \mc{C}^{\op}}^{\ms{IndSch}} \mc{X}(g(\{n\})), 
		\]
		where this coproduct takes place in $\ms{IndSch}_{/Y}$. For the same reason as before, we have an equivalence of $\D$-module $\infty$-categories: 
		\[
			\D\big( \colim_{\mb{\Delta}^{\op}} \on{LKE}_{\mr{oblv}^{\op}} (\mc{X} \circ \ell^{\op}) \big) \simeq \D\big( \colim_{\mb{\Delta}^{\op}} \wh{\mc{X}} \big). 
		\]
		
		The argument in the first paragraph also shows that, for any geometric point $\eta : \Spec F \to Y$, we have 
		\e{
			\colim_{\mc{C}^{\op}} \Hom_Y(\Spec F, -) \circ \mc{X} &\simeq \colim_{\mb{\Delta}_{/\mc{C}}^{\op}} \Hom_Y(\Spec F, -) \circ \mc{X} \circ \ell^{\op} \\
			&\simeq \colim_{\mb{\Delta}^{\op}} \Hom_Y(\Spec F, -) \circ \mr{LKE}_{\mr{oblv}^{\op}}(\mc{X} \circ \ell^{\op}) \\
			&\simeq \colim_{\mb{\Delta}^{\op}} \Hom_Y(\Spec F, -) \circ \wh{\mc{X}}
		} 
		where the last line uses that $\Spec F$ is connected. By our hypothesis $(\textnormal{H}2'')$, the left hand side is contractible. Therefore, the right hand side is contractible, so the augmented simplicial indscheme $\wh{\mc{X}} \to Y$ satisfies $(\textnormal{H}2')$. Applying Proposition~\ref{fib0}, we conclude that 
		\[
			\mc{D}(Y) \simeq \D(\colim \wh{\mc{X}}). 
		\]
		The proof is completed by putting this together with the previous equivalences.
	\end{myproof}
	
	\section{Applications} \label{s-app} 
	
	We will prove the results which were mentioned in~\ref{intro2} and~\ref{intro5}. 
	
	Let $k$ be an algebraically closed field of characteristic zero, and let $G$ be a semisimple, simply connected algebraic group over $k$. Let $\mc{P}_{I, \mr{fin}}$ be the poset of finite type subsets of $I$. For each $J \subset I$, let $\mb{P}_J \subset \LG$ be the corresponding parahoric subgroup (with $\mb{I} := \mb{P}_{\emptyset}$). 
	
	See~\ref{main-remarks} for some ideas on how to work in a more general setting than this one. 
	
	\begin{rmk*}
	A theme in this section is that the following two difficulties cancel out: 
	\begin{enumerate}[label=(\roman*)]
		\item The $\infty$-category of $\D$-modules on an infinite type scheme is defined to be the limit or colimit of $\D$-module $\infty$-categories on finite type quotients, see especially~\cite[Prop.\ 3.6.1]{raskin}. In particular, if $X$ is an infinite type scheme, then
		\[
		\D(X) \not\simeq \lim_{S \in (\ms{Sch}^{\mr{aff,f.t.}}_{/X})^{\op}} \D(S), 
		\]
		where the superscript `f.t.' stands for `finite type.' 
		\item Our hyperdescent result (Theorem~\ref{simp0}) includes a finite type hypothesis. 
	\end{enumerate}
	The applications involve infinite type situations, but we will only apply hyperdescent after taking a finite type quotient. The (crossed out) equivalence in (i) was freely used in Section~\ref{s-hyper} because all (ind-)schemes which appeared therein were (ind-)finite type. 
	\end{rmk*}
	
	\subsection{Varshavsky's theorem on \texorpdfstring{$\LG$}{LG}-invariants} The results in this subsection are due to Varshavsky, see~\cite{varshavsky}. See~\ref{intro5} for more discussion. 
	
	\begin{thm} \label{var0}
		We have $\lim_{J \in \mc{P}_{I, \mr{fin}}} \D(\LG/\mb{P}_J) \simeq \ms{Vect}_k$. 
	\end{thm}
	The following proof is different from the one which Varshavsky gave in~\cite{varshavsky}. 
	\begin{myproof}
		The map for the equivalence arises from the diagram $(\mc{P}_{I, \mr{fin}})^{\triangleright} \to \ms{IndSch}_k$ which sends $J \in \mc{P}_{I, \mr{fin}}$ to $\LG / \mb{P}_J$ and $I \mapsto \pt$. We can apply Theorem~\ref{simp0} because partial affine flag varieties are ind-proper. Thus, we only need to show that, for each field extension $k \subset F$ where $F$ is algebraically closed, we have a homotopy equivalence of spaces
		\[
			\colim_{J \in \mc{P}_{I, \mr{fin}}} (\LG/\mb{P}_J)(F) \simeq \mr{pt}. 
		\]
		If $G_F$ is the semisimple group obtained from $G$ by base change along $k \hra F$, then the left hand side is the Bruhat--Tits building for $\mc{L}G_F$. Therefore, it is~contractible. 	
	\end{myproof}
	 
	\begin{cor} \label{var1}
		For any $\mc{C} \in \D(\LG)\mod$, the map 
		\[
			\lim_{J \in \mc{P}_{I, \mr{fin}}} \mc{C}^{\mb{P}_J} \to \mc{C}^{\LG}
		\]
		is an equivalence. 
	\end{cor}
	
	For completeness, we will sketch the proof of this corollary which was given in~\cite{varshavsky}. Note, however, that this corollary is also a direct consequence of our main theorem (Theorem~\ref{colim-thm}). 
	\begin{myproof}
		For each $J \in \mc{P}_{I, \mr{fin}}$, the functor of $\mb{P}_J$-invariants can be characterized as follows: 
		\[
			\mc{C}^{\mb{P}_J} \simeq \Hom_{\D(\LG)\mod}(\D(\LG/\mb{P}_J), \mc{C}). 
		\]
		Similarly, the functor of $\LG$-invariants can be characterized as follows: 
		\[
			\mc{C}^{\LG} \simeq \Hom_{\D(\LG)\mod}(\mr{Vect}, \mc{C}), 
		\]
		where $\mr{Vect} \simeq \D(\pt)$ is the trivial $\D(\LG)$-module. Therefore it suffices to show that the map 
		\[
			\colim_{J \in \mc{P}_{I, \mr{fin}}} \D(\LG/\mb{P}_J) \to \mr{Vect}
		\]
		is an equivalence, where the colimit is computed in $\D(\LG)\mod$. Using the 2-categorical structure of $\D(\LG)\mod$, we may pass to right adjoints, so it is equivalent to show that 
		\[
			\mr{Vect} \to \lim_{J \in \mc{P}_{I, \mr{fin}}} \D(\LG/\mb{P}_J)
		\]
		is an equivalence. Finally, this follows from Theorem~\ref{var0} and the fact that the forgetful functor $\D(\LG)\mod \to \cat$ reflects limits. 
	\end{myproof}
	
	\subsection{\texorpdfstring{$\D$}{D}-modules on the simplicial affine Springer resolution} 
	The next result did not appear in~\cite{varshavsky}, but it is inspired by Varshavsky's cell-by-cell analysis of the affine Springer map. In particular, the same cell-by-cell analysis yields some finite type quotients to which we can apply hyperdescent. 
	
	\subsubsection{The affine Springer maps} 
	We define a functor 
	\[
		\mc{N} : (\mc{P}_{I, \mr{fin}})^{\triangleright} \to \ms{IndSch}_k
	\]
	as follows. On objects, we have 
	\[
		\mc{N}(J) := \LG \overset{\mb{P}_J}{\times} \mb{P}_J,
	\]
	where the upper $\mb{P}_J$ acts on $\LG$ from the right and on $\mb{P}_J$ by conjugation. The value on the cone point is $\LG$. The behavior on morphisms is defined using the obvious embeddings $\mb{P}_J \hra \mb{P}_{J'}$ whenever $J \subseteq J'$. 
	
	The maps $\pi_J : \mc{N}(J) \to \LG$ are ind-proper. Let $\LG' \hra \LG$ be the closed embedding given by the union of their images. 
	
	\begin{thm} \label{spr0} 
		We have $\lim_{J \in \mc{P}_{I, \mr{fin}}} \mc{N}(J) \simeq \D(\LG')$. 
	\end{thm}
	
	We reformulate the definition of the affine Springer maps and then prove the theorem. 
	
	\subsubsectiona \label{spr1} 
	For any $J$, the map $\pi_J$ factors as 
	\begin{cd}
		\mc{N}(J) \ar[r, hookrightarrow, "{(p, \pi_J)}"] \ar[rd, swap, "\pi_J"] & (\LG / \mb{P}_J) \times \LG \ar[d, "\pr_2"] \\
		& \LG
	\end{cd}
	where the horizontal map is a closed embedding whose first factor is the projection
	\[
		p : \mc{N}(J) \sra LG \overset{\mb{P}_J}{\times} \pt \simeq \LG / \mb{P}_J. 
	\]
	For any $g \in \LG$, the restriction
	$
		p : \pi_J^{-1}(g) \hra \LG / \mb{P}_J
	$
	identifies $\pi_J^{-1}(g)$ with the fixed points of $g$ acting on $\LG / \mb{P}_J$. (This uses that $G$ is semisimple and simply connected.) 
	
	\subsubsection{Cells} \label{spr2} 
	Each simple affine root $\alpha \in I$ corresponds to an affine-linear function $f_\alpha : \mf{t}_{i\BR} \to \BR$, so that the (closed) fundamental chamber $C_0 \subset \mf{t}_{i\BR}$ is defined by 
	\[
		C_0 = \{x\in \mf{t}_{i\BR} \, | \, f_\alpha(x) \ge 0 \text{ for all } \alpha \in I\}. 
	\]
	Let $n > 0$ be any integer, and define $W_{I, \le n} \subset W_I$ to be the subset which indexes chambers contained in the set 
	\[
		C_{\le n} := \{x\in \mf{t}_{i\BR} \, | \, n+1 \ge f_\alpha(x) \ge -n \text{ for all } \alpha \in I\}. 
	\]
	In particular, we have $C_{\le 0} = C_{0}$, and $W_{I, \le 0} = \{1\}$. 
	
	In brief, the next lemma says that $W_{I, \le n}$ is convex, combinatorially convex, and closed. In this subsection, the word `gallery' means what it usually does in building theory, namely a sequence of chambers such that each pair of consecutive chambers share a wall (\cite[Def.\ 5.15]{building}). This contrasts with the more general notion of `gallery' introduced in~\ref{gal-def}, which is a technical device to be used in the proof of Theorem~\ref{colim-thm}. 

	\begin{lem*} The following are true: 
		\begin{enumerate}[label=(\roman*)] 
			\item The set $C_{\le n} \subset \mf{t}_{i \BR}$ is convex. 
			\item If $\mc{G}$ is a minimal gallery in the affine hyperplane arrangement which goes from one chamber in $C_{\le n}$ to another, then each chamber of $\mc{G}$ lies in $C_{\le n}$. 
			\item The set $W_{I, \le n}$ is downwards-closed with respect to the (two-sided) Bruhat partial order on $W_I$. 
		\end{enumerate}
	\end{lem*} 
	\begin{myproof}
		Point (i) is obvious. 
		
		For any affine root $(\alpha, m) \in I \times \BZ$, there is the root hyperplane 
		\[
			h_{(\alpha, m)} := \{x \in \mf{t}_{i \BR} \, | \, f_\alpha(x) = m\}. 
		\] 
		For any $(\alpha, m)$ as above, if $\mc{G}$ is a minimal gallery which goes between two chambers which lie on the same side of $h_{(\alpha, m)}$, then every chamber of $\mc{G}$ lies on that side of $h_{(\alpha, m)}$. Applying this to each boundary wall of $C_{\le n}$ proves (ii). 
		
		For (iii), we have to show that, for any $w \in W_{i, \le n}$ and $(\alpha, m)$ as above, if $h_{(\alpha, m)}$ separates the chamber of $w$ from $C_0$, then the reflection of the chamber of $w$ across $h_{(\alpha, m)}$ lies in $C_{\le n}$. If $m \ge 1$, then the chamber of $w$ lies in 
		\[
		C_{\le n} \cap \{x \in \mf{t}_{i\BR} \, | \, f_\alpha(x) \ge m\}, 
		\] 
		and the reflection of this set across $h_{(\alpha, m)}$ lies in $C_{\le n}$. Similarly, if $m \le 0$, then the chamber of $w$ lies in 
		\[
		C_{\le n} \cap \{x \in \mf{t}_{i\BR} \, | \, f_\alpha(x) \le m\}, 
		\] 
		and the reflection of this set across $h_{(\alpha, m)}$ also lies in $C_{\le n}$. 
	\end{myproof}
	
	\begin{rmk*}
		If we had instead defined $W_{I, \le n}$ using the convex set 
		\[
		\{x\in \mf{t}_{i\BR} \, | \, f_\alpha(x) \ge -n \text{ for all } \alpha \in I\}, 
		\]
		then $W_{I, \le n}$ would not have been downwards-closed with respect to the Bruhat partial order for any $n \ge 1$, so the analogue of (iii) would have failed. If we had instead defined $W_{I, \le n}$ to consist of all elements of length $\le n$, then the analogues of (i) and (ii) would have failed. 
	\end{rmk*}
	
	\subsubsectiona \label{spr21} 
	Lemma~\ref{spr2}(iii) implies that taking the union of all Bruhat cells indexed by $W_{I, \le n}$ yields a closed subvariety 
	\[
		(\LG / \mb{P}_J)_{\le n} \hra \LG / \mb{P}_J. 
	\]
	for each $J \subset I$. Let $\mc{B}$ be the Bruhat--Tits building of $\LG$. The $k$-points of the varieties $(\LG / \mb{P}_J)_{\le n}$, for various $J$, define a simplicial subcomplex $\mc{B}_{\le n} \subset \mc{B}$ which is combinatorially characterized as follows: a chamber $C$ lies in $\mc{B}_{\le n}$ if and only if the type of some (equivalently, any) minimal gallery from $C_0$ to $C$ lies in $W_{I, \le n}$. 
	
	\begin{cor*}
		If $\mc{G}$ is a minimal gallery in $\mc{B}$ joining two chambers of $\mc{B}_{\le n}$, then $\mc{G}$ is contained in $\mc{B}_{\le n}$. This implies that $\mc{B}_{\le n} \subset \mc{B}$ is convex.  
	\end{cor*}
	\begin{myproof}
		If $C_1, C_2$ are two chambers of $\mc{B}$, let $\delta(C_1, C_2) \in W_I$ denote the type of a minimal gallery joining $C_1$ to $C_2$, see~\cite[Prop.\ 4.81]{building}. 	
		
		Let $C_1, C_2$ be the start and end chambers of $\mc{G}$, and let $C_3$ be any chamber of $\mc{G}$. Let $\Sigma$ be any apartment containing $C_1, C_2$. Then $\Sigma$ contains $C_3$ because apartments are combinatorially convex (\cite[Prop.\ 4.40]{building}). Let $\rho : \mc{B} \to \Sigma$ be the canonical retraction (\cite[Def.\ 4.38]{building}) associated to $C_3$ and $\Sigma$. Since $\rho$ preserves adjacency of chambers, we have 
		\e{
			\delta(\rho(C_0), \rho(C_1)) &\le \delta(C_0, C_1) \\
			\delta(\rho(C_0), \rho(C_2)) &\le \delta(C_0, C_2)
		} 
		with respect to the Bruhat partial order. Since $\delta(C_0, C_1), \delta(C_0, C_2) \in W_{I, \le n}$ by hypothesis, Lemma~\ref{spr2}(iii) implies that $\delta(\rho(C_0), \rho(C_1)), \delta(\rho(C_0), \rho(C_2)) \in W_{I, \le n}$ as well. 		
		
		Let $C'_{\le n} \subset \Sigma$ be the subset of chambers in $\Sigma$ which is defined analogously to $C_{\le n}$, using $\rho(C_0)$ as the `fundamental chamber.' The previous paragraph implies that $\rho(C_1), \rho(C_2) \subset C'_{\le n}$. Now Lemma~\ref{spr2}(ii) implies that $\rho(C_3) = C_3$ is contained in $C'_{\le n}$, so $\delta(\rho(C_0), C_3) \in W_{I, \le n}$. A key property of the minimal retraction $\rho$ is that it preserves minimal galleries from $C_3$, see the proof of~\cite[Prop.\ 4.39(2)]{building}. Therefore, $\delta(C_0, C_3) \in W_{I, \le n}$. This shows that $C_3 \subset \mc{B}_{\le n}$, as desired. 
	\end{myproof}
	
	The proof of the corollary is inspired by the proof that $\mc{B}$ has the CAT(0) property, see~\cite[Thm.\ 11.16(2)]{building}. 
	
	\subsubsectiona \label{spr22} 
	For any integer $n \ge 0$, there is an integer $D_n$ such that any congruence subgroup $\bm{G}_d \subset \LG$ with $d > D_n$ acts trivially on $(\LG / \mb{P}_J)_{\le n}$ for each $J \subset I$. This is because each $(\LG / \mb{P}_J)_{\le n}$ is a finite type variety. We also define 
	\[
		\mc{N}(J)_{\le n} := p^{-1}((\LG / \mb{P}_J)_{\le n}). 
	\]
	\begin{lem*}
		Fix $n, d$ such that $d > D_n$. Then the closed subscheme given by the image of 
		\begin{cd}[column sep = 0.4in]
			\mc{N}(J)_{\le n} \ar[r, hookrightarrow, "{(p, \pi_{J, n})}"] & (\LG / \mb{P}_J)_{\le n} \times \LG
		\end{cd}
		is invariant under the action of $\bm{G}_d$ on $\LG$ from the right. 
	\end{lem*}
	\begin{myproof}
		This follows from the interpretation of this closed subscheme as parameterizing fixed points of $\LG$ acting on $\LG / \mb{P}_J$, see~\ref{spr1}. 
	\end{myproof}
	
	As a matter of notation, we let 
	$
		\mc{N}(J)_{\le n} / \bm{G}_d
	$
	denote the quotient of $\bm{G}_d$ acting on the closed subscheme in the lemma. We can now express the affine Springer map in terms of ind-finite type data: 
	\[
		\left( 
		\begin{tikzcd}
			\mc{N}(J) \ar[d, "\pi_J"] \\
			\LG'
		\end{tikzcd} \right) 
		\simeq 
		\colim_n \lim_{d > D_n} \left( 
		\begin{tikzcd}
			\mc{N}(J)_{\le n} / \bm{G}_d \ar[d, "\pi_{J, n, d}"] \\
			\LG' / \bm{G}_d
		\end{tikzcd}\right) 
	\]
	
	This discussion also implies that there is an ind-scheme presentation $\LG = \colim_m \LG'_{(m)}$ where each closed subscheme $\LG'_{(m)}$ is invariant under some congruence subgroup (acting from the right). Indeed, take each $\LG'_{(m)}$ to be the schematic image of $\pi_{\emptyset, n}: \mc{N}(\emptyset)_{\le n} \to \LG$ for some $n$, which is invariant under $\bm{G}_d$ for $d > D_n$.

	\subsubsection{Proof of Theorem~\ref{spr0}} \label{spr0-proof} 
		For each $m$, there exists an integer $N_m$ such that, for all $n > N_m$, the map
		\[
			\colim_{J \in \mc{P}_{I, \mr{fin}}} \pi_{J, n}^{-1}(\LG'_{(m)}) \to \LG'_{(m)}
		\]
		is surjective on geometric points. It suffices to show that, for all $m, n, d$ such that $n > N_m$, $d > D_n$, and $d$ is large enough so that $\LG'_{(m)}$ is invariant under $\bm{G}_d$, we have 
		\[
			\D \big( \pi_{J, n, d}^{-1}(\LG'_{(m)}  / \bm{G}_d) \big)  \simeq \D(\LG'_{(m)} / \bm{G}_d). 
		\]
		This statement involves ind-schemes of ind-finite type, so Theorem~\ref{simp0} applies. To conclude, we need to show that, for each geometric point $g : \Spec F \to \LG'_{(m)}$, the space 
		\[
			\colim_{J \in \mc{P}_{I, \mr{fin}}} \Hom_{\LG'_{(m)} / \bm{G}_d}\big(\Spec F, \pi_{J, n, d}^{-1}(\LG'_{(m)}  / \bm{G}_d)\big)
		\]
		is contractible. In view of the cartesian diagrams 
		\begin{cd}
			\pi_{J, n}^{-1}(\LG'_{(m)}) \ar[r] \ar[d] & \pi_{J, n, d}^{-1}(\LG'_{(m)}  / \bm{G}_d) \ar[d] \\
			\LG'_{(m)} \ar[r] & \LG'_{(m)} / \bm{G}_d
		\end{cd}
		the space in question is homotopy equivalent to 
		\[
			\colim_{J \in \mc{P}_{I, \mr{fin}}} \Hom_{\LG}( \Spec F, \mc{N}(J)_{\le n}). 
		\]
		
		Let $\mc{B}$ be the Bruhat--Tits building for $\LG_F$ (note the change of base field from $k$ to $F$). Let $\mc{B}^g \subset \mc{B}$ be the subcomplex fixed by $g$, and let $\mc{B}_{\le n} \subset \mc{B}$ be the subcomplex defined in~\ref{spr21}. The fixed points interpretation of~\ref{spr1} makes it clear that 
		\[
			\colim_{J \in \mc{P}_I^{\mr{fin}}} \Hom_{\LG}(\Spec F, \mc{N}(J)_{\le n}) \simeq \mc{B}^g \cap \mc{B}_{\le n}. 
		\]
		Since $\mc{B}$ is a Euclidean building, it is a $\text{CAT}(0)$ space, so geodesics are unique~\cite[Prop.\ 11.5]{building}. If $x, y \in \mc{B}^g$, then the unique geodesic joining $x$ to $y$ must lie in $\mc{B}^g$. (If not, then acting by $g$ gives a different geodesic, which contradicts uniqueness.) Therefore $\mc{B}^g$ is convex. Corollary~\ref{spr21} says that $\mc{B}_{\le n}$ is convex, so $\mc{B}^g \cap \mc{B}_{\le n}$ is convex as well. It is also nonempty since $n > N_m$, so it is contractible. \hfill $\square$ 
		
	\begin{rmk*}
		Note that this proof uses something a bit stronger than the mere contractibility of $\mc{B}^g$. Rather, $\mc{B}^g$ appears intersected with $\mc{B}_{\le n}$ because $\D$-modules on infinite type schemes must be defined via pullbacks from the finite type situation. This is a common issue: to apply hyperdescent in an infinite type setting, one sometimes needs to establish contractibility of all filtered pieces of a filtered topological space. (In a sense, the majority of the work in~\cite{hecke} is done to prove a `filtered contractibility' result.) As we will see in the next subsection, however, no filtered contractibility is needed for our main theorem. 
	\end{rmk*}
	
	\subsection{Colimit presentation of \texorpdfstring{$\D(\LG)$}{D(LG)}} 
	Here is the main theorem of this paper: 
	
	\begin{thm} \label{colim-thm} 
		We have 
		\[
			\colim_{J \in \mc{P}_{I, \mr{fin}}} \D(\mb{P}_J) \simeq \D(\LG), 
		\]
		where the colimit takes place in $\alg(\cat)$. 
	\end{thm}
	\begin{myproof}
		The first step is to reformulate this monoidal colimit as a non-monoidal colimit over a larger category. This is completely analogous to~\cite[5.2]{hecke}, but we write out the details for the reader's convenience and to show that we are not eliding anything difficult. 
		
		Let $\mc{C}$ be any monoidal $\infty$-category. As explained in~\cite[5.2]{hecke}, the underlying object of the colimit of a functor 
		\[
			F : \mc{P}_{I, \mr{fin}} \to \alg(\mc{C}) 
		\]
		is equivalent to the colimit of a functor 
		\[
			\wt{F} : (\mc{P}_{I, \mr{fin}})^{\amalg}_{\mr{act}} \to \mc{C}. 
		\]
		Here the category $(\mc{P}_{I, \mr{fin}})^{\amalg}_{\mr{act}}$ is defined as follows: 
		\begin{itemize}
			\item The objects are sequences $(J_1, J_2, \ldots, J_n)$ of objects in $\mc{P}_{I, \mr{fin}}$. 
			\item A morphism 
			\[
			(J_1, J_2, \ldots, J_n) \to (J_1', J_2', \ldots, J_m')
			\]
			is a weakly increasing map $f : \{1, \ldots, n\} \to \{1, \ldots, m\}$ which satisfies the following: 
			\begin{itemize}
				\item For each $j \in \{1, \ldots, m\}$, we have $J_j' \supseteq \bigcup_{i \in f^{-1}(j)} J_i$. 
			\end{itemize}
		\end{itemize}
		Tracing through the construction of~\cite[Cor.\ 5.2.7]{hecke}, one sees that the behavior of $\wt{F}$ on objects is given by 
		\[
			\wt{F}((J_1, \ldots, J_n)) \simeq F(J_1) \otimes \cdots \otimes F(J_n). 
		\]
		The behavior on morphisms is defined using the monoidal structures of the $F(J)$'s. Finally, the desired equivalence 
		\[
			\mr{oblv}(\colim F) \simeq \colim \wt{F} 
		\]
		is stated in~\cite[Prop.\ 5.2.9]{hecke}. 
		
		Next, we choose $\mc{C}$ judiciously, following~\cite[5.2.10]{hecke}. Since $\emptyset \in \mc{P}_{I, \mr{fin}}$ is initial, it is equivalent to evaluate $\colim_{J\in \mc{P}_{I, \mr{fin}}} \D(\mb{P}_J)$ in the $\infty$-category $\alg(\cat)_{\D(\mb{I})/}$. By~\cite[Cor.\ 3.4.1.7]{ha}, there is an equivalence of $\infty$-categories 
		\[
			\alg(\cat)_{\D(\mb{I})/} \simeq \alg(\prescript{}{\D(\mb{I})}{\bmod}_{\D(\mb{I})}(\cat)), 
		\]
		where the monoidal structure on the bimodule $\infty$-category is given by relative tensor product over $\D(\mb{I})$. When we take $\mc{C}$ to be this bimodule $\infty$-category, the preceding paragraph reduces us to computing $\colim \wt{F}$ where $\wt{F}$ is a functor whose behavior on objects is as follows: 
		\[
			\wt{F}((J_1, \ldots, J_n)) \simeq \D(\mb{P}_{J_1}) \underset{\D(\mb{I})}{\otimes} \cdots \underset{\D(\mb{I})}{\otimes} \D(\mb{P}_{J_n}). 
		\]
		By~\cite[Cor.\ 3.4.4.6]{ha}, it suffices to compute $\colim \wt{F}$ in the $\infty$-category $\cat$, forgetting the bimodule structure. More precisely, we want to show that the functor 
		\[
			\colim \wt{F} \to \D(\LG) 
		\]
		is an equivalence in $\cat$. 
		
		The next step is to interpret $\colim \wt{F}$ as the $\infty$-category of $\D$-modules on a simplicial scheme. A standard argument involving descent for $\D$-modules implies that 
		\[
			\wt{F}((J_1, \ldots, J_n)) \simeq \D\big(\mb{P}_{J_1} \overset{\mb{I}}{\times} \cdots \overset{\mb{I}}{\times}  \mb{P}_{J_n}\big), 
		\]
		and, in this reformulation, the behavior of $\wt{F}$ on morphisms is given by $*$-pushforward along the three kinds of morphisms displayed in~\ref{outline}. In other words, if we define the functor 
		\e{
			\mc{X} : (\mc{P}_{I, \mr{fin}})^{\amalg}_{\mr{act}} & \to \ms{Sch}_k \\
			(J_1, \ldots, J_n) &\mapsto \mb{P}_{J_1} \overset{\mb{I}}{\times} \cdots \overset{\mb{I}}{\times}  \mb{P}_{J_n}
		}
		then we have $\wt{F} \simeq \D(-) \circ \mc{X}$. Moreover, we can extend $\mc{X}$ to a functor 
		\[
			((\mc{P}_{I, \mr{fin}})^{\amalg}_{\mr{act}})^{\triangleright} \to \ms{IndSch}_k
		\]
		by sending the cone point to $\LG$. 
		
		We want to finish by applying Theorem~\ref{simp0}. We can reduce to a finite type situation by taking quotients on the right by any congruence subgroup $\bm{G}_d \subset \LG$; one recovers the original situation by taking a limit over $d$. (This quotient does not affect the subsequent arguments involving geometric fibers.) To show that $\mc{X}$ sends every morphism to a proper map, take quotients on the right by $\mb{I}$ and recall that Bott--Samelson varieties are proper. Finally, to check the condition $(\textnormal{H}2'')$, we will fix a geometric point $\eta : \Spec F \to \LG$ and show that
		\[
			\colim_{(J_1, \ldots, J_n) \in (\mc{P}_{I, \mr{fin}})^{\amalg}_{\mr{act}}} \Hom_{\LG}(\Spec F, \mc{X}((J_1, \ldots, J_n)))
		\]
		is contractible. By working with the group $G_F$ obtained by base-change of $G$ along $k \hra F$, we may assume that $F = k$. 
		
		First, note that $\Hom_{\LG}(\Spec F, \mc{X}((J_1, \ldots, J_n)))$ identifies with the closed points of the fiber of the Bott--Samelson map 
		\[
			\mb{P}_{J_1} \overset{\mb{I}}{\times} \cdots \overset{\mb{I}}{\times}  \mb{P}_{J_n} / \mb{I} \to \LG / \mb{I}
		\]
		over the point $\Spec k \xra{\eta} \LG \to \LG / \mb{I}$, which we denote $p$. It is well-known that this set is in bijection with the set of galleries 
		\[
			(C_0, F_1, C_1, F_2, \ldots, F_n, C_n)
		\]
		satisfying the following properties: 
		\begin{itemize}
			\item $C_0$ is the fundamental chamber and $C_n$ is the chamber labeled by $p$. 
			\item For each $i \ge 1$, the face $F_i \subseteq C_{i-1}$ has type $J_i$. 
		\end{itemize}
		(The standard reference is~\cite{gallery}, although this observation goes back to Contou-Carr\`ere's 1983 thesis.) With this interpretation in mind, the unstraightening of the functor 
		\[
			\Hom_{\LG}(\Spec F, \mc{X}(-)) : (\mc{P}_{I, \mr{fin}})^{\amalg}_{\mr{act}} \to \set 
		\]
		identifies with the obvious functor 
		\[
			\gal_K(\text{fundamental chamber}, \text{chamber labeled by } p) \to (\mc{P}_{I, \mr{fin}})^{\amalg}_{\mr{act}}. 
		\]
		Hence, by Thomason's theorem on homotopy colimits~\cite{thomason}, it suffices to show that this category of galleries is contractible. The latter is true because the Bruhat--Tits building is contractible, and the category of galleries is a model for its topological path space, by Theorem~\ref{stone-thm} and Corollary~\ref{gal00}. 
	\end{myproof}

	%\newpage

\end{document}